\documentclass[12pt,a4paper,draft]{article}

\usepackage{amsfonts,amssymb,amsmath,amscd}
\newcommand{\rr}{\mathbb R}
\newcommand{\zz}{\mathbb Z}
\newcommand{\nn}{\mathbb N}
\renewcommand{\ss}{\mathcal S}
\renewcommand{\aa}{\mathcal R}
\newcommand{\kk}{\mathcal K}
\newcommand{\Fcal}{\mathcal{F}}
\newcommand{\Acal}{\mathcal{A}}
\newcommand{\Tfrak}{\mathfrak T}
\newcommand{\BB}{\mathfrak B}

\newcommand{\hh}{\mathfrak H}
\newcommand{\LL}{\mathfrak L}
\newcommand{\LU}{\mathbf{LU}}
\newcommand{\LD}{\mathbf{LD}}

\renewcommand{\mod}{\mathop{\mathrm{mod}}\nolimits}

\newcommand{\Orb}{\mathop{\mathrm{Orb}}\nolimits}
\newcommand{\Cl}[1]{\mathop{\overline{#1}}\nolimits}
\renewcommand{\Bar}[1]{\mathop{\overline{#1}}\nolimits}
\newcommand{\EssP}[1]{\mathop{{\mathcal P}(#1)}\nolimits}
\newcommand{\Equiv}[3]{#1 \equiv #2 \; (\mod\,#3)}
\newcommand{\nEquiv}[3]{#1 \not\equiv #2 \; (\mod\,#3)}
\renewcommand{\gcd}[2]{\mathop{\mbox{\rm gcd}}\nolimits (#1, #2)}
\newcommand{\lcm}[2]{\mathop{\mbox{\rm lcm}}\nolimits (#1, #2)}
\newcommand{\LimInv}{\mathop{\mathrm{lim} \;\, \mathrm{inv}}\limits}
\newcommand{\Tau}{\mathcal T}
\newcommand{\dist}{\mathrm d}
\newcommand{\Int}{\mathrm{Int}}
\newcommand{\Diam}{\mathrm{Diam}}
\newcommand{\zer}{\mathop{\mathrm{zer}}\nolimits}
\newcommand{\fact}{\mathop{\mathrm{fact}}\nolimits}

\newcommand{\Ob}[1]{\mathop{\mathrm{Ob}}#1}
\newcommand{\Mor}[1]{\mathop{\mathrm{Mor}}#1}
\newcommand{\Iso}[1]{\mathop{\mathrm{Iso}}#1}

\newcounter{Stepcount}
{\setcounter{Stepcount}{0}\medskip}

\newenvironment{proof}[1][\hspace{-0.6ex}]{\noindent{\bf Proof #1.} { }}{$\square$ \medskip}

\newtheorem{defn}{Definition}[section]
\newtheorem{rem}{Remark}[section]
\newtheorem{prop}{Proposition}[section]
\newtheorem{lemma}{Lemma}[section]
\newtheorem{theorem}{Theorem}[section]
\newtheorem{corr}{Corollary}[section]
\newtheorem{example}{Example}[section]

\title{On projections onto odometers of dynamical systems
with the compact phase space.
\thanks{AMS classification 2000: 37B05, 37B20}
\thanks{Keywords: compact, dynamical system, projection, odometer}
}
\author{Eugene Polulyakh}
\date{e-mail address: \sl polulyah@imath.kiev.ua}
\begin{document}
\maketitle
\begin{abstract}
We investigate projections to odometers (group rotations over adic groups)
of topological invertible dynamical systems with discrete time and compact
Hausdorff phase space.

For a dynamical system $(X, f)$ with a compact phase space we consider
the category of its projections onto odometers. We examine the connected
partial order relation on the class of all objects of a skeleton of this
category. We claim that this partially ordered class always have maximal
elements and characterize them. It is claimed also, that this class have
a greatest element and is isomorphic to some characteristic for the dynamical
system $(X, f)$ subset of the set $\Sigma$ of ultranatural numbers
if and only if the dynamical system $(X, f)$ is indecomposable (the space
$X$ could not be decomposed into two proper disjoint closed
invariant subsets).

\end{abstract}
\section*{Introduction}

An important role in analysis of convertible dynamical systems (d. s.) with discrete
time (cascades) plays the information
\begin{itemize}
    \item[---] on what minimal dynamical systems there exists
        projection of the given dynamical system $(X, f)$;
    \item[---] how these projections are arranged;
    \item[---] how are dependent different projections of d. s. $(X, f)$ on
        minimal d. s.; in particular, whether for two projections
        $h_{1}: (X, f) \rightarrow (Y_{1}, g_{1})$ and
        $h_{2}: (X, f) \rightarrow (Y_{2}, g_{2})$ there exists a morphism
        $\psi: (Y_{1}, g_{1}) \rightarrow (Y_{2}, g_{2})$ of dynamical systems,
        such that $h_{2} = \psi \circ h_{1}$.
\end{itemize}

It is  is highly nontrivial problem to receive the answers
on these questions in general case. To approach to its solution, modern contributors
consider projections of given d. s. not onto all minimal d. s., but on
some "convenient" classes of minimal d. s. (distal and
equicontinuous minimal d. s., uniquely ergodic minimal d. s. etc.).

Suppose we consider some family $\mathfrak{A}$ of minimal dynamical systems
and explore properties of projections of d. s. $(X, f)$ onto elements
of this family. In some cases elements of a class of all
projections of d. s. $(X, f)$ onto d. s. from $\mathfrak{A}$ can be put in order
in the following sense.

Let $h_{1}: (X, f) \rightarrow (Y_{1}, g_{1}) $ and $h_{2}: (X, f)
\rightarrow (Y_{2}, g_{2})$ are projections. We shall say that $h_{1} \sim h_{2}$,
if there is an isomorphism of dynamical systems $\psi : (Y_{1}, g_{1})
\rightarrow (Y_{2}, g_{2})$, such that $h_{2} = \psi \circ h_{1}$.
Let us designate by $\mathfrak{B}$ a factor-set of all projections from $(X, f)$
onto elements of $\mathfrak{A}$ on this equivalence relation.
We introduce binary relation $\preceq$ on $\mathfrak{B}$. Let $B_{1}$,
$B_{2} \in \mathfrak{B}$. We say that $B_{1} \preceq B_{2}$ if there exist
representatives $h_{1}: (X, f) \rightarrow (Y_{1}, g_{1})$ of the class
$B_{1}$ and
$h_{2}: (X, f) \rightarrow (Y_{2}, g_{2})$ of the class $B_{2}$ and also
morphism $ \psi: (Y_{1}, g_{1}) \rightarrow (Y_{2}, g_{2})$ such that
$h_{2} = \psi \circ h_{1}$. The relation $\preceq$ is easily checked to
be defined correctly (i. e. it does not depend on the choice of representatives from $B_{1}$
and $B_{2}$).

It is important to know, whether $\preceq$ is the partial order relation on the class
$\mathfrak{B}$.%
\footnote {%
It is easy to test, that if the relation $\preceq$ is not partial order
on $\mathfrak{B}$ then for some d. s. $(Y, g)$ from
$\mathfrak{A}$ there is a morphism $\alpha : (Y, g) \rightarrow (Y, g)$
such that the map $\alpha : Y \rightarrow Y$ is not injective.
The problem on existence of such minimal dynamical systems is interesting by itself.
}
In the case of the positive answer to this question there appears a problem to describe
properties of the class $\mathfrak{B}$ with the partial order $\preceq$, in particular
to determine classes of all its maximal and minimal elements and to find
the greatest and least element (if they exist).

In what follows we consider the class $\Acal$ of all odometers
(group rotations over adic groups). This class is known to
coincide with the class of all minimal distal dynamical systems
with phase space homeomorphic to the Cantor set or some finite set.
It is known also, that elements of the class $\Acal$ are classified
(up to topological conjugacy) by means of the lattice of so-called
ultranatural numbers $(\Sigma, \leq)$.

We consider dynamical systems $(X, f)$ with Hausdorff compact
phase space and their projections to elements of $\Acal$ (mark
that always there exists a trivial projection $(X, f) \rightarrow (\{ pt \},
Id)$ on dynamical system, which phase space consists from
one point).

It appears that an existence of nontrivial projections of d. s. $(X, f)$ on
elements of $\Acal$ is interconnected with an existence of so-called
periodic partitions of the d. s. $(X, f)$ (finite closed partitions
of $X$ which elements are cyclically rearranged under the action
of $f : X \rightarrow X $).

Let $ \EssP{X, f} \subseteq \nn $ is a set of cardinalities of all
periodic partitions of $(X, f)$. Then $\EssP{X, f}$ is the topological
invariant of $(X, f)$ and the existence of nontrivial projections of $ (X,
f)$ on elements of $\Acal$ is equivalent to the inequality $\EssP{X, f}
\neq \{ 1 \}$.

Designate by $\Acal (X, f)$ a class of all elements of $\Acal$ on which
we can project dynamical system $(X, f)$. Let $\Sigma (X, f)$
be a subset of the set of ultranatural numbers corresponding to
$\Acal (X, f)$. Let $\BB (X, f)$ be a set of all projections of
$(X, f)$ on elements of $\Acal$ and $\BB' (X, f)$ be a factor-class of
$\BB (X, f)$ under the relation $\sim$ (see above).

Among the main results obtained in the paper we can rank
the following statements.

Let $(X, f)$ be a dynamical system with compact Hausdorff phase
space and $\EssP{X, f} \neq \{ 1 \}$. Then
\begin{itemize}
    \item[---] binary relation $\preceq$ on $\BB' (X, f)$
        is the relation of the partial order;
    \item[---] there exists a surjective map $\Lambda_{0} :
        (\BB' (X, f), \preceq) \rightarrow (\Sigma (X, f), \leq)$ which
        preserves order relation and such that a class of all maximal
        elements from $(\BB' (X, f), \preceq)$ coincides with a pre-image
        of the greatest element of $(\Sigma (X, f), \leq)$;
    \item[---] the ordered class $(\BB' (X, f), \preceq)$ is isomorphic to
        the ordered set $(\Sigma (X, f), \leq) $ if and only
        if d. s. $(X, f)$ is indecomposable (that is the space $X$
        can not be presented as a union of two disjoint proper
        closed invariant subsets).
\end{itemize}

To receive these results we explore in detail properties of
periodic partitions, odometers and ultranatural numbers.

In the last section we extract corollaries from the main results
which relates to so-called almost one-to-one expansions of
odometers, class of dynamical systems, which is intensively explored in
last time (see~\cite{Williams, G-J, Down_2, Down_3, Down_1, Down_4, Skau_1,
Skau_2, Skau_3, Block}).

Finally, author would like to thank
Igor~Vlasenko,
Sergey~Kolyada,
Vladimir~Lubashenko,
Sergey~Maksimenko,
Mark~Pankov,
Aleksander~Prishlyak,
Vladimir~Sergeychuk,
Vladimir~Sharko and
Aleksander~Sharkovskiy
for discussion of results on seminars and series of valuable notes.
The separate gratitude I want to express to Sergey~Kolyada who
has acquainted me with modern results on expansions of
odometers.

\section*{Preliminaries}

\subsection*{Quotient spaces and factor-maps.}

Let $A$ be a certain set.

\begin{defn} \label{defn_1}
{\em Partition} of set $A$ is a family $\{ A_{\alpha}
\}_{\alpha \in \Lambda}$ of nonempty subsets of $A$
which complies with the following requirements:
\begin{itemize}
    \item[1)] $A = \bigcup_{\alpha \in \Lambda} A_{\alpha}$;
    \item[2)] $A_{\alpha} \cap A_{\beta} = \emptyset$ for all $\alpha$, $\beta
        \in \Lambda$, $\alpha \neq \beta$.
\end{itemize}
\end{defn}

\begin{defn}
Partition $\{ \widetilde {A}_{\gamma} \}_{\gamma \in \Sigma}$ of $A$
is called the {\em refinement} of partition $\{ A_{\alpha} \}_{\alpha \in
\Lambda}$ if for every $\gamma \in \Sigma$ there exists $\alpha \in
\Lambda$ such that $\widetilde {A}_{\gamma} \subseteq A_{\alpha}$.
\end{defn}

\begin{rem}
Let partition $\{ \widetilde {A}_{\gamma} \}_{\gamma \in \Sigma}$ is
the refinement of partition $\{ A_{\alpha} \}_{\alpha \in \Lambda}$ of
$A$. From property 2) of definition~\ref{defn_1} it easily follows that for any
$\alpha \in \Lambda$ and $\gamma \in \Sigma$ either $\widetilde {A}_{\gamma}
\subseteq A_{\alpha}$ or $\widetilde {A}_{\gamma} \cap A_{\alpha} =
\emptyset$.
\end{rem}

\begin{rem} \label{rem_2}
There exists a bijective correspondence between partitions of the set $A$ and
equivalence relations on $A$.
\begin{itemize}
    \item[1)] With any partition $\{ A_{\alpha} \}_{\alpha \in \Lambda}$
        we can associate an equivalence relation $\rho$
        with the help of relation
        $$
            \left( a_{1} \; \rho \; a_{2} \right) \Longleftrightarrow
            \left( \exists \, \alpha \in \Lambda \;: a_{1}, a_{2}
            \in A_{\alpha} \right) \,;
        $$
    \item[2)] conversely, a partition on equivalence classes corresponds to
        any equivalence relation $\sigma$ on $A$.
\end{itemize}
\end{rem}

Let $A$ is a set, $\mathfrak{A} = \{ A_{\alpha} \}_{\alpha \in
\Lambda}$ is partition of $A$.

\begin{defn}
Set $A/\mathfrak{A}$ which elements are the elements of
partition $\mathfrak{A}$ is called the {\em factor set} of $A$
on the partition $\mathfrak{A}$.

The map $pr : A \rightarrow \mathfrak{A}$ which associate to every
element $a \in A$ an element $A_{\alpha} \in A/\mathfrak{A}$ such that
$a \in A_{\alpha}$ is called {\em projection}.
\end{defn}

By analogy, it is possible to define a factor set $A/\rho$ under the
equivalence relation $\rho$ (see remark~\ref{rem_2}).

Let $X$ is a topological space, $\mathfrak{H} =
\{ H_{\alpha} \}_{\alpha \in \Lambda}$ is a partition on $X$.

Define topology on the set $X/\mathfrak{H}$ by the following rule:
say that subset $B \subseteq X/\mathfrak{H}$ is open if and
only if its pre-image $pr^{-1} (B)$ is open in $X$. This
topology is named {\em quotient topology} and it is the weakest
topology on $X$ in which the map $pr : X
\rightarrow X/\mathfrak{H}$ is continuous.

Let $X$ and $Y$ are topological spaces, $\mathfrak{H}$ is a
partition on $X$ and $\mathfrak{T}$ is a partition on $Y$.
Let $f : X \rightarrow Y$ be a continuous map, which translates
elements of the partition $\mathfrak{H}$ into elements of the partition $\mathfrak{T}$.
Then it is defined a continuous {\em factor-map} $\fact f :
X/\mathfrak{H} \rightarrow Y/\mathfrak{T}$ such that
the following diagram is commutative
$$
\begin{CD}
X @>{f}>> Y \\
@V{pr_{X}}VV @VV{pr_{Y}}V \\
X/\mathfrak{H} @>>{\fact f}> Y/\mathfrak{T}
\end{CD}
$$

Let again $f : X \rightarrow Y$ is a continuous map.
Designate by $\zer f$ a partition
of $X$ which elements are pre-images of points of
$Y$ under map $f$. Let $\mathfrak{T}$ be a
partition of $Y$, which elements are
points of $Y$. It is clear that $pr_{Y} = Id : Y \rightarrow Y/\mathfrak{T} $ is
identical map.

\begin{defn}
Map $\fact f : X/\zer f \rightarrow Y$ for which
the following diagram
$$
\begin{CD}
X @>{f}>> Y \\
@V{pr_{X}}VV @| \\
X/\zer f @>>{\fact f}> Y
\end{CD}
$$
is commutative is called {\em one-to-one factor} of $f$.
\end{defn}

That the one-to-one factor is injective it is checked immediately.

\begin{defn}
A continuous map $ f : X \rightarrow Y$ is referred as {\em factorial},
If $f (X) = Y$ and one-to-one factor $\fact f: X/\zer f
\rightarrow Y$ is a homeomorphism.
\end{defn}

\begin{prop}[see.~\cite{Engelking}]\label{prop_1}
Suppose that the following requirements are fulfilled
for a continuous map $f : X \rightarrow Y$:
\begin{itemize}
    \item[(1)] $f (X) = Y$;
    \item[(2)] map $f$ is open (is closed).
\end{itemize}
Then $f$ is the factorial map.
\end{prop}

In what follows we will need
\begin{lemma} \label{lemma_1}
Let $X$, $Y_{1}$, $Y_{2}$ are topological spaces, $\varphi_{1} :
X \rightarrow Y_{1}$ and $\varphi_{2}: X \rightarrow Y_{2}$ are continuous
maps.

If the map $\varphi_{1}$ is factorial then following conditions are equivalent:
\begin{itemize}
    \item[(1)] partition $\zer \varphi_{1}$ of $X $ is
        the refinement of partition $\zer \varphi_{2}$;
    \item[(2)] there exists a continuous map $\psi: Y_{1} \rightarrow
        Y_{2}$ such that $\varphi_{2} = \psi \circ \varphi_{1}$.
\end{itemize}
\end{lemma}

\begin{proof}
{\bf 1.} {}
Let partition $\zer \varphi_{1}$ is a refinement of partition $\zer
\varphi_{2}$. Then the map $\varphi_{2}$ translates elements of the partition $\zer
\varphi_{1}$ in points of space $Y_{2}$ and the
factor-map $\pi = \fact \varphi_{2} : X/\zer \varphi_{1} \rightarrow
Y_{2}$ is well defined, for which the diagram is commutative
$$
\begin{CD}
X @>{\varphi_{2}}>> Y_{2} \\
@V{pr_{1}}VV @| \\
X/\zer \varphi_{1} @>>{\pi}> Y_{2}
\end{CD}
$$

Let $\chi = \fact \varphi_{1} : X/\zer \varphi_{1} \rightarrow Y_{1}$ be
one-to-one factor of map $\varphi_{1}$, that is
the diagram is commutative
$$
\begin{CD}
X @>{\varphi_{1}}>> Y_{1} \\
@V{pr_{1}}VV @| \\
X/\zer \varphi_{1} @>>{\chi}> Y_{1}
\end{CD}
$$
Since the map $\varphi_{1}$ is factorial then $\chi$ is homeomorphism
of $X/\zer \varphi_{1}$ onto $Y_{1}$.

Consider the continuous map
$$
\psi = \pi \circ \chi^{-1} : Y_{1} \rightarrow Y_{2} \;.
$$
We have
$$
\psi \circ \varphi_{1} = \pi \circ \chi^{-1} \circ \varphi_{1} =
\pi \circ \chi^{-1} \circ \chi \circ pr_{1} = \pi \circ pr_{1} = \varphi_{2}
\;,
$$
as it was required.

\medskip
{\bf 2.} {}
Suppose that there exists a continuous map $\psi : Y_{1} \rightarrow Y_{2}$
to comply the equality $\varphi_{2} = \psi \circ \varphi_{1}$.

We fix an element $H_{1}$ of partition $\zer \varphi_{1}$. By definition
there exists $y_{1} \in Y_{1}$ such that $H_{1} = \varphi_{1}^{-1} (y_{1})$.
Let $y_{2} = \psi (y_{1}) \in Y_{2}$. Then $H_{2} = \varphi_{2}^{-1} (y_{2}) =
( \psi \circ \varphi_{1})^{-1} (y_{2}) = \varphi_{1}^{-1} (\psi^{-1} (y_{2})) \supseteq
\varphi_{1}^{-1} (y_{1}) = H_{1}$. Again by definition $H_{2}$ is the element of
partition $\zer \varphi_{2}$.

Due to arbitrariness in a choice of element $H_{1}$ of partition $\zer \varphi_{1}$
we conclude that the partition $\zer \varphi_{1}$ is refinement of the partition
$\zer \varphi_{2}$.
\end{proof}

\subsection*{Categories and functors.}

\begin{defn}
Category $\kk$ consists of a class of {\em objects} $\Ob \kk$ and class of
{\em morphisms} $\Mor \kk$, which are linked by following
conditions:
\begin{itemize}
    \item[1)] certain set $H_{\kk} (A, B)$ of morphisms of the category $\kk$
    is associated to each ordered pair $A, B \in \Ob \kk$;

    \item[2)] each morphism of a category $\kk$ belongs to one and only
    one of sets $H_{\kk} (A, B)$;

    \item[3)] in the class $\Mor \kk$ the partial binary relation of multiplication
    is defined as follows: product $\beta \circ \alpha$ of morphisms $\alpha
    \in H_{\kk} (A, B)$ and $\beta \in H_{\kk} (C, D)$ is defined if and
    only if $B = C$ and in this case $\beta \circ \alpha
    \in H_{\kk} (A, D)$;

    the partial multiplication is associative:
    $ \gamma \circ (\beta \circ \alpha) = (\gamma \circ \beta) \circ
    \alpha $ for any $ \alpha \in H_{\kk} (A, B) $, $ \beta \in H_{\kk} (B,
    C) $ and $ \gamma \in H_{\kk} (C, D) $;

    \item[4)] for every $A \in \Ob \kk$ the set $H_{\kk} (A, A)$
    contains a {\em unit morphism} $1_{A}$ such that $\beta \circ
    1_{A} = \beta$ and $1_{A} \circ \alpha = \alpha$ for any morphisms
    $\alpha \in H_{\kk} (B, A)$ and $\beta \in H_{\kk} (A, C)$.
\end{itemize}
\end{defn}

\begin{defn}
Category $\LL$ is called {\em subcategory} of category $\kk$ if
\begin{itemize}
    \item[a)] $\Ob \LL \subseteq \Ob \kk$;

    \item[b)] $\Mor \LL \subseteq \Mor \kk$;

    \item[c)] unit morphisms of $\LL$ are unit
    morphisms of $ \kk $;

    \item[d)] composition $\beta \circ \alpha$ of morphisms $\alpha$,
    $\beta \in \Mor \LL$ coincides with composition of these morphisms in
    $\kk$.
\end{itemize}
\end{defn}

\begin{defn}
Subcategory $\LL$ of category $\kk$ is named {\em complete} subcategory
if $H_{\LL} (A, B) = H_{\kk} (A, B)$ for every $A$, $B \in \Ob \LL$.
\end{defn}

\begin{defn}
Morphism $\sigma: A \rightarrow B$ is referred as {\em monomorphism}
of category $\kk$ ($\sigma \in \mathrm{Mon} \, \kk$) if for any two
morphisms $\alpha$, $\beta: X \rightarrow A$ the equality $\sigma \circ
\alpha = \sigma \circ \beta$ implies $\alpha = \beta$.
\end{defn}

\begin{defn}
Morphism $\nu: A \rightarrow B$ is called {\em epimorphism} ($\nu \in
\mathrm{Ep} \, \kk$) if for any $\alpha$, $\beta: B \rightarrow Y$ from
equality $\alpha \circ \nu = \beta \circ \nu$ follows that $\alpha =
\beta$.
\end{defn}

\begin{defn}
Morphism $ \rho: A \rightarrow B$ is named {\em bimorphism} ($\rho \in
\mathrm{Bim} \, \kk$) if $\rho \in \mathrm{Mon} \; \kk \cap \mathrm{Ep} \;
\kk$.
\end{defn}

\begin{defn}
Morphism $\varphi : A \rightarrow B$ is called {\em isomorphism}
($ \varphi \in \Iso \kk$) if there exists such morphism
$\psi: B \rightarrow A$ that $\psi \circ \varphi = 1_{A}$ and
$\varphi \circ \psi = 1_{B}$.
\end{defn}

\begin{defn}
Two objects $A$, $B \in \Ob \kk$ are referred as {\em isomorphic} if
$H_{\kk} (A, B) \cap \Iso \kk \neq \emptyset$.
\end{defn}

\begin{defn}
Complete subcategory $\mathfrak Z$ of category $\kk$ which contains exactly
one representative from each class of isomorphic objects of category
$\kk$ is called {\em skeleton} of category $\kk$.
\end{defn}

\begin{defn}
Object $0_{r}$ of category $\kk$ is named {\em right zero} of category
$\kk$ if for every $A \in \Ob \kk$ there exists a unique morphism
$\alpha_{A} : A \rightarrow 0_{r}$.
\end{defn}

\begin{defn}
Object $0_{l}$ is called {\em left zero} of category
$\kk$ if for every $A \in \Ob \kk$ there exists a unique morphism
$\beta_{A} : 0_{l} \rightarrow A$.
\end{defn}

\begin{defn}
By a {\em (monadic) covariant functor} from a category $\kk$ to a category
$\LL$ we shall mean correspondence $F : \kk \rightarrow \LL$ which satisfies
the following requirements:
\begin{itemize}
    \item[1)] $F(A) \in \Ob \LL$ for every $A \in \Ob \kk$;

    \item[2)] $F(\alpha) \in H_{\LL} (F (A), F (B))$ for every $\alpha
    \in H_{\kk} (A, B)$;

    \item[3)] $F(1_{A}) = 1_{F (A)}$ for any unit morphism
    of category $\kk$;

    \item[4)] if $\alpha \in H_{\kk} (A, B)$, $\beta \in H_{\kk} (B, C)$
    then $F(\beta \circ \alpha) = F (\beta) \circ F (\alpha)$.
\end{itemize}
\end{defn}

\begin{defn}
Monadic covariant functor which bijectively maps
category $\kk$ on a category $\LL$ is called {\em isomorphism}
of categories.
\end{defn}

\subsection*{Dynamical systems.}

\begin{defn}
{\em Dynamical system with discrete time} is a pair $(X, f)$,
where $X$ is a topological space and $f : X \rightarrow X$ is a homeomorphism.
Space $X$ is called {\em phase space} of this dynamical
system.
\end{defn}

Consider a category $\kk$ which objects are dynamical
systems and morphisms of dynamical systems $(X, f)$ and $(Y, g)$ are
continuous maps $h : X \rightarrow Y$ of their phase spaces, for
which the diagram is commutative
\begin{equation}
\begin{CD}
X @>{f}>> X \\
@V{h}VV @VV{h}V \\
Y @>{g}>> Y
\end{CD}
\end{equation}

Further we shall designate a morphism $h$ of object $(X, f)$ in $(Y, g)$
as follows
$$
h : (X, f) \rightarrow (Y, g) \;.
$$

\begin{defn}
Morphism $h : (X, f) \rightarrow (Y, g)$ is called {\em imbedding} of
dynamical system $(X, f)$ in $(Y, g)$ if the map $h$ is injective.
In this case $(X, f)$ is named {\em subsystem} of the dynamical system
$(Y, g)$.
\end{defn}

\begin{defn}
Morphism $h : (X, f) \rightarrow (Y, g)$ is called {\em projection} if
$h (X) = Y$.

The dynamical system $(Y, g)$ is named {\em factor-system} of the
dynamical system $(X, f)$.

The dynamical system $(X, f)$ is named {\em expansion} of the
dynamical system $(Y, g)$.
\end{defn}

\begin{defn}
Dynamical systems $(X, f)$ and $(Y, g)$ are {\em topologically
conjugate} if there exists such morphism $h : (X, f) \rightarrow (Y,
g)$ that the map $h : X \rightarrow Y$ is homeomorphism of the
space $X$ on $Y$.
\end{defn}

In all further considerations we shall restrict ourselves to the complete subcategory
$\kk_{0}$ of $\kk$, which objects are {\em dynamical
systems with Hausdorff compact phase spaces}. We shall name
them dynamical systems or {\em flows}.

\begin{defn}
Let $(X, f)$ is a dynamical system (space $X$ is Hausdorff and
compact). Subset $A \subseteq X$ is called {\em invariant
set} of $(X, f)$ if $f(A) = A$.
\end{defn}

With each point $x \in X$ of phase space of the dynamical system $(X,
f)$ it is usual to associate following invariant sets:
\begin{itemize}
    \item[--] trajectory of the point $x$
        $$
            \Orb_{f} (x) = \bigcup_{n \in \zz} f^{n} (x) \,;
        $$
    \item[--] closure $\Cl{\Orb_{f} (x)}$ of the trajectory of $x$;
    \item[--] $\alpha$ and $\omega$--limit sets of the point $x$
        $$
            \alpha (x) = \bigcap_{n < 0} \Cl{\bigcup_{k \leq n} f^{k} (x)} \;,
            \quad
            \omega (x) = \bigcap_{n > 0} \Cl{\bigcup_{k \geq n} f^{k} (x)} \;.
        $$
\end{itemize}

\begin{defn}
Point $x$ is called {\em stable by Poisson in negative
(positive) direction} if $\alpha (x) = \Cl{\Orb_{f} (x)}$
( if $\omega (x) = \Cl{\Orb_{f} (x)}$).

Point $x$ is called {\em stable by Poisson}, if $\alpha (x) =
\omega (x) = \Cl{\Orb_{f} (x)}$.
\end{defn}

\begin{defn}
Point $x$ is called {\em recurrent}, if for any neighborhood $U$
of $x$ there exists such $n(U) \in \nn$ that for every $k \in \zz$
inequality is fulfilled
$$
U \cap \bigcup_{i = k}^{k + n(U) - 1} f^{i} (x) \neq \emptyset \;.
$$
\end{defn}

\begin{defn}
Point $x$ is named {\em almost-periodic}, if for any neighborhood
$U$ of $x$ there exists such $n(U) \in \nn$ that
$$
\bigcup_{k \in \zz} f^{k n (U)} (x) \subseteq U \;.
$$
\end{defn}

\begin{rem}
Last definition is in no way conventional.

There is an other nomenclature (see~\cite{Gottshalk, Furstenberg}) in
which points stable by Poisson are called recurrent, and
recurrent points are called almost-periodic.
\end{rem}

\begin{defn}
Nonempty closed invariant set $A \subseteq X$ is called {\em
minimal set} of dynamical system $(X, f)$ if $A$ does not
contain any proper closed invariant subset of this dynamical
system.
\end{defn}

It is easy to see that for an object $(X, f)$ of category $\kk_{0}$ any minimal
set $A$ is characterized by the following property: $\Cl{\Orb_{f} (x)} = A$
for every $x \in A$.

Further we will need the following theorem (see~\cite{Alekseev,
Birkgoff, Furstenberg})
\begin{theorem}[Birkhoff]
Each object $(X, f)$ of $\kk_{0}$ complies with the following statements:
\begin{itemize}
    \item[--] for every $x \in X$ the sets $ \alpha(x)$ and $\omega(x)$
        contain some minimal subsets of dynamical system
        $(X, f)$;
    \item[--] for any recurrent point $x \in X $ the set
        $\Cl{\Orb_{f} (x)}$ is minimal;
    \item[--] each point $x \in A$ of an arbitrary minimal set
        $A$ is recurrent.
\end{itemize}
\end{theorem}

\begin{defn}
Dynamical system $(X, f)$ is called {\em minimal}, if its
phase space $X$ is a minimal set.
\end{defn}

In what follows we will take an advantage from the following
\begin{lemma} \label{lemma_2}
Let $(X, f)$, $(Y_{1}, g_{1})$, $(Y_{2}, g_{2}) \in \Ob \kk_{0}$,
$\varphi_{1} : (X, f) \rightarrow (Y_{1}, g_{1})$ and $\varphi_{2}
: (X, f) \rightarrow (Y_{2}, g_{2})$ are morphisms.

Suppose the map $\varphi_{1} : X \rightarrow Y_{1}$ is surjective. Then
the following conditions are equivalent:
\begin{itemize}
    \item[(1)] partition $\zer \varphi_{1}$ of $X$ is
        refinement of the partition $\zer \varphi_{2}$;
    \item[(2)] there exists a morphism $\psi: (Y_{1}, g_{1}) \rightarrow
        (Y_{2}, g_{2})$ such that $\varphi_{2} = \psi \circ \varphi_{1}$.
\end{itemize}
\end{lemma}

\begin{proof}
{\bf 1.} {}
Suppose that partition $\zer \varphi_{1}$ of space $X $ is the refinement
of partition $\zer \varphi_{2}$.

It is known, that any continuous mapping of a compact set into Hausdorff space is closed
(see~\cite{Engelking}). It is known also that continuous surjective
closed map is factorial (see proposition~\ref{prop_1}).

Thus, the surjective map of compact sets $\varphi_{1}: X \rightarrow
Y_{1}$ is factorial and we are in the conditions of lemma~\ref{lemma_1}.

Hence, there exists a continuous map $\psi : Y_{1}
\rightarrow Y_{2}$ such that $\varphi_{2} = \psi \circ \varphi_{1}$.

Let us check commutability of the diagram
$$
\begin{CD}
Y_{1} @>{g_{1}}>> Y_{1} \\
@V{\psi}VV @VV{\psi}V \\
Y_{2} @>>{g_{2}}> Y_{2}
\end{CD}
$$
Fix $y_{1} \in Y_{1}$. Since the map $\varphi_{1}$ is surjective
on the condition of Lemma then there exists $x \in \varphi_{1}^{-1} (y_{1}) \subseteq X$.
Therefore $g_{2} \circ \psi (y_{1}) = g_{2} \circ \psi \circ \varphi_{1} (x) =
g_{2} \circ \varphi_{2} (x) = \varphi_{2} \circ f (x) = \psi \circ
\varphi_{1} \circ f (x) = \psi \circ g_{1} \circ \varphi_{1} (x) =
\psi \circ g_{1} (y_{1})$.

Due to arbitrariness in the choice of $y_{1} \in Y_{1}$ we conclude that $\psi \circ
g_{1} = g_{2} \circ \psi $ and $ \psi \in \Mor \kk_{0}$.

\medskip
{\bf 2.} {}
Suppose there exists a morphism $\psi : (Y_{1}, g_{1}) \rightarrow (Y_{2},
g_{2})$, such that $\varphi_{2} = \psi \circ \varphi_{1}$.

Then $\varphi_{2} = \psi \circ \varphi_{1} : X \rightarrow Y_{2}$ and
the further proof is an exact repetition of the second part of proof of
Lemma~\ref{lemma_1}.
\end{proof}

\section{Periodic partitions.}

\subsection{Definition of periodic partition.}

Suppose we have a compact Hausdorff space $X$ and
a homeomorphism $f: X \rightarrow X$.

\begin{defn} \label{defn_1.1}
We call a finite family $W^{(m)} = \{ W_{i}^{(m)} \}_{i=0}^{m-1}$ of subsets
of space $X$ {\em periodic partition} of dynamical system $(X,
f)$ {\em of length} $m$, if it satisfies to the following requirements:
\begin{itemize}
    \item[(i)] all $W_{i}^{(m)}$ are open-closed subsets
        of $X$;
    \item[(ii)] $W^{(m)}_{i} = f (W^{(m)}_{i-1})$, $i = 1, \ldots, m-1$
        and $W^{(m)}_{0} = f (W^{(m)}_{m-1})$;
    \item[(iii)] $W^{(m)}_{i} \cap W^{(m)}_{j} = \emptyset$ when $i \neq j$;
    \item[(iv)] $X = \bigcup_{i=0}^{m-1} W^{(m)}_{i}$.
\end{itemize}
\end{defn}

\begin{defn} \label{defn_1.2}
Set of all lengths of all possible periodic partitions of dynamical
system $(X, f)$ we call {\em a set of periods} of the dynamical system
$(X, f)$ and designate it $\EssP{X, f}$.
\end{defn}

\begin{rem} \label{rem_1.1}
For any dynamical system $(X, f)$ the set $\EssP{X, f}$ is not
empty. Really, always exists trivial periodic
partition $W^{(1)} = \{ W_{0}^{(1)} = X \}$ of the dynamical system
$(X, f)$ of length $1 \in \EssP{X, f}$.
\end{rem}

\begin{rem} \label{rem_1.2}
Let $W^{(m)} = \{ W^{(m)}_{i} \}_{i=0}^{m-1}$ is a periodic partition of
dynamical system $(X, f)$ of length $m$. From properties (ii) and (iii) of
Definition~\ref{defn_1.1} it immediately follows that for every $n \in \zz$ we have
\begin{itemize}
    \item[] $f^{n} (W^{(m)}_{0}) = W^{(m)}_{i}$ when $\Equiv{n}{i}{m}$ and
    \item[] $f^{n} (W^{(m)}_{0}) \cap W^{(m)}_{i} = \emptyset$ when
        $\nEquiv{n}{i}{m}$.
\end{itemize}
More commonly
\begin{itemize}
    \item[] $f^{n} (W^{(m)}_{i}) = W^{(m)}_{j}$ when $\Equiv{n}{j-i}{m}$ and
    \item[] $f^{n} (W^{(m)}_{i}) \cap W^{(m)}_{j} = \emptyset$ when
        $\nEquiv{n}{j-i}{m}$.
\end{itemize}
\end{rem}

The base properties of the set $\EssP{X, f}$ are described by two
following statements

\begin{prop} \label{prop_1.1}
Let $m \in \EssP{X, f}$ and $m$ is divided by $d \in \nn$. Then $d \in
\EssP{X, f}$.
\end{prop}

\begin{proof}
Let $\{ W^{(m)}_{i} \}_{i=0}^{m-1}$ be a periodic partition of length $m$.
Let us present $m$ as $m=ad$, $a \in \nn$. Consider a family of sets
$$
V^{(d)}_{j} = \bigcup_{k=0}^{a-1} W^{(m)}_{j+kd} =
\bigcup_{%
\begin{subarray}{c}
    s \in \{ 0, \ldots, m-1 \} \;, \\
    \Equiv{s}{j}{d}
\end{subarray}%
}
f^{s} (W^{(m)}_{0}) \;,
\quad j = 0, \ldots, d-1 \;.
$$

It is obvious that the family $ \{ V^{(d)}_{j} \}_{j=0}^{d-1} $ defined in this way
complies with properties~(i), (iii) and (iv) of definition~\ref{defn_1.1}.
Let us check that it satisfies to property~(ii) of this definition.

Since $\Equiv{m}{0}{d}$ then congruences $\Equiv{s}{j}{d}$ and
$\Equiv{s}{j+tm}{d}$ are equivalent for all $t \in \zz$. On the other hand,
$f^{tm} (W^{(m)}_{s}) = W^{(m)}_{s}$, $t \in \zz$, $s = 0, \ldots, m-1$.
Therefore
$$
V^{(d)}_{j} = \bigcup_{%
\begin{subarray}{c}
    s \in \{ 0, \ldots, m-1 \} \;, \\
    \Equiv{s}{j}{d}
\end{subarray}%
}
f^{s} (W^{(m)}_{0}) =
\bigcup_{t \in \zz}
\bigcup_{%
\begin{subarray}{c}
    s \in \{ 0, \ldots, m-1 \} \;, \\
    \Equiv{s}{j}{d}
\end{subarray}%
}
f^{tm+s} (W^{(m)}_{0}) =
$$
$$
= \bigcup_{t \in \zz}
\bigcup_{%
\begin{subarray}{c}
    r \in \{ tm, \ldots, (t+1) m-1 \} \;, \\
    \Equiv{r}{j}{d}
\end{subarray}%
}
f^{r} (W^{(m)}_{0}) =
\bigcup_{%
\begin{subarray}{c}
    r \in \zz \;, \\
    \Equiv{r}{j}{d}
\end{subarray}%
}
f^{r} (W^{(m)}_{0}) \;.
$$
The validity of property~(ii) of definition~\ref{defn_1.1} is the obvious
corollary of this sequence of equalities.

Proposition is proved.
\end{proof}

\begin{prop} \label{prop_1.2}
Let $m_{1}, m_{2} \in \EssP{X, f}$ and $D$ is the least common multiple
of $m_{1}$ and $m_{2}$. Then $D \in \EssP{X, f}$.
\end{prop}

To prove Proposition~\ref{prop_1.2} we need
some additional inspection which will be done in the following
subsection.

\subsection{Main properties of periodic partitions.}

It is clear that for any $m \in \EssP{X, f}$, $m > 1$, there exist more than
one periodic partition of dynamical system $(X, f)$ of length $m$.
Really, fix a partition $W^{(m)} =
\{ W^{(m)}_{i} \}_{i=0}^{m-1}$.
With the help of cyclical permutation
of indexes in the partition $W^{(m)}$
it is possible to construct periodic partition $W^{(m)} (k) =
\{ W^{(m)}_{i} (k) \}_{i=0}^{m-1}$,
$$
W^{(m)}_{j} (k) = W^{(m)}_{i} \quad \mbox {when} \; \Equiv{j}{i+k}{m} \qquad
j = 0, \ldots, m-1 \;,
$$
for arbitrary $k \in \{ 1, \ldots, m-1 \}$.

\begin{defn} \label{defn_1.3}
Let $m \in \EssP{X, f}$. Two periodic partitions of dynamical
system $(X, f)$ are called {\em equivalent} if one partition could
be obtained from the other by means of cyclical permutation of indexes.
\end{defn}

We ask a question: if $\EssP{X, f} \neq \{ 1 \}$
then under what conditions on $(X, f)$ and $m \in \EssP{X, f}$ every two
periodic partitions of length $m$ are equivalent?

\begin{defn} \label{defn_1.4}
We say that dynamical system $(X, f)$ is {\em indecomposable} if it
satisfies to the following property:
\begin{itemize}
    \item[(A)] If $X = X_{1} \cup X_{2}$, $X_{1} \cap X_{2} = \emptyset$
        and $X_{1}$, $X_{2}$ are closed invariant subsets of
        $(X, f)$, then either $X_{1} = \emptyset$ or
        $X_{2} = \emptyset$.
\end{itemize}
\end{defn}

\begin{rem} \label{rem_1.3}
Assume $K$ is a closed invariant set of $(X, f)$ and
$W^{(m)} = \{ W^{(m)}_{i} \}_{i=0}^{m-1}$ is a periodic partition of length $m$.

For every $i = 0, \ldots, m-1$ we have
$$
f (W^{(m)}_{i} \cap K) = f (W^{(m)}_{i}) \cap f (K) = f (W^{(m)}_{i}) \cap K \;,
$$
therefore, in particular $W^{(m)}_{i} \cap K \neq \emptyset$, $i = 0, \ldots,
m-1$, and if $K$ is open-closed in $X$ then the family
$\{ V^{(m)}_{i} = W^{(m)}_{i} \cap K \}_{i=0}^{m-1}$ satisfies
to properties (i) -- (iii) of Definition~\ref{defn_1.1}.
\end{rem}

\begin{prop} \label{prop_1.3}
Let $m \in \EssP{X, f}$, $m > 1$. Dynamical system $(X, f)$
is indecomposable if and only if there exists unique up to
the cyclical permutation of indexes periodic partition
$W^{(m)}$ of length $m$.
\end{prop}

\begin{proof}
{\bf 1.}
Assume that $W^{(m)}$ and $\widetilde{W}^{(m)}$ are two nonequivalent periodic
partitions of dynamical system $(X, f)$ of length $m$.

From property (iv) of Definition~\ref{defn_1.1} it follows that with the help of cyclical
permutation of indexes in one of partitions we can achieve that
$V^{(m)}_{0} = W^{(m)}_{0} \cap \widetilde{W}^{(m)}_{0} \neq \emptyset$.
Under our supposition $W^{(m)}_{0} \neq \widetilde{W}^{(m)}_{0}$. Let,
for instance, $K = W^{(m)}_{0} \setminus \widetilde{W}^{(m)}_{0} \neq
\emptyset$.

Designate $V^{(m)}_{i} = f^{i} (V^{(m)}_{0})$, $i = 1, \ldots, m-1$.
Remark, that $V^{(m)}_{i} \subset W^{(m)}_{i}$, therefore $V^{(m)}_{i} \cap
W^{(m)}_{0} = \emptyset$ when $i = 1, \ldots, m-1$. This follows from the
requirement~(iii) of definition~\ref{defn_1.1}.

On the other hand,
$$
f (V^{(m)}_{m-1}) = f^{m} (V^{(m)}_{0}) = f^{m} (W^{(m)}_{0} \cap
\widetilde{W}^{(m)}_{0}) = f^{m} (W^{(m)}_{0}) \cap
f^{m} (\widetilde{W}^{(m)}_{0}) = W^{(m)}_{0} \cap \widetilde{W}^{(m)}_{0}
= V^{(m)}_{0} \;.
$$
The third equality is valid since $f$ is the homeomorphism,
the penultimate equality follows from the requirement~(ii) of definition~\ref{defn_1.1}.

Hence,
$$
X_{1} = \bigcup_{i=0}^{m-1} V^{(m)}_{i} = \bigcup_{i=0}^{m-1}
f^{i} (V^{(m)}_{0})
$$
is the invariant subset of $(X, f)$ (then also $X_{2} =
X \setminus X_{1}$ is invariant).
In this case $X_{1} \neq \emptyset$ on the construction and $X_{2} \neq \emptyset$
since $X_{1} \cap W^{(m)}_{0} = V^{(m)}_{0}$ and $W^{(m)}_{0} \setminus
V^{(m)}_{0} = K \subset X \setminus X_{1}$.

From the requirement~(i) of Definition~\ref{defn_1.1} it follows, that the set
$V^{(m)}_{0}$ is open-closed (and then all $V^{(m)}_{i}$ are open-closed).
Therefore, the sets $X_{1}$ and $X_{2}$ are open-closed in $X$.

So, the dynamical system $(X, f)$ is not indecomposable.

The case $\widetilde{W}^{(m)}_{0} \setminus W^{(m)}_{0} \neq \emptyset$
is considered similarly.

{ \bf 2.}
Backwards, we shall assume that dynamical system $(X, f)$ is not
indecomposable. We fix a partition $X = X_{1} \cup X_{2}$ of $X$
on two proper disjoint invariant closed subsets.
Mark, that the subsets $X_{1}$ and $X_{2}$ are open in $X$ as well.

We fix periodic partition $W^{(m)} =
\{ W^{(m)}_{i} \}_{i=0}^{m-1}$ of $(X, f)$ of length $m$.

Nonempty families of sets $\{ V^{(m), \, 1}_{i} = W^{(m)}_{i} \cap
X_{1} \}_{i=0}^{m-1}$ and $\{ V^{(m), \, 2}_{i} = W^{(m)}_{i} \cap
X_{2} \}_{i=0}^{m-1}$ comply with properties (i) -- (iii) of
Definition~\ref{defn_1.1} (see Remark~\ref{rem_1.3}).

It is easy to see that
$$
X_{j} = X_{j} \cap \bigcup_{i=0}^{m-1} W^{(m)}_{i} = \bigcup_{i=0}^{m-1}
( W^{(m)}_{i} \cap X_{j}) = \bigcup_{i=0}^{m-1} V^{(m), \, j}_{i} \;, \quad j
= 1, 2 \;.
$$
Since $X_{1} \cap X_{2} = \emptyset$, then $V^{(m),\,1}_{r} \cap
V^{(m),\,2}_{s} = \emptyset$ for every $r, s \in \{ 0, \ldots, m-1 \}$.

We set
$$
\widetilde{W}^{(m)}_{i} = V^{(m), \, 1}_{i} \cup V^{(m), \, 2}_{i-1} \;, \quad
i = 1, \ldots, m-1 \;,
$$
$$
\widetilde{W}^{(m)}_{0} = V^{(m), \, 1}_{0} \cup V^{(m), \, 2}_{m-1} \;.
$$
The immediate check shows that the family
 $\{\widetilde{W}^{(m)}_{i}\}_{i=0}^{m-1}$ is the periodic
partition of $(X, f)$ of length $m$. In addition $W^{(m)}_{0}
\setminus \widetilde{W}^{(m)}_{0} = V^{(m), \, 1}_{0} \neq \emptyset$ and
$\widetilde{W}^{(m)}_{0} \setminus W^{(m)}_{0} = V^{(m), \, 2}_{m-1} \neq
\emptyset$.

Therefore, the family $\{ \widetilde{W}^{(m)}_{i} \}$ can not be obtained
from $\{ W^{(m)}_{i} \}$ by cyclical permutation of indexes.
\end{proof}

Let $m_{1}, m_{2} \in \EssP{X, f}$. Suppose $d$ and $D$ are the greatest
common divisor and respectively the least common multiple of numbers
$m_{1}$ and $m_{2}$.

We consider periodic partitions $\{ W^{(m_{1})}_{i} \}_{i=0}^{m_{1} -1}$ and
$\{ W^{(m_{2})}_{j} \}_{j=0}^{m_{2} -1}$ of $(X, f)$ of lengths $m_{1}$ and
$m_{2}$.

\begin{prop} \label{prop_1.4}
Let for some $k \in \{ 0, \ldots, m_{1} -1 \}$, $l \in \{ 0, \ldots,
m_{2} -1 \}$ intersection $W_{k}^{1} \cap W_{l}^{2}$ is not empty.

Then the family $\{ V^{(D)}_{s} = f^{s} (W^{(m_{1})}_{k} \cap
W^{(m_{2})}_{l}) \}_{s=0}^{D-1}$ complies with the requirements (i) --- (iii)
of Definition~\ref{defn_1.1}.
\end{prop}

\begin{proof}
The set $V^{(D)}_{0} = W^{(m_{1})}_{k} \cap W^{(m_{2})}_{l}$ is open-closed
in $X$ by definition~\ref{defn_1.1}. Since $f$ is the homeomorphism,
then all $V^{(D)}_{s}$ are open-closed in $X$ and the family $\{ V^{(D)}_{s} \}$
satisfies to the requirement~(i) of Definition~\ref{defn_1.1}.

Again, taking into account the bijectivity of $f$ we shall receive
\begin{equation} \label{eq_2}
f^{s} (W^{(m_{1})}_{k} \cap W^{(m_{2})}_{l}) = f^{s} (W^{(m_{1})}_{k}) \cap
f^{s} (W^{(m_{2})}_{l}) \;, \quad s \in \zz \;,
\end{equation}
In particular
$$
f (V^{(D)}_{D-1}) = f^{D} (V^{(D)}_{0}) = f^{D} (W^{(m_{1})}_{k} \cap
W^{(m_{2})}_{l}) = f^{D} (W^{(m_{1})}_{k}) \cap f^{D} (W^{(m_{2})}_{l}) =
W^{(m_{1})}_{k} \cap W^{(m_{2})}_{l} = V^{(D)}_{0} \;,
$$
since by definition $\Equiv{D}{0}{m_{r}}$, $r = 1, 2$. By this
the fulfillment of property~(ii) of Definition~\ref{defn_1.1} is proved.

Taking into account property~(ii) of Definition~\ref{defn_1.1} and equality~(\ref{eq_2}),
for the proof of property~(iii) it is enough to us now to show that
$V^{(D)}_{0} \cap V^{(D)}_{s} = \emptyset$, $s = 1, \ldots, D-1$.

Assume, that $V^{(D)}_{0} \cap f^{s} (V^{(D)}_{0}) \neq \emptyset$ for
certain $s \in \zz$. Then, in particular
$$
W^{(m_{1})}_{k} \cap f^{s} (W^{(m_{1})}_{k}) \neq \emptyset \;, \quad
W^{(m_{2})}_{l} \cap f^{s} (W^{(m_{2})}_{l}) \neq \emptyset \;,
$$
and it is possible by Remark~\ref{rem_1.2} only if
$\Equiv{s}{0}{m_{1}}$ and $\Equiv{s}{0}{m_{2}}$, that is, only when
$s$ is the common multiple of numbers $m_{1}$ and $m_{2}$. Hence $V^{(D)}_{0} \cap
V^{(D)}_{s} = \emptyset$ for $s = 1, \ldots, D-1$ and the family
$\{ V^{(D)}_{s} \}$ complies with the condition~(iii) of Definition~\ref{defn_1.1}.
\end{proof}

We designate
$$
V^{(D)}_{s} (k, l) = f^{s} (W^{(m_{1})}_{k} \cap W^{(m_{2})}_{l}) \;, \quad
s = 0, \ldots, D-1 \;,
$$
$$
A(k, l) = \bigcup_{s=0}^{D-1} V^{(D)}_{s} (k, l) \;.
$$
From Proposition~\ref{prop_1.4} follows, that $A (k, l)$ is open-closed
invariant subset of $(X, f)$, and if
$V^{(D)}_{0} (k, l) \neq \emptyset$, then the family $\{ V^{(D)}_{s} (k,
l) \}_{s=0}^{D-1}$ is periodic partition of the dynamical system
$(A (k, l), f |_{A (k, l)})$ of length $D$.

\begin{rem} \label{rem_1.4}
If dynamical system $(X, f)$ is indecomposable, then either $A (k, l) =
\emptyset$ or $A (k, l) = X$, and in this case $\{ V^{(D)}_{s} (k, l) \}_{s=0}^{D-1}$
is periodic partition of the dynamical system $(X, f)$ of length $D$ and $D
\in \EssP{X, f}$.

From property (iv) of Definition~\ref{defn_1.1} it follows that there exist $k, l$, for
which the set $A (k, l)$ is not empty. According to what has been said we conclude that
Proposition~\ref{prop_1.2} is valid for indecomposable dynamical systems.
\end{rem}

Generally speaking, if dynamical system $(X, f)$ is not indecomposable,
it is not necessary that $A (k, l) \in \{ \emptyset \} \cup \{ X \}$.

\begin{defn} \label{defn_1.5}
Let $m_{1}, m_{2} \in \EssP{X, f}$. Periodic partitions
$W^{(m_{1})}$ and $W^{(m_{2})}$ of dynamical system $(X, f)$ are called
{\em compatible}, if for any $k \in \{ 0, \ldots, m_{1} -1 \}$, $l
\in \{ 0, \ldots, m_{2} -1 \}$ either $A (k, l) = \emptyset$ or $A (k, l) =
X$.
\end{defn}

\begin{rem} \label{rem_1.5}
It easily follows from Proposition~\ref{prop_1.4}, that if $m_{2}
= m_{1}$, then compatibility of partitions $\{ W_{i}^{(m_{1})} \}_{i=0}^{m_{1} -1}$
and $\{ W_{j}^{(m_{2})} \}_{j=0}^{m_{2} -1}$ means, that these two partitions
are equivalent.
\end{rem}

\begin{rem} \label{rem_1.6}
The immediate corollary of definitions of equivalence and compatibility
of periodic partitions is the following statement.

Let periodic partitions $W^{(m_{1})}$ and $W^{(m_{2})}$ are compatible.
If partition $\widetilde{W}^{(m_{2})}$ is equivalent
to the partition $W^{(m_{2})}$, then periodic partitions $W^{(m_{1})}$ and
$\widetilde{W}^{(m_{2})}$ are compatible.
\end{rem}

If for some $m_{1}$, $m_{2} \in \EssP{X, f}$ there exist
compatible partitions of space $X$ of lengths $m_{1}$ and $m_{2}$,
then iterating argument from Remark~\ref{rem_1.4} we can claim, that
the least common multiple $D$ of numbers $m_{1}$ and $m_{2}$ belongs to $\EssP{X,
f}$. Hence, the proof of Proposition~\ref{prop_1.2} is reduced
to verification that for any pair $m_{1}$, $m_{2} \in \EssP{X, f}$
there exist relevant compatible periodic partitions of the
dynamical system $(X, f)$.

In the rest of subsection we shall prove a somewhat more common
\begin{prop} \label{prop_1.5}
Let $m_{1}$, $m_{2} \in \EssP{X, f}$. For any periodic
partition $W^{(m_{1})}$ of dynamical system $(X, f)$ of length $m_{1}$
there exists a periodic partition $W^{(m_{2})}$ of length
$m_{2}$, which is compatible with the partition $W^{(m_{1})}$.
\end{prop}

\begin{corr} \label{corr_1.1}
Let dynamical system $(X, f)$ is indecomposable. Then any two
periodic partitions of $(X, f)$ are compatible.
\end{corr}

In order to prove Proposition~\ref{prop_1.5} we shall study some
properties of the constructions given above.

From Remark~\ref{rem_1.3} it follows that the family $\{ W^{(m_{r})}_{i}
\cap A (k, l) \}_{i=0}^{m_{r}}$, $r = 1, 2 $, is periodic
partition of length $m_{r}$ for a dynamical system $(A (k, l), f |_{A (k,
l)})$.

Following statement gives the answer to a question about correlation of periodic
partitions $\{ V^{(D)}_{s} (k, l) \}_{s=0}^{D-1}$ and $\{ W^{(m_{r})}_{i} \cap
A (k, l) \}_{i=0}^{m_{r}}$, $r = 1, 2$ of the dynamical system $(A (k, l),
f |_{A (k, l)})$.

\begin{lemma} \label{lemma_1.1}
Let $V^{(D)}_{0} (k, l) \neq \emptyset$, $i \in \{ 0, \ldots, m_{1} -1 \}$,
$j \in \{ 0, \ldots, m_{2} -1 \}$.

If $\Equiv{j-i}{l-k}{d}$, then $W^{(m_{1})}_{i} \cap W^{(m_{2})}_{j}
\subseteq A (k, l)$. Moreover $W^{(m_{1})}_{i} \cap W^{(m_{2})}_{j} =
V^{(D)}_{s} (k, l)$, where $s \in \{ 0, \ldots, D-1 \}$ is a solution
of the following system of congruences
\begin{equation} \label{eq_3}
\left\{
    \begin{array}{c}
        \Equiv{s}{i-k}{m_{1}} \,, \\
        \Equiv{s}{j-l}{m_{2}} \,.
    \end{array}
\right.
\end{equation}

If $\nEquiv{j-i}{l-k}{d}$, then $(W^{(m_{1})}_{i} \cap W^{(m_{2})}_{j})
\cap A (k, l) = \emptyset$.
\end{lemma}

\begin{rem}
We remind, that the system~(\ref{eq_3}) is compatible if and only if
$\Equiv{j-i}{l-k}{d}$ (see~\cite{Vinogradov}).
\end{rem}

\begin{proof}[of Lemma~\ref{lemma_1.1}]
{\bf 1.}{}
Let $\Equiv{j-i}{l-k}{d}$. Then there exists unique
($\mathrm{mod} \; D$) solution $s$ of system~(\ref{eq_3}) (see~\cite{Vinogradov}).
Therefore,
$$
\begin{array}{rcl}
    A (k, l) \supseteq V^{(D)}_{s} (k, l) & = & f^{s} (W^{(m_{1})}_{k} \cap
        W^{(m_{2})}_{l}) = f^{s} (W^{(m_{1})}_{k}) \cap
        f^{s} (W^{(m_{2})}_{l}) = \\
    & & f^{i-k} (W^{(m_{1})}_{k}) \cap f^{j-l} (W^{(m_{2})}_{l}) = \\
    & & f^{i-k} \circ f^{k} (W^{(m_{1})}_{0}) \cap f^{j-l} \circ
        f^{l} (W^{(m_{2})}_{0}) = \\
    & & f^{i} (W^{(m_{1})}_{0}) \cap f^{j} (W^{(m_{2})}_{0})
        = W^{(m_{1})}_{i} \cap W^{(m_{2})}_{j} \,. \\
\end{array}
$$
Here we have taken advantage of  that $f$ is the homeomorphism and, in particular,
maps $f$ and $f^{-1}$ are one-to-one.

{\bf 2.}{} We shall assume now that $(W^{(m_{1})}_{i} \cap W^{(m_{2})}_{j})
\cap A (k, l) \neq \emptyset$. Then there exists $s \in \{ 0, \ldots, D-1 \}$,
such that
$$
\emptyset \neq (W^{(m_{1})}_{i} \cap W^{(m_{2})}_{j}) \cap V^{(D)}_{s} (k,
l) = (W^{(m_{1})}_{i} \cap W^{(m_{2})}_{j}) \cap (f^{s} (W^{(m_{1})}_{k})
\cap f^{s} (W^{(m_{2})}_{l})) \,.
$$
Hence, $W^{(m_{1})}_{i} \cap f^{s} (W^{(m_{1})}_{k}) \neq
\emptyset$ and $W^{(m_{2})}_{j} \cap f^{s} (W^{(m_{2})}_{l}) \neq \emptyset$.
From Remark~\ref{rem_1.2} it follows now, that $s$ satisfies
to system~(\ref{eq_3}), and in particular $\Equiv{i-k}{j-l}{d}$.
\end{proof}

\begin{corr} \label{corr_1.2}
Let for some $k_{1}, k_{2} \in \{ 0, \ldots, m_{1} -1 \}$, $l_{1},
l_{2} \in \{ 0, \ldots, m_{2} -1 \}$ sets $A (k_{1}, l_{1})$, $A (k_{2},
l_{2})$ are not empty.

If $\Equiv{k_{1} -k_{2}}{l_{1} -l_{2}}{d}$, then $A (k_{1}, l_{1}) = A (k_{2},
l_{2})$. Otherwise $A (k_{1}, l_{1}) \cap A (k_{2}, l_{2}) =
\emptyset$.
\end{corr}

\begin{proof}
We already know that $A (k_{1}, l_{1})$ and $A (k_{2}, l_{2})$ are closed
invariant subsets of dynamical system $(X, f)$.

{\bf 1.}{}
Suppose, $\Equiv{k_{1} -k_{2}}{l_{1} -l_{2}}{d}$. Then according to Lemma~\ref{lemma_1.1}
we have $(W^{(m_{1})}_{k_{1}} \cap W^{(m_{2})}_{l_{1}}) \subseteq A (k_{2}, l_{2})$. Hence,
$$
A (k_{1}, l_{1}) = \bigcup_{s=0}^{D-1} f^{s} (W^{(m_{1})}_{k_{1}} \cap
W^{(m_{2})}_{l_{1}}) \subseteq \bigcup_{s=0}^{D-1} f^{s} (A (k_{2}, l_{2}))
= A (k_{2}, l_{2}) \,.
$$
Changing roles of $A (k_{1}, l_{1})$ and $A (k_{2}, l_{2})$, we shall receive the inverse
inclusion.

{\bf 2.}{}
Suppose now that $\nEquiv{k_{1} -k_{2}}{l_{1} -l_{2}}{d}$. Then
$$
\emptyset = \bigcup_{s=0}^{D-1} f^{s} (W^{(m_{1})}_{k_{1}} \cap
W^{(m_{2})}_{l_{1}}) \cap f^{s} (A (k_{2}, l_{2})) =
$$
$$
= \bigcup_{s=0}^{D-1} f^{s} (W^{(m_{1})}_{k_{1}} \cap W^{(m_{2})}_{l_{1}})
\cap A (k_{2}, l_{2}) = A (k_{1}, l_{1}) \cap A (k_{2}, l_{2}) \,.
$$

Corollary is proved.
\end{proof}

\begin{corr} \label{corr_1.3}
If for some $k \in \{ 0, \ldots, m_{1} -1 \}$ and $l \in \{ 0, \ldots,
m_{2} -1 \}$ the equality $A (k, l) = X$ is valid, then periodic partitions
$W^{(m_{1})}$ and $W^{(m_{2})}$ are compatible.
\end{corr}

\begin{corr} \label{corr_1.4}
If for some $k \in \{ 0, \ldots, m_{1} \}$ and $l \in \{ 0, \ldots,
m_{2} \}$ there exists a set $K \subseteq W^{(m_{1})}_{k} \cap
W^{(m_{2})}_{l}$, such that
$$
\bigcup_{n \in \zz} f^{n} (K) = X \;,
$$
then periodic partitions $W^{(m_{1})}$ and $W^{(m_{2})}$ are compatible.
\end{corr}

\begin{proof}
The statement of Corollary follows from the chain of equalities
$$
\begin{array}{rcl}
    X = \bigcup_{n \in \zz} f^{n} (K) & \subseteq & \bigcup_{n \in \zz}
        f^{n} (W^{(m_{1})}_{k} \cap W^{(m_{2})}_{l}) = \bigcup_{m \in \zz}
        \bigcup_{r=0}^{D-1} f^{r+Dm} (V^{(D)}_{0} (k, l)) = \\
    & & \bigcup_{m \in \zz} f^{Dm}
        \left(
        \bigcup_{r=0}^{D-1} f^{r} (W^{(m_{1})}_{k} \cap W^{(m_{2})}_{l})
        \right) =
        \bigcup_{m \in \zz} f^{Dm} (A (k, l)) = \\
    & & A (k, l) \\
\end{array}
$$
and from Corollary~\ref{corr_1.3}
\end{proof}

\begin{corr} \label{corr_1.5}
Let $W^{(m_{1})}_{k} \supseteq W^{(m_{2})}_{l}$ for some
$k \in \{ 0, \ldots, m_{1} -1 \}$ and $l \in \{ 0, \ldots, m_{2} -1 \}$.

Then
\begin{itemize}
    \item[1)] partitions $W^{(m_{1})}$ and $W^{(m_{2})}$ are compatible;
    \item[2)] $m_{1}$ divides $m_{2}$.
\end{itemize}
\end{corr}

\begin{proof}
{\bf 1.}{}
Set $K = W^{(m_{2})}_{l}$. From property~(iv) of periodic partitions and from
Corollary~\ref{corr_1.4} it follows, that periodic partitions $ \{
W^{(m_{1})}_{i} \}$ and $\{ W^{(m_{2})}_{j} \}$ are compatible
and $A (k, l) = X$.

{\bf 2.}{}
From Lemma~\ref{lemma_1.1} follows, that $W^{(m_{1})}_{i} \cap
W^{(m_{2})}_{l} \neq \emptyset$ if and only if
$\Equiv{i}{k}{d}$.

Assume, that $d \neq m_{1}$. Then there exists $\tau \in \{ 0, \ldots,
m_{1} -1 \}$, $ \tau \neq k$, such that $\Equiv{\tau}{k}{d}$. Hence,
$W^{(m_{1})}_{\tau} \cap W^{(m_{2})}_{l} \neq \emptyset$. From the other
side, we have $W^{(m_{2})}_{l} \subseteq W^{(m_{1})}_{k}$ on condition of Corollary and
$W^{(m_{1})}_{k} \cap W^{(m_{1})}_{\tau} = \emptyset$ by
Definition~\ref{defn_1.1}.

The obtained contradiction proves, that $d = m_{1}$ and $m_{1}$ divides
$m_{2}$.
\end{proof}

\begin{corr} \label{corr_1.6}
Let $m_{1}, m_{2} \in \EssP{X, f}$ and $m_{2} $ is divided by $m_{1}$.

Periodic partitions $W^{(m_{1})}$ and $W^{(m_{2})}$ are
compatible if and only if
the partition $\{ W^{(m_{2})}_{j} \}$ of space $X$ is
refinement of the partition $\{ W^{(m_{1})}_{i} \}$.
\end{corr}

\begin{proof}
{\bf 1.}{} Necessity. Suppose periodic partitions $W^{(m_{1})}$ and
$W^{(m_{2})}$ are compatible.

We find $k \in \{ 0, \ldots, m_{1} -1 \}$ and $l \in \{ 0, \ldots, m_{2} -1
\}$, for which $W^{(m_{1})}_{k} \cap W^{(m_{2})}_{l} \neq \emptyset$.
Then $A (k, l) = X$.

Since $m_{1}$ divides $m_{2}$, then $d = m_{1}$.

We fix $j \in \{ 0, \ldots, m_{2} -1 \}$. From Lemma~\ref{lemma_1.1}
it follows, that $W^{(m_{1})}_{i} \cap W^{(m_{2})}_{j} \neq \emptyset$ if and
only if $ \Equiv{i}{j-l+k}{m_{1}} $. There exists
a unique $\tau \in \{ 0, \ldots, m_{1} -1 \}$, such that
$\Equiv{\tau}{j-l+k}{m_{1}}$. Since
$$
X = \bigcup_{i=0}^{m_{1} -1} W^{(m_{1})}_{i}
$$
and $W^{(m_{1})}_{i} \cap W^{(m_{2})}_{j} = \emptyset$ when $i \in \{ 0,
\ldots, m_{1} -1 \}$, $i \neq \tau $, then $W^{(m_{2})}_{j} \subseteq
W^{(m_{1})}_{\tau}$.

By virtue of arbitrariness in the choice of $j \in \{ 0, \ldots, m_{2} -1 \}$ we conclude, that
the partition $\{ W^{(m_{2})}_{j} \}$ of space $X$ is refinement
of the partition $\{ W^{(m_{1})}_{i} \}$.

{ \bf 2.}{} Sufficiency follows from Corollary~\ref{corr_1.5}.
\end{proof}

 %****************************************************************************

\begin{proof}[of Proposition~\ref{prop_1.5}]
Assume that $W^{(m_{1})}$ is periodic partition of length
$m_{1}$. We fix some periodic partition
$W^{(m_{2})}$ of length $m_{2}$.

Consider the following sets
$$
A (0, j) = \bigcup_{s=0}^{D-1}
f^{s} \left( W^{(m_{1})}_{0} \cap W^{(m_{2})}_{j} \right) \;,
\quad j = 0, \ldots, m_{2} -1 \;.
$$
Obviously,
$$
W^{(m_{1})}_{0} = W^{(m_{1})}_{0} \cap
\bigcup_{j=0}^{m_{2} -1} W^{(m_{2})}_{j} =
\bigcup_{j=0}^{m_{2} -1}
\left( W^{(m_{1})}_{0} \cap W^{(m_{2})}_{j} \right) \subseteq
\bigcup_{j=0}^{m_{2} -1} A (0, j) \;.
$$

Since $A (0, j)$, $j=0, \ldots, m_{2} -1$, are the invariant subsets
of $(X, f)$, then
$$
\bigcup_{j=0}^{m_{2} -1} A (0, j) = \bigcup_{j=0}^{m_{2} -1}
\bigcup_{i=0}^{m_{1} -1} f^{i} (A (0, j)) = \bigcup_{i=0}^{m_{1} -1} f^{i}
\left( \bigcup_{j=0}^{m_{2} -1} A (0, j) \right) \supseteq
\bigcup_{i=0}^{m_{1} -1} f^{i} (W^{(m_{1})}_{0}) = X \;.
$$

We know from Corollary~\ref{corr_1.2} that $A (0, j) = A (0, k)$ when
$\Equiv{j}{k}{d}$ and $A (0, j) \cap A (0, k) = \emptyset$ if
$\nEquiv{j}{k}{d}$, therefore
$$
X = \bigcup_{j=0}^{d-1} A (0, j) \;.
$$
In the last equality all sets $A (0, j)$ are pairwise disjoint.
Some of these sets can be empty. Let $A_{1} = A (0, k_{1}),
\ldots, A_{l} = A (0, k_{l})$ are all nonempty subsets from the family
$\{ A (0, j) \}_{j=0}^{d-1}$.

Designate
$$
V^{(D)}_{s} (j) = V^{(D)}_{s} (0, k_{j}) =
f^{s} \left( W^{(m_{1})}_{0} \cap W^{(m_{2})}_{k_{j}} \right)
\;, \quad j = 1, \ldots, l \;, \quad s = 0, \ldots, D-1 \;.
$$
Then
\begin{equation} \label{eq_4}
X = \bigcup_{j=1}^{l} A_{j} = \bigcup_{j=1}^{l} \bigcup_{s=0}^{D-1}
V^{(D)}_{s} (j) = \bigcup_{s=0}^{D-1}
\left( \bigcup_{j=1}^{l} V^{(D)}_{s} (j) \right) \;,
\end{equation}

Remark, that $V^{(D)}_{s} (i) \cap V^{(D)}_{r} (j) = \emptyset$ if $s
\neq r$ or $i \neq j$.

Really, $V^{(D)}_{s} (i) \subset A_{i}$, $V^{(D)}_{r} (j) \subset
A_{j}$, therefore on construction $V^{(D)}_{s} (i) \cap V^{(D)}_{r} (j) \subset
A_{i} \cap A_{j} = \emptyset$ when $i \neq j$.

Let now $i = j$. From Proposition~\ref{prop_1.4} we know, that the family
$\{ V^{(D)}_{s} (i) \}_{s=0}^{D-1}$ satisfies
to properties~(i)--(iii) of Definition~\ref{defn_1.1}. Hence,
$V^{(D)}_{s} (i) \cap V^{(D)}_{r} (i) = \emptyset$ when $s \neq r$.

We set
$$
\widetilde{W}^{(D)}_{0} = \bigcup_{j=1}^{l} V^{(D)}_{0} (j) \;, \quad
\widetilde{W}^{(D)}_{s} =
f^{s} \left( \widetilde{W}^{(D)}_{0} \right) =
\bigcup_{j=1}^{l} V^{(D)}_{s} (j) \;, \quad
s = 1, \ldots, D-1 \;.
$$

According to what has been said, we have $\widetilde{W}^{(D)}_{s} \cap
\widetilde{W}^{(D)}_{r} = \emptyset$ when $s \neq r$, that is the family
$\{\widetilde{W}^{(D)}_{s}\}_{s=0}^{D-1}$ complies with the
requirement~(iii) of Definition~\ref{defn_1.1}. From the formula~(\ref{eq_4}) it follows, that
this family satisfies to property~(iv) of the indicated definition as well.

We recollect, that all systems of sets $\{ V^{(D)}_{s} (i) \}_{s=0}^{D-1}$,
$i=1, \ldots, l$, satisfy to properties~(i)--(iii) of
Definition~\ref{defn_1.1}. From this in the first place it follows, that all sets
of family $\{ \widetilde{W}^{(D)}_{s} \}$ are open-closed in $X$, secondly
$$
f \left( \widetilde{W}^{(D)}_{D-1} \right) =
\bigcup_{j=1}^{l} f (V^{(D)}_{D-1} (j)) =
\bigcup_{j=1}^{l} V^{(D)}_{0} (j) = \widetilde{W}^{(D)}_{0} \;.
$$

So, the system of sets $\{\widetilde{W}^{(D)}_{s}\}_{s=0}^{D-1}$ is
periodic partition of $(X, f)$ of length $D$. By this we
completely have proved Proposition~\ref{prop_1.2}.

Periodic partitions $W^{(m_{1})}$ and $\widetilde{W}^{(D)}$ are compatible
since $\widetilde{W}^{(D)}_{0} \subseteq W^{(m_{1})}_{0}$ on construction
(see corollary~\ref{corr_1.5}).

We designate
$$
\widetilde{W}^{(m_{2})}_{k} = \bigcup_{%
    \begin{subarray}{c}
        s \in \{ 0, \ldots, D-1 \} \;, \cr
        \Equiv{s}{k}{m_{2}} \cr
    \end{subarray}%
} \widetilde{W}^{(D)}_{s} \;, \quad k =
0, \ldots, m_{2} -1 \;.
$$
The easy immediate check shows, that
$\{\widetilde{W}^{(m_{2})}_{k}\}_{k=0}^{m_{2}-1}$ is the periodic partition
of length $m_{2}$ (see proof of Proposition~\ref{prop_1.1}).

It follows from Corollary~\ref{corr_1.4} that periodic partitions
$W^{(m_{1})}$ and $\widetilde{W}^{(m_{2})}$ are compatible since
$W^{(m_{1})}_{0} \cap \widetilde{W}^{(m_{2})}_{0} \supseteq
\widetilde{W}^{(D)}_{0}$.

Proposition~\ref{prop_1.5} is completely proved.
\end{proof}

\begin{rem} \label{rem_1.8}
Generally speaking, if dynamical system $(X, f)$ is not
indecomposable, then starting from an arbitrary fixed partition
$W^{(m_{2})}$ it is possible to construct more than one periodic partition of length
$m_{2}$ compatible with a given partition $W^{(m_{1})}$ of length
$m_{1}$.

Namely, using notation introduced in the proof of
Proposition~\ref{prop_1.5} we assume
$$
\widetilde{W}^{(D)}_{0} (t_{1}, \ldots, t_{l}) =
\bigcup_{j=1}^{l} V^{(D)}_{t_{j}} (j) \;,
$$
where $\Equiv{t_{j}}{0}{m_{1}}$, $t_{j} \in \{ 0, \ldots, D-1 \}$, $j = 1,
\ldots, l$. Then, from Remark~\ref{rem_1.2} we conclude that
$$
\widetilde{W}^{(D)}_{0} (t_{1}, \ldots, t_{l}) \subseteq
W^{(m_{1})}_{0} \;.
$$

Designate
\begin{equation} \label{eq_5}
\widetilde{W}^{(m_{2})}_{k} (t_{1}, \ldots, t_{l}) =
\bigcup_{%
\begin{subarray}{c}
    s \in \{ 0, \ldots, D-1 \} \;, \\
    \Equiv{s}{k}{m_{2}}
\end{subarray}%
} f^{s} (\widetilde{W}^{(D)}_{0} (t_{1}, \ldots, t_{l})) \;,
\quad k = 0, \ldots, m_{2} -1 \;.
\end{equation}

Now we can show, having applied the same argument as in the proof
of Proposition~\ref{prop_1.5}, that the family of sets
$\{ \widetilde{W}^{(m_{2})}_{k} (t_{1}, \ldots, t_{l}) \}_{k=0}^{m_{2} -1}$
is periodic partition of length $m_{2}$, compatible with
the partition $W^{(m_{1})}$.

Let $ \widetilde{W}^{(m_{2})}_{k} (t_{1}, \ldots, t_{l})$ and
$\widetilde{W}^{(m_{2})}_{r} (\tau_{1}, \ldots, \tau_{l})$ are two
periodic partitions of form~(\ref{eq_5}). Substituting them on equivalent, without
breaking the compatibility relation it is possible to achieve, that $t_{1} = \tau_{1}
= 0$ (see Remark~\ref{rem_1.5}). We can easily see, that
$$
\bigcup_{s=0}^{m_{2} -1} \left(
\widetilde{W}^{(m_{2})}_{s} (0, t_{2}, \ldots, t_{l}) \cap
\widetilde{W}^{(m_{2})}_{s} (0, \tau_{2}, \ldots, \tau_{l})
\right) =
A_{1} \cup
\bigcup_{%
\begin{subarray}{c}
    t_{j} = \tau_{j} \;, \\
    j \in \{ 2, \ldots, l \}
\end{subarray}%
}
A_{j} \;.
$$
Hence, periodic partitions
$\{ \widetilde{W}^{(m_{2})}_{s} (0, t_{2}, \ldots, t_{l}) \}$ and
$\{ \widetilde{W}^{(m_{2})}_{s} (0, \tau_{2}, \ldots, \tau_{l}) \}$
are compatible if and only if $t_{j} = \tau_{j}$ for all $j \in
\{ 2, \ldots, l \}$.
\end{rem}

\begin{rem} \label{rem_1.9}
Obviously, relation of the compatibility of two periodic partitions
is reflexive and is symmetrical. However previous note shows, that
generally speaking this relation is not transitive.
\end{rem}

\subsection{Sequences of periodic partitions.}

In this subsection we shall prove a number of statements in order to
analyze the transitivity problem of the compatibility relation
of periodic partitions. The statements we are going to discuss will be used in
further constructions.

Later on the following objects will be necessary for us.

\begin{defn} \label{defn_1.6}
Let a sequence of numbers $\{ n_{i} \in \EssP{X, f} \}_{i \in
\nn} $ is given.

We call a sequence $\{ W^{(n_{i})} \}_{i \in \nn}$ of
periodic partitions of dynamical system $(X, f)$ {\em
regular}, if it satisfies to the following conditions
\begin{itemize}
    \item[1)] $n_{k}$ divides $n_{k+1}$, $k \in \nn$;
    \item[2)] partitions $\{ W^{(n_{k})}_{s_{k}} \}$ and
        $\{ W^{(n_{k+1})}_{s_{k+1}} \}$ are compatible for all $k \in \nn$.
\end{itemize}
\end{defn}

\begin{rem} \label{rem_1.10}
Using Corollary~\ref{corr_1.6}, it is easy to verify that every regular
sequence $\{ W^{(n_{k})} \}_{k \in \nn}$ of periodic
partitions of dynamical system $(X, f)$ complies with the requirement
\begin{itemize}
    \item[3)] periodic partitions $\{ W^{(n_{k})}_{s_{k}} \}$ and
        $\{ W^{(n_{l})}_{s_{l}} \}$ are compatible for all $k, l \in \nn$.
\end{itemize}
Really, let $k$, $l \in \nn$, $k < l$. In accord with
Corollary~\ref{corr_1.6}, partition $\{ W^{(n_{i})}_{s_{i}} \}$ of space
$X$ is the refinement of partition $\{ W^{(n_{i+1})}_{s_{i+1}} \}$ for
every $i \in \nn$. Hence, partition $\{ W^{(n_{l})}_{s_{l}} \}$
is the refinement of partition $\{ W^{(n_{k})}_{s_{k}} \}$.
Using again Corollary~\ref{corr_1.6} we conclude that partitions
$\{ W^{(n_{l})}_{s_{l}} \}$ and $\{ W^{(n_{k})}_{s_{k}} \}$ are compatible.
\end{rem}

\begin{prop} \label{prop_1.6}
Let a sequence of numbers $\{ n_{k} \in \EssP{X, f} \}_{k \in
\nn}$ is given, which complies with requirement 1) of Definition~\ref{defn_1.6}.

There exists a regular sequence $\{ W^{(n_{k})} \}_{k \in \nn}$ of
periodic partitions of dynamical system $(X, f)$.
\end{prop}

\begin{proof}
This statement is easily proved by the inductive application of
Proposition~\ref{prop_1.5}.
\end{proof}

\begin{prop} \label{prop_1.7}
Let a regular sequence $\{ W^{(n_{k})} \}_{k \in \nn}$ of
periodic partitions of dynamical system $(X, f)$  is given.

Let $m_{1}, m_{2} \in \EssP{X, f}$ and $m_{1}$ divides $m_{2}$.

Then
\begin{itemize}
    \item[a)] there exists periodic partition
        $\{ W^{(m_{1})}_{i} \}_{i=0}^{m_{1} -1}$ of length $m_{1}$, which is
        compatible with each of partitions
        $\{ W^{(n_{k})}_{s_{k}} \}$, $k \in \nn$;
    \item[b)] for any periodic partition $\{ W^{(m_{1})}_{i} \}$ complying with
        item a) of the proposition there exists periodic partition
        $\{ W^{(m_{2})}_{j} \}_{j=0}^{m_{2} -1}$ of length $m_{2}$, which is compatible
        with $\{ W^{(m_{1})}_{i} \}$ and with each of partitions
        $\{ W^{(n_{k})}_{s_{k}} \}$, $k \in \nn$.
\end{itemize}
\end{prop}

Proof of Proposition~\ref{prop_1.7} is based on three lemmas.

Let $m_{1}, m_{2}, n \in \EssP{X, f}$, and $m_{1}$ divides $m_{2}$.

Designate by $d_{i}$ the greatest common divisor of numbers $n$ and $m_{i}$,
$i=1, 2$. Let also $D_{i}$ be the least common multiple of numbers $n$ and
$m_{i}$, $i=1, 2$.

\begin{lemma} \label{lemma_1.2}
Let $\{ W^{(m_{1})}_{i} \}_{i=0}^{m_{1} -1}$,
$\{ W^{(m_{2})}_{j} \}_{j=0}^{m_{2} -1}$ and $\{ W^{(n)}_{k} \}_{k=0}^{n-1}$
are periodic partitions of dynamical system $(X, f)$ of lengths $m_{1}$,
$m_{2}$ and $n$, accordingly. Suppose, partitions
$\{ W^{(n)}_{k} \}$ and $\{ W^{(m_{2})}_{j} \}$
are compatible.

If the partitions $\{ W^{(m_{1})}_{i} \}$ and $\{ W^{(m_{2})}_{j} \}$ are compatible,
then the partitions $\{ W^{(m_{1})}_{i} \}$ and $\{ W^{(n)}_{k} \}$ are compatible too.
\end{lemma}

\begin{lemma} \label{lemma_1.3}
Let $d_{1} = d_{2}$.

Let $\{ W^{(m_{1})}_{i} \}_{i=0}^{m_{1} -1}$ and
$\{ W^{(n)}_{k} \}_{k=0}^{n-1}$ are compatible periodic partitions of
dynamical system $(X, f)$ of lengths $m_{1}$ and $n$, respectively.

Then any periodic partition $\{ W^{(m_{2})}_{j} \}_{j=0}^{m_{2} -1}$
of length $m_{2}$, which is compatible with the partition $\{ W^{(m_{1})}_{i} \}$,
also is compatible with the partition $\{ W^{(n)}_{k} \}$.
\end{lemma}

\begin{lemma} \label{lemma_1.4}
Let $D_{1} = D_{2}$.

Let $\{ W^{(m_{1})}_{i} \}_{i=0}^{m_{1} -1}$ and
$\{ W^{(n)}_{k} \}_{k=0}^{n-1}$ are compatible periodic partitions
of dynamical system $(X, f)$ of lengths $m_{1}$ and $n$, respectively.

Then there exists periodic partition
$\{ W^{(m_{2})}_{j} \}_{j=0}^{m_{2} -1}$ of length $m_{2}$, which is
compatible both with the partition $\{ W^{(m_{1})}_{i} \}$ and
with the partition $\{ W^{(n)}_{k} \}$.
\end{lemma}

\begin{proof}[of Lemma~\ref{lemma_1.2}]
Substituting partition $\{ W^{(m_{2})}_{j} \}$ on equivalent, we can regard
that $W^{(m_{2})}_{0} \cap W^{(m_{1})}_{0} \neq \emptyset$ (see.
Remark~\ref{rem_1.6}). Then, applying Corollary~\ref{corr_1.6}, we shall receive
inclusion $W^{(m_{2})}_{0} \subseteq W^{(m_{1})}_{0}$.

Similarly, substituting partition $\{ W^{(n)}_{k} \}$ on equivalent, we shall
regard that $K = W^{(n)}_{0} \cap W^{(m_{2})}_{0} \neq \emptyset$. From
Definition~\ref{defn_1.5} it follows, that
$$
\bigcup_{t \in \zz} f^{t} (K) = \bigcup_{t=0}^{D_{2} -1}
f^{t} \left( W^{(m_{2})}_{0} \cap W^{(n)}_{0} \right) = X \;.
$$

On the other hand, $K = W^{(m_{2})}_{0} \cap W^{(n)}_{0} \subseteq
W^{(m_{1})}_{0} \cap W^{(n)}_{0}$, therefore from Corollary~\ref{corr_1.4} it follows,
that partitions $W^{(m_{1})}_{i}$ and $W^{(n)}_{k}$ are compatible.
\end{proof}

\begin{proof}[of Lemma~\ref{lemma_1.3}]
As $d_{1} =d_{2}$ and $d_{1}$ divides $m_{1}$, then $d_{2}$ divides $m_{1}$.

The replacement of a periodic partition on an equivalent does not
affect the relation of compatibility (see note~\ref{rem_1.6}), therefore we can
assume, that $W^{(n)}_{0} \cap W^{(m_{2})}_{0} \neq \emptyset$ and
$W^{(m_{1})}_{0} \cap W^{(m_{2})}_{0} \neq \emptyset$. Since partitions
$\{ W^{(m_{1})}_{i} \}$ and $\{ W^{(m_{2})}_{j} \}$ are compatible, then
$W^{(m_{2})}_{0} \subseteq W^{(m_{1})}_{0}$ (see Corollary~\ref{corr_1.6}).
Hence, $W^{(n)}_{0} \cap W^{(m_{1})}_{0} \neq \emptyset$.

We designate
$$
A = \bigcup_{s=0}^{D_{2} -1} f^{s} (W^{(m_{2})}_{0} \cap W^{(n)}_{0}) \;.
$$
Taking into account Corollary~\ref{corr_1.3}, to complete proof of lemma
it is enough to check the equality $A = X$.

Consider the pair of compatible partitions $\{ W^{(m_{1})}_{i} \}$ and
$\{ W^{(m_{2})}_{j} \}$. From Definition~\ref{defn_1.5}, Lemma~\ref{lemma_1.1} and
Corollary~\ref{corr_1.6} we receive
$$
W^{(m_{1})}_{0} =
W^{(m_{1})}_{0} \cap \bigcup_{s=0}^{m_{2} -1}
f^{s} \left( W^{(m_{1})}_{0} \cap W^{(m_{2})}_{0} \right) =
\bigcup_{%
\begin{subarray}{c}
    s \in \{ 0, \ldots, m_{2} -1 \} \;, \\
    \Equiv{s}{0}{m_{1}}
\end{subarray}%
}
W^{(m_{2})}_{s} \;.
$$

Consider now the pair $\{ W^{(n)}_{k} \}$ and
$\{ W^{(m_{2})}_{j} \}$ of periodic partitions.

From Lemma~\ref{lemma_1.1} we get
$$
W^{(n)}_{0} \cap A = W^{(n)}_{0} \cap \bigcup_{%
\begin{subarray}{c}
    r \in \{ 0, \ldots, m_{2} -1 \} \;, \\
    \Equiv{r}{0}{d_{2}}
\end{subarray}%
} W^{(m_{2})}_{r} \;.
$$

As $d_{2}$ divides $m_{1}$, then congruence $\Equiv{r}{0}{d_{2}}$ is
the consequence of the congruence $\Equiv{r}{0}{m_{1}}$ and
$$
W^{(m_{1})}_{0} =
\bigcup_{%
\begin{subarray}{c}
    r \in \{ 0, \ldots, m_{2} -1 \} \;, \\
    \Equiv{r}{0}{m_{1}}
\end{subarray}%
}
W^{(m_{2})}_{r}
\subseteq
\bigcup_{%
\begin{subarray}{c}
    r \in \{ 0, \ldots, m_{2} -1 \} \;, \\
    \Equiv{r}{0}{d_{2}}
\end{subarray}%
}
W^{(m_{2})}_{r} \;.
$$

Hence, $W^{(m_{1})}_{0} \cap W^{(n)}_{0} \subseteq W^{(n)}_{0} \cap
A \subset A$. But the set $A$ is $f$--invariant,
therefore
$$
A = \bigcup_{s=0}^{D_{1} -1} f^{s} (A) \supseteq \bigcup_{s=0}^{D_{1} -1}
f^{s} \left( W^{(m_{1})}_{0} \cap W^{(n)}_{0} \right) = X \;.
$$
Last equality is valid, since periodic partitions
$\{ W^{(m_{1})}_{i} \}$ and $\{ W^{(n)}_{k} \}$ are compatible.
\end{proof}

\begin{proof}[of Lemma~\ref{lemma_1.4}]
Taking into account Remark~\ref{rem_1.6}, we shall regard
that $W^{(n)}_{0} \cap W^{(m_{1})}_{0} \neq \emptyset$. Designate
$$
V_{s} = f^{s} \left( W^{(n)}_{0} \cap W^{(m_{1})}_{0} \right) \;, \quad
s = 0, \ldots, D_{1} -1 \;.
$$
Since periodic partitions $\{ W^{(n)}_{k} \}$ and $\{ W^{(m_{1})}_{i} \}$
are compatible and $D_{2} = D_{1}$ on a condition of Lemma, then
$$
X = \bigcup_{s=0}^{D_{1} -1} V_{s} = \bigcup_{s=0}^{D_{2} -1} V_{s}
$$
and $\{ V_{s} \}_{s=0}^{D_{2} -1}$ is a periodic partition of
dynamical system $(X, f)$ of length $D_{2}$.

We set
$$
W^{(m_{2})}_{j} = \bigcup_{%
\begin{subarray}{c}
    s \in \{ 0, \ldots, D_{2} -1 \} \;, \\
    \Equiv{s}{j}{m_{2}}
\end{subarray}%
}
V_{s} \;, \quad j = 0, \ldots, m_{2} -1 \;.
$$
We receive periodic partition $\{ W^{(m_{2})}_{j} \}_{j=0}^{m_{2} -1}$
of dynamical system $(X, f)$ of length $m_{2}$ (see proof of
Proposition~\ref{prop_1.1}).

We shall prove now, that this periodic partition is compatible with each of
partitions $\{ W^{(m_{1})}_{i} \}$ and $\{ W^{(n)}_{k} \}$.

On construction we have $V_{0} = W^{(m_{1})}_{0} \cap W^{(n)}_{0} \subseteq
W^{(m_{2})}_{0}$, therefore $V_{0} = W^{(m_{1})}_{0} \cap W^{(n)}_{0}
\subseteq W^{(m_{2})}_{0} \cap W^{(n)}_{0}$. Now from property~(iv) of
Definition~\ref{defn_1.1} and from Corollary~\ref{corr_1.4} it
follows, that periodic partitions $\{ W^{(m_{2})}_{i} \}$ and
$\{ W^{(n)}_{k} \}$ are compatible.

Similarly, on construction $V_{0} = W^{(m_{1})}_{0} \cap W^{(n)}_{0}
\subseteq W^{(m_{1})}_{0}$ and $V_{0} \subseteq W^{(m_{2})}_{0}$.
Hence, $V_{0} \subseteq W^{(m_{1})}_{0} \cap W^{(m_{2})}_{0}$ and from
Corollary~\ref{corr_1.4} we receive, that partitions $\{ W^{(m_{1})}_{i} \}$ and
$\{ W^{(m_{2})}_{j} \}$ are compatible.
\end{proof}

\begin{proof}[of Proposition~\ref{prop_1.7}]

{\bf a)}
Let $d^{1}_{k}$ is the greatest common divisor of numbers $n_{k}$ and
$m_{1}$. From condition 1) of Proposition it follows, that $n_{k+l}$
is divided by $d^{1}_{k}$ for all $l \in \nn$. Therefore
$$
d^{1}_{1} \leq d^{1}_{2} \leq \ldots \leq d^{1}_{k} \leq \ldots \;.
$$

On the other hand, $d^{1}_{k} \leq m_{1}$ for all $k \in \nn$.
Hence, there exists $l \in \nn$, such that
$d^{1}_{k} = d^{1}_{l}$ when $k \geq l$.

By using Proposition~\ref{prop_1.5},
we find periodic partition $\{ W^{(m_{1})}_{i} \}_{i=0}^{m_{1} -1}$ of
length $m_{1}$, which is compatible with the partition
$\{ W^{(n_{l})}_{s_{l}} \}$. Show, that for every $k \in \nn$
this partition is compatible with the
partition $\{ W^{(n_{k})}_{s_{k}} \}$.

Let $k > l$. Then $d^{1}_{k} = d^{1}_{l}$. From Remark~\ref{rem_1.10} it
follows, that partitions $\{ W^{(n_{l})}_{s_{l}} \}$ and
$\{ W^{(n_{k})}_{s_{k}} \}$ are compatible. We apply Lemma~\ref{lemma_1.3}
to periodic partitions $\{ W^{(m_{1})}_{i} \}$,
$\{ W^{(n_{l})}_{s_{l}} \}$ and $\{ W^{(n_{k})}_{s_{k}} \}$, and conclude
that partitions $\{ W^{(m_{1})}_{i} \}$ and $\{ W^{(n_{k})}_{s_{k}} \}$
are compatible.

Let now $k < l$. Again from Remark~\ref{rem_1.10}
we derive, that partitions $\{ W^{(n_{l})}_{s_{l}} \}$ and
$\{ W^{(n_{k})}_{s_{k}} \}$ are compatible. We apply now Lemma~\ref{lemma_1.2}
to periodic partitions $\{ W^{(n_{k})}_{s_{k}} \}$,
$\{ W^{(n_{l})}_{s_{l}} \}$ and $\{ W^{(m_{1})}_{i} \}$
and conclude that partitions
$\{ W^{(m_{1})}_{i} \}$ and $\{ W^{(n_{k})}_{s_{k}} \}$ are compatible.

{\bf b)}
Let $d^{2}_{k}$ is the greatest common divisor of numbers $n_{k}$ and
$m_{2}$. Repeating argument from item~a), we consequence that there exists
$\tau \in \nn$, such that $d^{2}_{k} = d^{2}_{\tau}$ when $k \geq \tau$.

To complete the proof of item~b) it is enough to us now to find periodic
partition of length $m_{2}$, which is compatible both with
$\{ W^{(n_{\tau})}_{s_{\tau}} \}$ and $\{ W^{(m_{1})}_{i} \}$. Then by repetition of
argument from item~a) we shall prove, that it is compatible with every
$\{ W^{(n_{k})}_{s_{k}} \}$, $k \in \nn$.

Consider triple of numbers $m_{1}$, $m_{2}$, $n_{\tau} \in \EssP{X, f}$.
Designate for convenience by $d_{k}$ and $D_{k}$ the greatest
common divisor and the least common multiple of numbers $n_{\tau}$ and $m_{k}$,
$k = 1, 2$.

As $m_{1}$ divides $m_{2}$ on the condition of Proposition, then $d_{1}$
divides $d_{2}$ and $D_{1}$ divides $D_{2}$. Set
$$
m = \frac{m_{1}}{d_{1}} \cdot d_{2} \;.
$$

It is clear, that $m_{1}$ divides $m$, since the number $d_{2} /d_{1}$ is integer.
On the other hand,
$$
\frac{m_{2}}{m} = \frac{m_{2}}{m_{1}} \cdot \frac{d_{1}}{d_{2}} =
\left( \frac{m_{2} n_{\tau}}{d_{2}} \right)
\left( \frac{m_{1} n_{\tau}}{d_{1}} \right)^{-1} =
\frac{D_{2}}{D_{1}} \in \nn \;,
$$
that is $m$ divides $m_{2}$.

Let $d$ and $D$ are the greatest common divisor and the least common multiple
of numbers $m$ and $n_{\tau}$.

Obviously, $d_{1}$ divides $d$ and $d$ divides $d_{2}$. Since the number
$m_{1} /d_{1}$ is integer, then $m$ is divided by $d_{2}$ and $d_{2}$ is
the common divisor of numbers $m$ and $n_{\tau}$. Hence, $d_{2} = d$.

On the other hand,
$$
D = \frac{m n_{\tau}}{d} = \frac{m_{1} d_{2}}{d_{1}}
\cdot \frac{n_{\tau}}{d_{2}} = \frac{m_{1} n_{\tau}}{d_{1}}
= D_{1} \;.
$$

Proposition~\ref{prop_1.1} gives us the inclusion $m \in \EssP{X, f}$, since
$m$ divides $m_{2}$ and $m_{2} \in \EssP{X, f}$.

Applying Lemma~\ref{lemma_1.4} to numbers $m_{1}$, $m$, $n_{\tau} \in
\EssP{X, f}$ and to compatible periodic partitions
$\{ W^{(m_{1})}_{i} \}$ and $\{ W^{(n_{\tau})}_{s_{\tau}} \}$, we shall find
periodic partition $\{ W^{(m)}_{k} \}_{k=0}^{m-1}$ of length $m$, which is
compatible with each of partitions $\{ W^{(m_{1})}_{i} \}$ and
$\{ W^{(n_{\tau})}_{s_{\tau}} \}$.

Having used Proposition~\ref{prop_1.5}, we find periodic partition
$\{ W^{(m_{2})}_{j} \}_{j=0}^{m_{2} -1}$ of length $m_{2}$, compatible with
partition $\{ W^{(m)}_{k} \}$. We apply Lemma~\ref{lemma_1.3} to
numbers $m$, $m_{2}$, $n_{\tau}$ and periodic partitions
$\{ W^{(m)}_{k} \} $, $ \{ W^{(m_{2})}_{j} \}$,
$\{ W^{(n_{\tau})}_{s_{\tau}} \}$, and conclude that partitions
$\{ W^{(m_{2})}_{j} \}$ and $\{ W^{(n_{\tau})}_{s_{\tau}} \}$ are compatible.

In accord with Corollary~\ref{corr_1.6}, partition $\{ W^{(m_{2})}_{j} \}$
of space $X$ is the refinement of the partition $\{ W^{(m)}_{k} \}$.
Consequently it all the more is the refinement of the partition
$\{ W^{(m_{1})}_{i} \}$. Again applying Corollary~\ref{corr_1.6}
we conclude, that partitions $\{ W^{(m_{1})}_{i} \}$ and
$\{ W^{(m_{2})}_{j} \}$ are compatible.

Proposition~\ref{prop_1.7} is completely proved.
\end{proof}

\begin{rem} \label{rem_1.11}
Let a regular sequence
$ \{ W^{(n_{i})}_{s_{i}} \}_{s_{i} =0}^{n_{i} -1} $ of periodic partitions
is given and
$$
\bigcap_{i \in \nn} W^{(n_{i})}_{r_{i}} \neq \emptyset
$$
for certain sequence $\{ r_{i} \}_{i \in \nn}$. Then from
Corollary~\ref{corr_1.6} we receive
$$
W^{(n_{1})}_{r_{1}} \supseteq W^{(n_{2})}_{r_{2}} \supseteq
\ldots W^{(n_{i})}_{r_{i}} \supseteq \ldots \;.
$$

Next, if $K \subseteq X$ is the closed subset of a compactum $X$, such that
$K \cap W^{(n_{i})}_{r_{i}} \neq \emptyset$ for every $i \in \nn $,
then a sequence of nested compact sets
$$
\left(
K \cap W^{(n_{1})}_{r_{1}}
\right)
\supseteq \ldots \supseteq
\left(
K \cap W^{(n_{i})}_{r_{i}}
\right)
\supseteq \ldots
$$
has nonempty intersection, in other words
$$
K \cap
\left(
\bigcap_{i \in \nn} W^{(n_{i})}_{r_{i}}
\right)
\neq \emptyset \;.
$$
\end{rem}

\begin{defn} \label{defn_1.7}
Two regular sequences $\{W^{(n_{i})}_{s_{i}}\}_{s_{i}=0}^{n_{i}-1}$,
$i \in \nn$, and $\{W^{(m_{j})}_{\tau_{j}}\}_{\tau_{j}=0}^{m_{j}-1}$,
$j \in \nn$, of periodic partitions of dynamical system $(X, f)$ are called
{\em compatible}, if periodic partitions
$\{ W^{(n_{i})}_{s_{i}} \}$ and $\{ W^{(m_{j})}_{\tau_{j}} \}$ are compatible for
all $i$, $j \in \nn$.
\end{defn}

Immediately from Proposition~\ref{prop_1.7} we receive the following

\begin{prop} \label{prop_1.8}
Let $\{W^{(n_{i})}_{s_{i}}\}_{s_{i}=0}^{n_{i}-1}$, $i \in \nn$, is
regular sequence of periodic partitions of dynamical system
$(X, f)$.

Let $m_{j} \in \EssP{X, f}$, $j \in \nn$, and $m_{j}$ divides $m_{j+1}$ for
every $j \in \nn$.

Then there exists a regular sequence
 $\{W^{(m_{j})}_{\tau_{j}}\}_{\tau_{j}=0}^{m_{j}-1}$, $j \in \nn $, of
periodic partitions of dynamical system $(X, f)$, which is compatible with
sequence $\{ W^{(n_{i})}_{s_{i}} \}$, $i \in \nn$.
\end{prop}

\begin{rem} \label{rem_1.12}
It follows immediately from Corollary~\ref{corr_1.1}, that if
dynamical system $(X, f)$ is indecomposable, then any two regular
sequences of periodic partitions of this dynamical
system are compatible.
\end{rem}

\begin{rem} \label{rem_1.13}
Let a sequence of numbers $\{ n_{i} \in \EssP{X, f} \}_{i \in
\nn}$ is given, which satisfies to condition 1) of Definition~\ref{defn_1.6}.

If dynamical system $(X, f)$ is not indecomposable, then there exist
two incompatible periodic partitions $W^{(n_{1})}$ and
$\widetilde{W}^{(n_{1})}$ of dynamical system $(X, f)$
(see Proposition~\ref{prop_1.3} and Remark~\ref{rem_1.5}).

Obviously, regular sequences $\{ W^{(n_{i})} \}_{i \in \nn}$ and
$\{ \widetilde {W}^{(n_{i})} \}_{i \in \nn}$ of periodic partitions, built
from partitions $W^{(n_{1})}$ and $\widetilde {W}^{(n_{1})}$ accordingly
with the help of inductive application of Proposition~\ref{prop_1.5},
are not compatible.
\end{rem}

\subsection{Periodic partitions and returnability of trajectories of
dynamical system.}

Consider some dynamical system $(X, f)$.

\begin{lemma} \label{lemma_1.5}
Let for certain recurrent point $x \in X$ there exists a closed
neighborhood $U$, which satisfies to the following property: there exists $n \in \nn$, such
that
$$
\bigcup_{k \in \zz} f^{kn} (x) \subset U \;.
$$
Let
\begin{equation} \label{eq_6}
m = \min \{ n \in \nn \; | \; f^{nk} (x) \in U \; \forall k \in \zz \} \;.
\end{equation}

Then dynamical system $(\Cl{\Orb_{f} (x)}, f) $ has
periodic partition $\{ W_{i} \}_{i = 0}^{m-1}$ of length
$m$, such that $x \in W_{0} \subseteq U$.
\end{lemma}

\begin{proof}
Under Birkgoff Theorem the set $\Cl{\Orb_{f} (x)}$ is minimal set of d. s.
$(X, f)$. In particular, the space $\Cl{\Orb_{f} (x)}$ is compact.

We assume, that homeomorphism $f$ is given on space $\Cl{\Orb_{f} (x)}$ with
the topology induced from $X$ and we shall consider in this topology all sets
which will arise in the proof.

Designate
$$
W = \Cl{\bigcup_{k \in \zz} f^{mk} (x)} \,, \quad \widetilde{W} = \Int W \;.
$$
Consider two families of sets
$$
\begin{array}{lll}
    W_{0} = W \,, & W_{i} = f^{i} (W) = f (W_{i-1}) \,, & i = 1, \ldots, m-1
        \,; \\
    & \widetilde{W}_{i} = \Int W_{i} \,, & i = 0, \ldots, m-1
        \,. \\
\end{array}
$$
It is clear, that for the family $\{ W_{i} \}$ conditions~(ii) and~(iv)
of Definition~\ref{defn_1.1} are fulfilled. All sets $W_{i}$ are closed,
therefore condition~(i) is an immediate corollary of conditions~(ii) and~(iii).

So, in order to prove Lemma it is enough to verify condition~(iii) of
Definition~\ref{defn_1.1}.

At first we shall show, that $\widetilde{W}_{i} \neq \emptyset$, $i = 0, 1, \ldots,
m-1$. Really, under Baire category Theorem%
\footnote {the Baire category Theorem is applicable in spaces, complete on Cech,
and it is known that any compact set is space, complete on Cech
(see~\cite{Engelking})}
it follows from (iv) that at least one of the sets $\{ W_{i} \}$
is not the set of I-st category (and has nonempty interior in
$\Cl{\Orb_{f} (x)}$). Since $f$ is the homeomorphism, then all $W_{i}$ have
nonempty interior, that is $\widetilde{W}_{i} \neq \emptyset$.

Now we shall check the relation $\widetilde{W}_{i} \cap \widetilde{W}_{j} =
\emptyset$, $i \neq j$.

Assume, that it is not the case. Let $V = \widetilde{W}_{i} \cap
\widetilde{W}_{j} \neq \emptyset$ and  $i < j$ for a determinancy. Since the
set $\bigcup_{k \in \zz} f^{km+i} (x)$ is dense in $\widetilde{W}_{i}$ on
construction, then $f^{ms+i} (x) \in V $ for some $s \in \zz$.

$$
f^{ms+i} (x) \in V \subseteq \widetilde{W}_{j} \subseteq W_{j} =
\Cl{\bigcup_{k \in \zz} f^{km+j} (x)} \;,
$$
therefore there exists a sequence $\{ l_{r} \}_{r \in \nn}$ of integers,
such that $f^{ml_{r} +j} (x) \rightarrow f^{ms+i} (x)$ when $r \rightarrow
\infty$.
This implies, that for every $k \in \zz$ we have
$$
f^{m (l_{r} -s+k) + (j-i)} (x) \rightarrow f^{km} (x) \;.
$$
That is
$$
\bigcup_{k \in \zz} f^{km} (x) \subseteq \Cl{\bigcup_{k \in \zz}
f^{km + (j-i)} (x)} = W_{j-i} \;.
$$
Passing to closures, we shall receive $W_{0} \subseteq W_{j-i}$. Changing
roles of $i$ and $j$, we shall receive inverse inclusion.
Therefore $W_{0} = W_{j-i}$.

Remark, that
$$
f^{(j-i) k} (x) \in f^{(j-i) k} (W_{0}) = f^{(j-i) (k-1)} \circ f^{j-i} (W_{0})
= f^{(j-i) (k-1)} (W_{0}) = \cdots =W_{0} \subseteq U
$$
when $k \geq 0$ and
$$
\begin{array} {cclcc}
    f^{-|k|(j-i)}(x) & \in & f^{-|k|(j-i)} (W_{0}) =
        f^{-|k|(j-i)} (W_{j-i}) = & & \\
    & & = f^{(-|k|+1)(j-i)} \circ f^{-(j-i)} (W_{j-i})
        = f^{(-|k|+1)(j-i)}(W_{0}) = \cdots =W_{0} & \subseteq & U \\
\end{array}
$$
for $k < 0$ (we remind, that the set $U$ is closed on condition of Lemma). That
is $f^{k (j-i)}(x) \in U$ for all $k \in \zz$.

On our supposition $0<j-i<m$, so we have received contradiction with
the choice of $m$ (see relation~(\ref{eq_6})). Hence,
$\widetilde{W}_{i} \cap \widetilde{W}_{j} = \emptyset$ when $i \neq j$.

We shell check now equalities $\widetilde{W}_{i} = W_{i}$, $i = 0, 1, \ldots,
m-1$.

It is easy to see, that
$$
\widetilde{W}_{i} = \Int W_{i} = \Int (f (W_{i-1})) = f (\Int (W_{i-1})) =
f (\widetilde{W}_{i-1}) \, \quad i = 1, \ldots, m-1 \;,
$$
and $\widetilde{W}_{0} = f (\widetilde{W}_{m-1})$, since $f$ is
homeomorphism. Designate
$$
Q = \bigcup_{i=0}^{m-1} \widetilde{W}_{i} = \bigcup_{i=0}^{m-1}
f^{i} (\widetilde{W}_{0}) \;.
$$
The set $Q $ is open invariant subset of dynamical system
$(\Cl{\Orb_{f} (x)}, f)$. Let $K = \Cl{\Orb_{f} (x)} \setminus Q$.
Obviously, $K$ is the closed invariant subset of this
dynamical system.

The set $ \Cl{\Orb_{f} (x)}$ is minimal, therefore either $K = \emptyset$
or $K = \Cl{\Orb {f} (x)}$. We already have proved, that $Q \neq \emptyset$,
hence $K = \emptyset$ and
$$
\Cl{\Orb_{f} (x)} = \bigcup_{i=0}^{m-1} \widetilde{W}_{i} \;.
$$

Sets from the family $\{ \widetilde{W}_{i} \}$ are pairwise disjoint, so
all of them are open-closed and $\widetilde{W}_{i} = W_{i}$,
$i = 0, 1, \ldots, m-1$. On proved above from this immediately follows
condition~(iii) of definition~\ref{defn_1.1}.

Lemma is completely proved.
\end{proof}

\section{Ultranatural numbers and subsets of natural
numbers.}

\subsection{Ultranatural numbers.}

\begin{defn} \label{defn_2.1}
Let $\ss \subset \nn$ is the set of all prime numbers ordered by
increment. A sequence
$$
N = \left( N_{2}, N_{3}, \ldots, N_{p}, \ldots \right) \;, \qquad
N_{p} \in \zz_{+} \cup \{ \infty \} \;, \quad p \in \ss \;,
$$
is called {\em ultranatural number}.
\end{defn}

The set of all ultranatural numbers we shall designate by $ \Sigma $.

We introduce the relation of partial order on $\Sigma$. Say that
$$
M \leq N \;, \quad M, N \in \Sigma \;,
$$
if $M_{p} \leq N_{p}$ for every $p \in \ss$ (we shall regard that
$k \leq \infty$ for any $k \in \zz_{+}$). Elementary immediate
verification shows the correctness of this definition.

Next, we introduce  binary operation on $ \Sigma $. For $M = (M_{p})$ and $N =
( N_{p})$ we set
$$
M \cdot N = K = (K_{p}) \;;
$$
$$
K_{p} = \left\{
    \begin{array}{cl}
        M_{p} + N_{p} \;, & \mbox{if } M_{p} \neq \infty \mbox{ and } N_{p} \neq \infty \;, \\
        \infty \;, & \mbox{otherwise} \;. \\
    \end{array}
\right\} \;, \quad
p \in \ss \;.
$$
It is trivially checked, that $(\Sigma, \cdot)$ is a semigroup with unity $E =
(E_{p} =0)$.

\begin{rem}
Easy immediate verification shows, that the equation $M \cdot X =
N$ has a solution in $(\Sigma, \cdot)$ only when $M \leq
N$.

However, this equation can have more than one solution.
\end{rem}

\begin{example}
Let $M = N = (N_{p})$,
$$
N_{p} = \left\{
    \begin{array}{cl}
        \infty \;, & \mbox{when } p = 2 \;, \\
        0 \;, & \mbox{for } p \neq 2 \;. \\
    \end{array}
\right.
$$

Then $X^{(n)} = (X^{(n)}_{p})$,
$$
X^{(n)}_{p} = \left\{
    \begin{array}{cl}
        n \;, & \mbox{for } p = 2 \;, \\
        0 \;, & \mbox{when } p \neq 2 \;. \\
    \end{array}
\right.
$$
is a solution of the equation $M \cdot X = N$ for every $n \in \zz_{+} \cup
\{ \infty \}$.
\end{example}

In just the same way as in the semigroup $(\nn, \cdot)$ we can well define
the greatest common divisor and least common
multiple for any two $M$, $N \in
\Sigma$. It is easy to see, that
\begin{eqnarray}
\gcd{M}{N} & = & (d_{p}) \;, \quad d_{p} = \min (M_{p}, N_{p}) \;, \quad p \in
\ss \;, \label{eq_7} \\
\lcm{M}{N} & = & (D_{p}) \;, \quad D_{p} = \max (M_{p}, N_{p}) \;, \quad p \in
\ss \;. \label{eq_8}
\end{eqnarray}
Here we use the following agreements:
$$
\max (a, \infty) = \infty \;, \quad
\min (a, \infty) = a \;, \quad
a \in \zz_{+} \cup \{ \infty \} \;.
$$

We define monomorphism $ \Phi_{0} : (\nn, \cdot) \rightarrow (\Sigma,
\cdot)$.

Let $n \in \nn$. Consider factorization
$$
n = p_{1}^{\alpha_{1}} \ldots p_{k}^{\alpha_{k}}
$$
of the number $n$ on prime factors (we regard that $p_{i} \neq p_{j}$
when $i \neq j$). Set $\Phi_{0} (n) = (\Phi_{0} (n)_{p})$,
$$
\Phi_{0} (n)_{p} = \left\{
    \begin{array}{cl}
        \alpha_{i} \;, & \mbox{when } p = p_{i} \in \{ p_{1}, \ldots, p_{k} \} \;,
        \\
        0 \;, & \mbox{otherwise} \;.
    \end{array}
\right.
$$

\begin{rem}
Mark that for all $m$, $n \in \nn$ we have
\begin{eqnarray*}
    \Phi_{0} (\gcd{m}{n}) & = & \gcd{\Phi_{0} (m)}{\Phi_{0} (n)} \;, \\
    \Phi_{0} (\lcm{m}{n}) & = & \lcm{\Phi_{0} (m)}{\Phi_{0} (n)} \;.
\end{eqnarray*}
\end{rem}

\subsection{Regular subsets of natural numbers.}

Let $A \subseteq \nn$. For every $p \in \ss$ let
$$
\Phi (A)_{p} = \sup \{ k \in \zz_{+} \; | \; \exists a \in A: p^{k} \mbox{ divides } a
\} = \sup_{a \in A} \Phi_{0} (a)_{p} \;.
$$

Designate by $\Phi$ the map $\Phi : A \mapsto \Phi (A) = (\Phi (A)_{p})$ from
class of all nonempty subsets of $ \nn $ to the set
$\Sigma$ of ultranatural numbers.

\begin{rem} \label{rem_2.3}
It is easily checked that order relation on $\Sigma$ defined above turns
map $\Phi$ into isotonic map, that is for any $A,
B \subseteq \nn$ inclusion $A \subseteq B$ implies $\Phi (A) \leq \Phi (B)$.
\end{rem}

\begin{example} \label{ex_2.2}
Let $A = \{ a \}$ is a singleton. From definition it easily follows,
that $ \Phi (\{ a \}) = \Phi_{0} (a)$.
\end{example}

\begin{example} \label{ex_2.3}
Let $A = \{ a_{1}, \ldots, a_{j} \} $ is a finite subset of $\nn$.
Consider factorizations of numbers $a_{1}, \ldots, a_{j}$ on prime factors
$$
a_{i} = \prod_{p \in \ss} p^{n_{p} (i)} \;, \quad i = 1, \ldots, j
$$
(here $n_{p} (i) \in \zz_{+} $, $p \in \ss$, $i = 1, \ldots, j$).

By definition
$$
\Phi (A)_{p} = \max \{ n_{p} (1), \ldots, n_{p} (j) \} \;, \quad p \in \ss \;,
$$
in other words $ \Phi (\{ a_{1}, \ldots, a_{j} \}) = \Phi (\{ D \}) $, where $D \in \nn$ is
the least common multiple of numbers $a_{1}, \ldots, a_{j}$.
\end{example}

\begin{rem} \label{rem_2.4}
Let $A \subseteq \nn$. Immediately from definition follows, that
$\Phi (\{ a \}) \leq \Phi (A)$ for every $a \in A$.
\end{rem}

\begin{rem} \label{rem_2.5}
Relations~(\ref{eq_7}) and~(\ref{eq_8}) imply, that for any nonempty
$A$, $B \subseteq \nn$
$$
\Phi (A \cup B) = \lcm{\Phi (A)}{\Phi (B)} \;.
$$
In addition if $A \cap B \neq \emptyset$, then
$$
\Phi (A \cap B) \leq \gcd{\Phi (A)}{\Phi (B)} \;.
$$
\end{rem}

\begin{defn} \label{defn_2.2}
We call a nonempty subset $A \subseteq \nn$ {\em regular}, if it
satisfies to the following conditions:
\begin{itemize}
    \item[(i)] if $a \in A$  and $d \in \nn$ divides $a$, then $d \in A$;
    \item[(ii)] for any $a, b \in A$ their least common multiple $D$
        also is contained in $A$.
\end{itemize}

We designate the family of all regular sets by $\aa$.
\end{defn}

\begin{rem} \label{rem_2.6}
It follows from Propositions~\ref{prop_1.1} and~\ref{prop_1.2}, that for any
dynamical system $(X, f)$ the set $\EssP{X, f}$ is regular.
\end{rem}

\begin{lemma} \label{lemma_2.1}
Let $A \in \aa$. Then
$$
A = \{ a \in \nn \; | \; \Phi (\{ a \}) \leq \Phi (A) \} =
\{ a \in \nn \; | \; \Phi_{0} (a) \leq \Phi (A) \} \;.
$$
\end{lemma}

\begin{proof}
Let $a \in \nn $ and $ \Phi (\{ a \}) \leq \Phi (A) $. Suppose
$$
a = p_{1}^{n_{1}} \ldots p_{k}^{n_{k}}
$$
is the factorization of $a$ on prime factors. By definition of
$\Phi$ there exist such $b_{1}, \ldots, b_{k} \in A$, that
$p_{i}^{n_{i}}$ divides $b_{i}$, $i = 1, \ldots, k$.
Assume $b$ is the least common multiple of numbers $b_{1}, \ldots, b_{k}$. Then $a$
divides $b$. But $b \in A$ by definition of regular set.
Hence and $a \in A$. That is
$$
A \supseteq \{ a \in \nn \; | \; \Phi (\{ a \}) \leq \Phi (A) \} \;.
$$

Inverse inclusion
$$
A \subseteq \{ a \in \nn \; | \; \Phi (\{ a \}) \leq \Phi (A) \} \;.
$$
follows from Remark~\ref{rem_2.4}.
\end{proof}

An immediate corollary of Lemma~\ref{lemma_2.1} is following
\begin{prop} \label{prop_2.1}
Mapping
$$
\Phi |_{\aa}: \aa \rightarrow \Sigma
$$
is bijective.
\end{prop}

\begin{rem}
Let $A, B \in \aa$. Obviously, $ 1 \in A \cap B \neq \emptyset$. From
Lemma~\ref{lemma_2.1} the relation immediately follows
$$
\Phi (A \cap B) = \gcd{\Phi (A)}{\Phi (B)} \;.
$$
\end{rem}

\begin{defn} \label{defn_2.3}
Let $A \subseteq \nn$. Call a sequence $\{ a_{i} \in A \}_{i
\in \nn}$ {\em regular}, if $a_{i}$ divides $a_{i+1}$ for every
$i \in \nn$.
\end{defn}

\begin{rem} \label{rem_2.8}
It follows from Remark~\ref{rem_2.3}, that for any $A \subseteq \nn$ and
any regular sequence $\{ a_{i} \in A \}_{i \in \nn}$ we have the inequality
$$
\Phi (\{ a_{i} \; | \; i \in \nn \}) \leq \Phi (A) \;.
$$
\end{rem}

\begin{prop} \label{prop_2.2}
Let $A \in \aa$. Then there exists a regular sequence $\{ b_{i}
\in A \}_{i \in \nn}$, such that $\Phi (\{ b_{i} \; | \; i \in \nn \}) = \Phi (A)$.
\end{prop}

\begin{proof}
Since $A \subseteq \nn$ is at most the enumerable set, we can
enumerate all elements of $A$ with the help of natural numbers,
$A = \{a_{1}, a_{2}, \ldots\}$. Let $b_{1} = a_{1}$, $b_{i}$
the least common multiple of numbers $a_{i}$ and
$b_{i-1}$ for $i > 1$.

It is clear, that
$$
\Phi (\{ b_{i} \}) \geq \Phi (\{ a_{1}, \ldots, a_{i} \}) \;, \quad i \in \nn \;,
$$
Therefore
$$
\Phi (\{ b_{i} \; | \; i \in \nn \}) \geq \Phi (\{ a_{i} \; | \; i \in \nn \}) = \Phi (A) \;.
$$

On the other hand, $b_{i} \in A$, $i \in \nn$, since $A$ is
regular. Consequently, $\{ b_{i} \}_{i \in \nn} \subseteq A$ and
$$
\Phi (\{ b_{i} \; | \; i \in \nn \}) \leq \Phi (A) \;.
$$

In order to complete the proof it suffices to note, that on construction $b_{i}$
divides $b_{i+1}$, $i \in \nn$, hence the sequence
$\{ b_{i} \}_{i \in \nn}$ is regular.
\end{proof}

\begin{prop} \label{prop_2.3}
Let $A \subseteq \nn$. Suppose two regular sequences
$\{ a_{i} \in A \}_{i \in \nn}$ and $\{ b_{j} \in A \}_{j \in \nn}$ are given.

The following conditions are equivalent:
\begin{itemize}
    \item[1)] $ \Phi (\{ a_{i} \; | \; i \in \nn \}) \leq \Phi (\{ b_{j} \; | \; j \in
        \nn \}) $;
    \item[2)] for every $i \in \nn$ there exists $j \in \nn$, such that
        $a_{i}$ divides $b_{j}$.
\end{itemize}
\end{prop}

\begin{proof}
{\bf 1)}
Let $\Phi (\{ a_{i} \; | \; i \in \nn \}) \leq \Phi (\{ b_{j} \; | \; j \in
\nn \})$.

Consider the factorization
$$
a = a_{i} = p_{1}^{\alpha_{1}} \ldots p_{k}^{\alpha_{k}}
$$
of $a_{i}$ on prime factors.

From Example~\ref{ex_2.2} and Remark~\ref{rem_2.4} it follows, that for $m = 1,
\ldots, k$
$$
\alpha_{m} = \Phi (\{ a \})_{p_{m}} \leq \Phi (\{ a_{i} \; | \; i \in
\nn \})_{p_{m}} \leq \Phi (\{ b_{j} \; | \; j \in \nn \})_{p_{m}} \;.
$$
Therefore, for every $p_{m}$, $m = 1, \ldots, k$, there exists $j_{m} \in
\nn$, such that $p_{m}^{\alpha_{m}}$ divides $b_{j_{m}}$.

Without loss of generality, we shall suppose that
$$
j_{1} \leq j_{2} \leq \ldots \leq j_{k} \;.
$$

Sequence $\{ b_{j} \}$ is regular, therefore $b_{j_{m}}$ divides
$b_{j_{k}}$, $m = 1, \ldots, k$. Hence, $p_{m}^{\alpha_{m}}$ divides
$b_{j_{k}}$, $m = 1, \ldots, k$.

Numbers $p_{1}, \ldots, p_{k}$ are relatively prime on construction, therefore
$p_{1}^{\alpha_{1}} \ldots p_{k}^{\alpha_{k}} = a$ divides
$b_{j_{k}} $.

{ \bf 2)}
Let now assume that condition 2) of Proposition~\ref{prop_2.3} is satisfied.

Suppose that for some $p \in \ss$ the inequality is fulfilled
$$
\Phi (\{ a_{i} \; | \; i \in \nn \})_{p} \gneqq \Phi (\{ b_{j} \; | \; j \in
\nn \})_{p} \;.
$$
Then $n = \Phi (\{ b_{j} \; | \; j \in \nn \})_{p} < \infty$ and on definition
of $ \Phi $
$$
b_{j} = p^{n_{j}} \widetilde {b}_{j} \, \quad n_{j} \leq n \, \quad
\gcd{p}{\widetilde {b}_{j}} = 1
$$
for every $j \in \nn$.

On the other hand, $\Phi (\{ a_{i} \; | \; i \in \nn \})_{p} \geq n + 1$.
Hence there exists $i_{0} \in \nn$, such that $p^{n+1}$ divides
$a_{i_{0}}$. We take an advantage now of condition~2) of Proposition~\ref{prop_2.3} and
find $j_{0} \in \nn$, such that $a_{i_{0}}$ divides $b_{j_{0}}$. All the more,
$p^{n+1}$ divides $b_{j_{0}}$.

However, on construction $\gcd{p^{n+1}}{b_{j_{0}}} = p^{n_{j_{0}}} \leq
p^{n}$.

The obtained contradiction proves that
$$
\Phi (\{ a_{i} \; | \; i \in \nn \})_{p} \leq \Phi (\{ b_{j} \; | \; j \in
\nn \})_{p} \, \quad p \in \ss \;,
$$
That is condition~1) of Proposition~\ref{prop_2.3} is valid.
\end{proof}

\begin{corr} \label{corr_2.1}
Let $A \subseteq \nn$. Assume $\{ a_{i} \in A \}_{i \in \nn}$ is regular
sequence.

For any subsequence $\{ b_{j} \}_{j \in \nn}$ of
$\{ a_{i} \}_{i \in \nn}$ the equality is valid
$$
\Phi (\{ a_{i} \; | \; i \in \nn \}) = \Phi (\{ b_{j} \; | \; j \in \nn \}) \;.
$$
\end{corr}

\section{Odometers and connected constructions.}

\subsection{Definition of odometer.}

We fix regular infinitely growing sequence
$ \{ n_{i} \in \nn \}_{i \in \nn} $.

Consider a sequence of finite cyclic groups $\zz_{n_{i}} =
\zz / n_{i} \zz$ and group homomorphisms
$$
\varphi_{i} : \zz_{n_{i+1}} \rightarrow \zz_{n_{i}} \;,
$$
$$
\varphi_{i} : 1 \mapsto 1 \;.
$$
Take the inverse limit $A = \LimInv_{i \rightarrow \infty}
\zz_{n_{i}}$ of this sequence of groups and homomorphisms.
We shall obtain an Abelian group $(A, +)$.

We endow each set $\zz_{n_{i}} = \{ 0, 1, \ldots, n_{i} -1 \}$
with the discrete topology. Each of maps $\varphi_{i}$ is
continuous in this topology. Space $A$ with the topology $\Tau$
of inverse limit is homeomorphic to Cantor set $\Gamma$.

It is easy to see, that in the group $(A, +)$ operation of addition and passage
to opposite element are continuous in the topology $\Tau$, thus $A$
turns to be the continuous group.

\begin{rem}
We remind, that the inverse limit $A = \LimInv_{i \rightarrow \infty}
\zz_{n_{i}}$ could be imagined as a subset
\begin{equation} \label{eq_9}
A = \{ \vec{a} = (a_{i} \in \zz_{n_{i}}) \; | \; \varphi_{i} (a_{i+1}) = a_{i}
\,, \; i \in \nn \}
\end{equation}
of the direct product
\begin{equation} \label{eq_10}
\prod_{i \in \nn} \zz_{n_{i}} \;.
\end{equation}

In this notation the operation of addition in $A$ is defined component-wise, that is
$\vec{a} + \vec{b} = (a_{i} + b_{i})$ for any $\vec{a} = (a_{i})$,
$\vec{b} = (b_{i}) \in A$.
\end{rem}

It is known, that the topology of direct product~(\ref{eq_10}) is set
with the help of basis, which consists of so-called cylindrical sets
$$
U (x_{i_{1}}, \ldots, x_{i_{k}}) =
\{ (a_{i}) \; | \; a_{i_{s}} = x_{i_{s}} \, \; s = 1, \ldots, k \} \;; \quad
x_{i_{s}} \in \zz_{n_{i_{s}}} \,, \; i_{1} < \ldots < i_{k} \,, \; k \in \nn \;.
$$

From definition of $A$ (see relation~(\ref{eq_9})) it is easy
to see, that
$$
U (x_{i_{1}}, \ldots, x_{i_{k}}) \cap A = U (x_{i_{k}}) \cap A
$$
for any $k \in \nn $, $i_{1} < \ldots < i_{k} $ and $x_{i_{s}} \in
\zz_{n_{i_{s}}}$.
So, the family of sets
\begin{equation} \label{eq_11}
\begin{array}{rcl}
    V_{x_{j}} & = & U (x_{j}) \cap A = \{ (a_{i}) \in A \; | \;
        a_{j} = x_{j} \} \quad = \\
    & = & \{ (a_{i}) \in A \; | \; a_{j} = x_{j} \,, \; a_{k} = \varphi_{k} \circ
        \ldots \circ \varphi_{j-1} (x_{j}) \mbox{ for } k < j \} \;; \quad
        j \in \nn \,, \; x_{j} \in \zz_{n_{j}}
\end{array}
\end{equation}
appears to be basis of topology of space $A$.

The {\em natural metric} $\dist : A \times A \rightarrow \rr_{+} $ on $A$,
associated with the sequence $\{ n_{i} \}$, is defined by the following
correlation
$$
\dist (\vec{x}, \vec{y}) = \frac{1}{n_{m}} \;, \quad
m = \min \{ i \in \nn \; | \; x_{k} = y_{k} \mbox{ when } k < i \mbox{ and }
x_{i} \neq y_{i} \} \;.
$$
The correctness of this definition is verified immediately.

Consider an element $\vec{e} = (1) = (1, \ldots, 1, \ldots) \in A$. It
is called {\em generator} of the group $A$ and
the cyclical subgroup $\langle \vec{e} \rangle$, generated by
this element is dense in $A$ in the topology $\Tau$.

The translation map
$$
g: A \rightarrow A \;,
$$
$$
g: \vec{x} \mapsto \vec{x} + \vec{e}
$$
obviously is homeomorphism.

\begin{defn}
Dynamical system $(A, g)$ is called {\em odometer}.
\end{defn}

\begin{rem} \label{rem_3.2}
From that fact the subgroup $\langle \vec{e} \rangle$ is dense in $A$ it immediately
follows, that each trajectory of dynamical system $(A, g)$ is dense in $A$.
In other words odometer is always minimal dynamical system.
\end{rem}

\begin{lemma} \label{lemma_3.1}
For any $k \in \nn$ and $x_{k} \in \zz_{n_{k}}$ a family of sets
$\{ W^{(n_{k})}_{j} = V_{x_{k} +j} \}_{j = 0, \ldots, n_{i} -1}$ is
periodic partition of dynamical system $(A, g)$ of length $n_{k}$.
\end{lemma}

\begin{proof}
Obviously,
$$
A = \bigcup_{s \in \zz_{n_{k}}} V_{s} = \bigcup_{j \in \zz_{n_{k}}}
V_{x_{k} +j} \;.
$$
Hence, for the family $\{ W^{(n_{k})}_{j} \}$ condition~(iv) of
Definition~\ref{defn_1.1} is fulfilled.

Since all sets $V_{x_{k}+j}$, $j \in \zz_{n_{k}}$, are open
by definition and pairwise disjoint, the collection $\{ W^{(n_{k})}_{j} \}$
satisfies also to conditions~(i) and~(iii) of definition of periodic partition.

To complete the proof it remains to verify that $g (V_{a_{k}}) =
V_{a_{k} +1}$ ($ 1 \in \zz_{n_{k}}$) for every $a_{k} \in \zz_{n_{k}}$.

Let $ \vec{b} = (b_{i}) \in V_{a_{k}}$. Then $b_{k} = a_{k}$ and
$g (\vec{b}) = \vec{b} + \vec{e} = (b_{i}+1) \in V_{a_{k}+1}$.
Therefore, $g (V_{a_{k}}) \subseteq V_{a_{k} +1}$.

Conversely, let $\vec{c} = (c_{i}) \in V_{a_{k} +1}$. Then $c_{k} =
a_{k} +1$ and $g^{-1} (\vec{c}) = \vec{c} - \vec{e} = (c_{i} -1) \in
V_{a_{k}}$. Hence, $g (V_{a_{k}}) \supseteq V_{a_{k} +1}$.
\end{proof}

\begin{rem} \label{rem_3.3}
It  immediately follows from relation~(\ref{eq_11}), that
$$
\dist (\vec{x}, \vec{y}) = \dist (g (\vec{x}), g (\vec{y}))
$$
for all $\vec{x} $, $ \vec{y} \in A$.
\end{rem}

\subsection{Regular sequences of periodic partitions and
associated partitions on a phase space of dynamical
system.}\label{subsect_3.2}

Let $(X, f)$ be a dynamical system with compact phase space,
$\{ n_{i} \in \EssP{X, f} \}_{i \in \nn}$ is an unlimited regular
sequence. Let $\{ W^{(n_{i})} \}$ be a regular
sequence of periodic partitions of dynamical system $(X, f)$.

Let $x \in X$. Remark, that in the strength of properties of periodic partitions
for every $i \in \nn$ there exists unique $\alpha_{i} (x) \in
\zz_{n_{i}}$, such that $x \in W^{(n_{i})}_{\alpha_{i} (x)}$. In other words
the map is correctly defined
$$
F : X \rightarrow \prod_{i \in \nn} \zz_{n_{i}} \,,
$$
$$
F: x \mapsto \left( \alpha_{i} (x) \right) \,, \quad x \in X \,.
$$

We associate with every $x \in X$ a subset
$$
H (x) = \bigcap_{i \in \nn} W^{(n_{i})}_{\alpha_{i} (x)} \ni x
$$
of the space $X$. It follows from Definitions~\ref{defn_1.1},~\ref{defn_1.6} and
Corollary~\ref{corr_1.6}, that
\begin{itemize}
    \item[1)] all $H (x)$ are nonempty closed sets;
    \item[2)] $H (x) = H (y)$ if $F (x) = F (y)$ and $H (x) \cap H (y) =
        \emptyset$ if $F (x) \neq F (y)$;
    \item[3)] $F (f^{\pm 1} (x)) = F (x) \pm \vec{e}$ for all $x \in X$
        (remind, that $\vec{e} = (1, 1, \ldots, 1, \ldots) \in A$).
\end{itemize}

For every $\vec{a} = (a_{i}) \in F (X)$ we fix $x \in F^{-1} (\vec{a})$
and designate $H_{\vec{a}} = H (x)$. It follows from 2), that the set
$H_{\vec{a}}$ does not depend on a choice of $x \in F^{-1} (\vec{a})$.

From 1) and 2) it immediately follows, that the family of sets
$\hh = \{ H_{\vec{a}} \}_{\vec{a} \in F (X)} = \zer (F)$ is partition
of the space $X$,  elements of which are pre-images of points
of space $F (X)$. Also diagram is commutative
$$
\begin{CD}
@. X @>{F}>> F (X) \\
@. @V{pr}VV   @| \\
X / \zer (F) @= X / \hh @>{\fact F}>> F (X)
\end{CD}
$$

\begin{prop} \label{prop_3.1}
$F$ is the continuous map.
\end{prop}

\begin{proof}
Consider a subbasis of topology
$$
U_{x_{j}} = \{ \vec{a} = (a_{i}) \in \prod_{i \in \nn} \zz_{n_{i}}
\; | \; a_{j} = x_{j} \} \,, \qquad
x_{j} \in \zz_{n_{j}} \,, \quad j \in \nn
$$
of the space $\prod_{i \in \nn} \zz_{n_{i}}$.

The easy immediate verification shows, that
$$
F^{-1} (U_{x_{j}}) = W^{(n_{j})}_{x_{j}} \,,
\qquad x_{j} \in \zz_{n_{j}} \,, \quad j \in \nn \,.
$$
To complete the proof it suffices to recollect, that all sets
$W^{(n_{j})}_{x_{j}}$ are open in $X$ by definition.
\end{proof}

$X$ is the compact set, $\fact F$ is continuous bijective
map of $X / \hh$ on $F (X)$ and space $F (X)$ is Hausdorff, therefore
$\fact F$ is the homeomorphism (see~\cite{Engelking}).

For every $\vec{a} \in F (X)$ from 2) and 3) the equality is easily received
$f(F^{-1} (\vec{a})) = F^{-1} (\vec{a} + \vec{e})$. Thus, if we
designate
$$
g : F (X) \rightarrow F (X) \,, \quad g : \vec{a} \mapsto \vec{a} + \vec{e} \,;
$$
$$
\overline{f} = \fact {f} : X / \hh \rightarrow X / \hh \,, \quad
\overline{f} : H_{\vec{a}} \mapsto H_{\vec{a} + \vec{e}} \,;
$$
then we receive the commutative diagram
\begin{equation} \label{eq_12}
\begin{CD}
( X, f) @>{F}>> (F (X), g) \\
@V{pr}VV         @| \\
( X / \hh, \overline{f}) @>{\fact F}>> (F (X), g)
\end{CD}
\end{equation}

We now ask the question: what is the set $F (X)$?

We fix $x \in X$ and consider the set
$$
F(\Orb_{f} (x)) = \{ F(x) + n \vec{e} \; | \; n \in \zz \} =
F(x) + \langle \vec{e} \rangle \,.
$$
It is clear, that $F(x) + \langle \vec{e} \rangle \subseteq F(\Cl{\Orb_{f} (x)})
\subseteq \Cl{F(x) + \langle \vec{e} \rangle} = F(x) + \Cl{\langle \vec{e}
\rangle} = F(x) + A$ ($A$ is adic group constructed on
the sequence $\{ n_{i} \}$, see above).

Since $\Cl{\Orb_{f} (x)}$ is compact set, then $F(\Cl{\Orb_{f} (x)})$
is closed in $\prod_{i \in \nn} \zz_{n_{i}}$. The set $F(x) + \langle
\vec{e} \rangle$ is dense in $F(x) + A$, therefore $F(\Cl{\Orb_{f} (x)})
= F (x) + A$.

Let now $y$ is another point of the space $X$. $\Cl{\Orb_{f} (x)}$
is the closed invariant subset of dynamical system $(X, f)$.
Hence $W^{(n_{i})}_{s_{i}} \cap \Cl{\Orb_{f} (x)} \neq \emptyset$
for all $i \in \nn$, $s_{i} \in \zz_{n_{i}}$ (see Remark~\ref{rem_1.3}).

Let $F(y) = (\beta_{i})$. It follows from Remark~\ref{rem_1.11} that
$H_{(\beta_{i})} \cap \Cl{\Orb_{f} (x)} \neq \emptyset$ and
$F(y) = F(H_{(\beta_{i})}) \in F(\Cl{\Orb_{f} (x)}) = F (x) + A$.

As a result we receive
$$
F (X) = F (x) + A \,,
$$
In other words $F (X)$ is the coset of group $\prod_{i \in \nn} \zz_{n_{i}}$
on the subgroup $A$.

\begin{rem} \label{rem_3.4}
Obviously, $F(X) = A$ if and only if
\begin{equation} \label{eq_13}
\bigcap_{i \in \nn} W^{(n_{i})}_{0} \neq \emptyset \,.
\end{equation}
\end{rem}

\begin{defn} \label{defn_3.2}
A regular sequence $\{ W^{(n_{i})} \}_{i \in \nn}$ of
periodic partitions of dynamical system $(X, f)$ is called {\em
coherent}, if it satisfies to the relation~(\ref{eq_13}).
\end{defn}

From Proposition~\ref{prop_1.6}, Remark~\ref{rem_1.5} and construction,
given above, we get

\begin{prop} \label{prop_A}
Let $(X, f)$ is a dynamical system with Hausdorff compact phase
space.

For any unlimited regular sequence
$\{ n_{i} \in \EssP{X, f} \}_{i \in \nn}$ there exists a projection
$\pi : (X, f) \rightarrow (A, g)$ onto odometer $(A, g)$, constructed on
the sequence $\{ n_{i} \}_{i \in \nn}$.
\end{prop}

Let now $\{ n_{i} \in \EssP{X, f} \}_{i \in \nn}$ is a regular
sequence and let $\{ W^{(n_{i})} \}$ and
$\{ \widetilde{W}^{(n_{i})} \}$ are two compatible regular
sequences of periodic partitions of dynamical system $(X, f)$.
Then (see Remark~\ref{rem_1.5} and Definition~\ref{defn_1.3})
periodic partitions $W^{(n_{i})}$ and
$\widetilde{W}^{(n_{i})}$ are equivalent for every $i \in
\nn$. From this we immediately conclude, that
$$
\hh = \hh (\{ W^{(n_{i})} \}_{i \in \nn}) =
\hh (\{ \widetilde{W}^{(n_{i})} \}_{i \in \nn}) =
\widetilde{\hh} \,.
$$

\begin{prop} \label{prop_3.2}
Let $\{ W^{(n_{i})} \}_{i \in \nn}$ is a regular sequence of
periodic partitions of a dynamical system $(X, f)$.

The family of sets $\{ pr (W^{(n_{i})}_{s_{i}}) \; | \; i \in \nn \, \; s_{i} \in
\zz_{n_{i}} \}$ satisfies to the following properties:
\begin{itemize}
    \item[1)] it is the regular sequence of periodic partitions of
        the dynamical system $(X / \hh, \overline{f})$;
    \item[2)] it is the basis of topology on space $X / \hh$.
\end{itemize}
\end{prop}

\begin{proof}
It immediately follows from what has been said that the sequence
$\{W^{(n_{i})} \}$ could be considered as coherent.

Now proposition follows from the relations
$$
W^{(n_{i})}_{s_{i}} = F^{-1} (V_{s_{i}}) \,, \qquad s_{i} \in \zz_{n_{i}}
\,, \quad i \in \nn
$$
(see formula~(\ref{eq_11}), Lemma~\ref{lemma_3.1} and
Corollary~\ref{corr_1.6}) and from commutative diagram~(\ref{eq_12}), the lower
arrow in which is homeomorphism.
\end{proof}

\begin{prop} \label{prop_3.3}
Let $\{ n_{i} \in \EssP{X, f} \}_{i \in \nn}$ and
$\{ m_{j} \in \EssP{X, f} \}_{j \in \nn}$ are two regular
sequences, such that
$\Phi (\{ n_{i} \; | \; i \in \nn \}) \leq
\Phi (\{ m_{j} \; | \; j \in \nn \})$.

Let regular sequences
$\{ W^{(n_{i})} \}_{i \in \nn}$ and
$\{ \widetilde{W}^{(m_{j})} \}_{j \in \nn}$ of periodic partitions of
dynamical system $(X, f)$ are compatible.

Then the partition $\widetilde{\hh}$ of space $X$, which is induced by
sequence $\{ \widetilde{W}^{(m_{j})} \}$, is the refinement of
the partition $\hh$, induced by the sequence
$\{ W^{(n_{i})} \}$.
\end{prop}

\begin{proof}
We fix $x \in X$. There exist such $(\alpha_{i}) \in \prod_{i \in \nn}
\zz_{n_{i}}$ and $(\beta_{j}) \in \prod_{j \in \nn} \zz_{m_{j}}$, that
$$
x \in
\bigcap_{i \in \nn} W^{(n_{i})}_{\alpha_{i}}
\cap
\bigcap_{j \in \nn} \widetilde{W}^{(m_{j})}_{\beta_{j}} \,.
$$
According to Proposition~\ref{prop_2.3}, for every $i \in \nn$ there exists
$k(i) \in \nn$, such that $n_{i}$ divides $m_{k (i)}$. Since periodic
partitions $W^{(n_{i})}$ and $\widetilde{W}^{(m_{k(i)})}$
are compatible, then by Corollary~\ref{corr_1.6} second of them is
the refinement of first one.

Thus, $\widetilde{W}^{(m_{k (i)})}_{\beta_{k (i)}} \subseteq
W^{(n_{i})}_{\alpha_{i}}$ for every $i \in \nn$.
Therefore, we have
$$
H_{(\alpha_{i})} =
\bigcap_{i \in \nn} W^{(n_{i})}_{\alpha_{i}}
\supseteq
\bigcap_{i \in \nn} \widetilde{W}^{(m_{k(i)})}_{\beta_{k(i)}}
\supseteq
\bigcap_{j \in \nn} \widetilde{W}^{(m_{j})}_{\beta_{j}} =
\widetilde{H}_{(\beta_{j})} \,.
$$

By virtue of arbitrariness in a choice of $x \in X$ we conclude from what was said above, that for
an arbitrary $(\alpha_{i}) \in F(X)$ and
$(\beta_{j}) \in \widetilde{F}(X)$ either $H_{(\alpha_{i})}
\cap \widetilde{H}_{(\beta_{j})} = \emptyset$ or
$H_{(\alpha_{i})} \supseteq \widetilde{H}_{(\beta_{j})}$.
\end{proof}

\begin{corr} \label{corr_3.1}
If in the conditions of Proposition~\ref{prop_3.2} the equality takes place
$\Phi (\{ n_{i} \; | \; i \in \nn \}) = \Phi (\{ m_{j} \; | \; j \in \nn \})$,
then partitions $\hh$ and $\widetilde{\hh}$ of the space $X$
coincide.
\end{corr}

\begin{prop} \label{prop_3.4}
Let $\{ n_{i} \in \EssP{X, f} \}_{i \in \nn}$ and
$\{ m_{j} \in \EssP{X, f} \}_{j \in \nn}$ are two regular
sequences.

Let $\{ W^{(n_{i})} \}_{i \in \nn}$ and
$\{ \widetilde{W}^{(m_{j})} \}_{j \in \nn}$ are regular sequences
of periodic partitions of dynamical system $(X, f)$. Assume that $\hh$ and
$\widetilde{\hh}$ are partitions of the space $X$, induced by these
sequences.

Let sequences $\{ W^{(n_{i})} \}$ and
$\{ \widetilde{W}^{(m_{j})} \}$ are not compatible.

Then
$H(x) \setminus \widetilde{H}(x) \neq \emptyset$
and
$\widetilde{H}(x) \setminus H(x) \neq \emptyset$
for every $x \in X$.
\end{prop}

\begin{proof}
Let $x \in X$. There exist such $(\alpha_{i}) \in \prod_{i \in
\nn} \zz_{n_{i}}$ and $(\beta_{j}) \in \prod_{j \in \nn} \zz_{m_{j}}$, that
$$
H(x) = H_{(\alpha_{i})} =
\bigcap_{i \in \nn} W^{(n_{i})}_{\alpha_{i}}
\,, \quad
\widetilde{H}(x) = \widetilde{H}_{(\beta_{j})} =
\bigcap_{j \in \nn} \widetilde{W}^{(m_{j})}_{\beta_{j}} \,.
$$
On definition $H (x) \cap \widetilde{H} (x) \neq \emptyset$.

According to conditions of Proposition there exist such $k$, $l \in \nn$, that
periodic partitions $W^{(n_{k})}$ and $\widetilde{W}^{(m_{l})}$
are not compatible. That is
$$
x \in A_{k, l} = \bigcup_{t \in \zz}
f^{t}
\Bigl(
W^{(n_{k})}_{\alpha_{k}} \cap
\widetilde{W}^{(m_{l})}_{\beta_{l}}
\Bigr)
\neq X \,.
$$

So, the space $X$ falls into the union of two disjoint
closed invariant sets $A_{k, l}$ and $B_{k, l} = X \setminus
A_{k, l}$ of dynamical system $(X, f)$ (see Proposition~\ref{prop_1.4}).

Obviously, $H (x) \cap \widetilde{H}(x) \subseteq
W^{(n_{k})}_{\alpha_{k}} \cap
\widetilde{W}^{(m_{l})}_{\beta_{l}} \subseteq A_{k, l}$.

However,
$$
H(x) \setminus \widetilde{H}(x) \supseteq H(x) \cap
B_{k, l} = \bigcap_{i \in \nn}
\Bigl(
W^{(n_{i})}_{\alpha_{i}} \cap B_{k, l}
\Bigr)
\neq \emptyset
$$
(see Remarks~\ref{rem_1.3} and~\ref{rem_1.11}). Similarly,
$\widetilde{H} (x) \setminus H (x) \neq \emptyset$.
\end{proof}

\begin{prop} \label{prop_3.5}
Let $\{ n_{i} \in \EssP{X, f} \}_{i \in \nn}$ and
$\{ m_{j} \in \EssP{X, f} \}_{j \in \nn}$ are two regular
sequences.

Let $\{ W^{(n_{i})} \}_{i \in \nn}$ and
$\{ \widetilde{W}^{(m_{j})} \}_{j \in \nn}$ are regular sequences of
periodic partitions of dynamical system $(X, f)$. Assume $\hh$ and
$\widetilde{\hh}$ are the partitions of space $X$, induced by these
sequences.

If there exists $x \in X$, such that
$$
H (x) = \bigcap_{i \in \nn} W^{(n_{i})}_{\alpha_{i} (x)} \supseteq
\bigcap_{j \in \nn} \widetilde{W}^{(m_{j})}_{\widetilde{\alpha}_{j}(x)} =
\widetilde{H} (x) \;,
$$
then the partition $\widetilde{\hh}$ is refinement of the partition $\hh$,
sequences $\{ W^{(n_{i})} \}$ and $\{ \widetilde{W}^{(m_{j})} \}$
are compatible and $\Phi (\{ n_{i} \; | \; i \in \nn \}) \leq
\Phi (\{ m_{j} \; | \; j \in \nn \})$.
\end{prop}

\begin{corr} \label{corr_3.2}
Let $\{ W^{(n_{i})} \}_{i \in \nn}$ and
$\{ \widetilde{W}^{(m_{j})} \}_{j \in \nn}$ are regular sequences of
periodic partitions of dynamical system $(X, f)$. Assume $\hh$ and
$\widetilde{\hh}$ are the partitions of space $X$, induced by these
sequences.

The following statements are equivalent:
\begin{itemize}
    \item[1)] the partition $\widetilde{\hh}$ is refinement of the partition
        $\hh$ (respectively, $\hh = \widetilde{\hh}$);
    \item[2)] sequences $\{ W^{(n_{i})} \}$ and
        $\{ \widetilde{W}^{(m_{j})} \}$ are compatible and $\Phi (\{ n_{i} \; | \;
        i \in \nn \}) \leq \Phi (\{ m_{j} \; | \; j \in \nn \})$
        (respectively, $\Phi (\{ n_{i} \; | \; i \in \nn \}) = \Phi (\{ m_{j}
        \; | \; j \in \nn \})$).
\end{itemize}
\end{corr}

In order to prove Proposition~\ref{prop_3.5} we will need following
almost obvious

\begin{lemma} \label{lemma_3.2}
Let $X$ be a Hausdorff space,
$$
K_{1} \supseteq K_{2} \supseteq \ldots \supseteq K_{i} \supseteq \ldots
$$
is a sequence of nonempty compact subsets of $X$.

For any open neighborhood $U$ of the set
$$
K = \bigcap_{i \in \nn} K_{i}
$$
there exists $n \in \nn$, such that $K_{i} \subseteq U$ for $i \geq n$.
\end{lemma}

\begin{proof}
Assume that there exist a neighborhood $U \supseteq K$ and a sequence
$\{ x_{i} \in K_{i} \setminus U \}_{i \in \nn}$.

Since $x_{i} \in K_{1}$, $i \in \nn$, on the construction and
$K_{1}$ is the compact set, then this sequence has at least one
limit point $x \in K_{1} \setminus U \subseteq X \setminus U$.

It follows from condition of Lemma that $x_{i} \in K_{i} \subseteq K_{m}$  for every
$m \in \nn$ and $i \geq m$. Hence, $x \in K_{m}$, $m \in \nn$ and $x \in K
\setminus U$.

The obtained contradiction proves Lemma.
\end{proof}

\begin{proof}[of Proposition~\ref{prop_3.5}]
{\bf 1.}
We shall prove, that for every $i \in \nn $ there exists $j (i) \in \nn$,
which comply with the following requirements:
\begin{itemize}
    \item[---] $n_{i}$ divides $m_{j (i)}$;
    \item[---] periodic partition $\widetilde{W}^{(m_{j (i)})}$
        is the refinement of partition $W^{(n_{i})}$.
\end{itemize}

We fix $i \in \nn$. Obviously, the open neighborhood
$W^{(n_{i})}_{\alpha_{i} (x)}$ of the set $ \widetilde{H} (x) \subseteq
H(x) \subseteq W^{(n_{i})}_{\alpha_{i}(x)}$ and the sequence of closed
sets
$$
\widetilde{W}^{(m_{1})}_{\widetilde{\alpha}_{1}(x)} \supseteq
\widetilde{W}^{(m_{2})}_{\widetilde{\alpha}_{2}(x)} \supseteq
\ldots \supseteq
\widetilde{W}^{(m_{j})}_{\widetilde{\alpha}_{j}(x)} \supseteq
\ldots
$$
satisfy to condition of Lemma~\ref{lemma_3.2}.
Hence, there exists $j (i) \in \nn$, for which
$\widetilde{W}^{(m_{j(i)})}_{\widetilde{\alpha}_{j(i)}(x)} \subseteq
W^{(n_{i})}_{\alpha_{i} (x)}$.

From Corollary~\ref{corr_1.5} we conclude, that
periodic partitions $\widetilde{W}^{(m_{j (i)})}$ and $W^{(n_{i})}$
are compatible and $n_{i}$ divides $m_{j (i)}$. Now it follows from Corollary~\ref{corr_1.6}
that the partition $\widetilde{W}^{(m_{j (i)})}$ is refinement of the
partition $W^{(n_{i})}$.

{\bf 2.}
Verify, that $H (y) \supseteq \widetilde{H}(y)$ for every $y \in X$,
in other words the partition $\widetilde{\hh}$ is refinement of the partition $\hh$.

Really, it follows from what we said above, that
$$
\widetilde{H}(y) =
\bigcap_{j \in \nn} \widetilde{W}^{(m_{j})}_{\widetilde{\alpha}_{j}(y)}
\subseteq
\bigcap_{i \in \nn} \widetilde{W}^{(m_{j(i)})}_{\widetilde{\alpha}_{j(i)}(y)}
\subseteq
\bigcap_{i \in \nn} W^{(n_{i})}_{\alpha_{i} (y)} = H(y)
$$
for every $y \in X$.

{\bf 3.}
The previous item and Proposition~\ref{prop_3.4} implies, that
the sequences of periodic partitions $\{ W^{(n_{i})} \}_{i \in \nn}$
and $\{ \widetilde{W}^{(m_{j})} \}_{j \in \nn}$ are compatible.

{\bf 4.}
We conclude from item 1 and Proposition~\ref{prop_2.3}, that
$\Phi (\{ n_{i} \; | \; i \in \nn \}) \leq
\Phi (\{ m_{j} \; | \; j \in \nn \})$.
\end{proof}

\subsection{Main properties of odometers} \label{subsect_3.3}

We show now on the example of odometers how the proved above statements
could be applied.

\begin{prop} \label{prop_3.6}
Let $(X, f)$ and $(Y, g)$ are dynamical systems, $p : (X, f) \rightarrow (Y,
g)$ is a projection. If $n \in \EssP{Y, h}$ and $W^{(n)} =
\{ W^{(n)}_{i} \}_{i \in \zz_{n}}$ is a periodic partition of dynamical
system $(Y, g)$, then $n \in \EssP{X, f}$ and $ \widetilde{W}^{(n)} =
\{\widetilde{W}^{(n)}_{i} = p^{-1} (W^{(n)}_{i}) \}_{i \in \zz_{n}}$ is the
periodic partition of dynamical system $(X, f)$.
\end{prop}

\begin{proof}
is the simple immediate verification.
\end{proof}

\begin{corr} \label{corr_3.3}
Let $(Y, h)$ is a factor--system of a dynamical system $(X, f)$. Then
$\EssP{Y, h} \subseteq \EssP{X, f}$.
\end{corr}

We take an advantage of Remark~\ref{rem_2.3} and obtain

\begin{corr} \label{corr_3.4}
In conditions of Corollary~\ref{corr_3.3} the inequality $\Phi (\EssP{Y,
h}) \leq \Phi (\EssP{X, f})$ is fulfilled.
\end{corr}

\begin{rem} \label{rem_3.5}
Thus, if dynamical systems $(X, f)$ and $(Y, g)$ are topologically
conjugate, then $\Phi (\EssP{X, f}) = \Phi (\EssP{Y, g})$.
Hence, $\Phi (\EssP{X, f}) \in \Sigma$ is
topological invariant of the dynamical system $(X, f) \in
\kk_{0}$.
\end{rem}

\begin{prop} \label{prop_3.7}
Let $(A, g)$ is an odometer constructed on a regular sequence $ \{
n_{i} \}_{i \in \nn}$.

Then $\Phi (\{\EssP{A, g} \}) = \Phi (\{ n_{i} \; | \; i \in \nn \})$.

Let $\{ W^{(m_{j})} \}$ is a regular sequence of periodic
partitions of dynamical system $(A, g)$.

The family of sets $\{ W^{(m_{j})}_{r_{j}} \; | \; r_{j} \in \zz_{m_{j}} \,, \; j
\in \nn \}$ is basis of topology of the space $A$ if and only if
$\Phi (\{ m_{j} \; | \; j \in \nn \}) = \Phi (\{\EssP{A, g} \})$.
\end{prop}

\begin{proof}
It follows from Lemma~\ref{lemma_3.1} and relation~(\ref{eq_11}) that the family
$$
W^{(n_{i})} = \{ W^{(n_{i})}_{s_{i}} = V_{s_{i}} \}_{s_{i} \in
\zz_{n_{i}}} \,, \quad i \in \nn \,,
$$
is the regular sequence of periodic partitions of
dynamical system $(A, g)$. We construct on this sequence the
associated partition $\hh$ of space $A$. Since the family~(\ref{eq_11}) is
basis of the topology of space $A$, then $H_{\vec{a}} = \{\vec{a} \}$ for
every $\vec{a} \in A$.

{\bf 1.} {}
The set $\EssP{A, g}$ is admissible (see Remark~\ref{rem_2.6}),
therefore there exists regular sequence
$\{ m_{j} \in \EssP{A, g} \}_{j \in \nn}$, such that
$$
\Phi (\{ m_{j} \; | \; j \in \nn \}) = \Phi (\EssP{A, g}) \geq
\Phi (\{ n_{i} \; | \; i \in \nn \})
$$
(see Proposition~\ref{prop_2.2} and Remark~\ref{rem_2.8}).

Construct on this sequence regular sequence of
periodic partitions $\{\widetilde{W}^{(m_{j})} \}_{j \in \nn}$, which is
compatible with the sequence $\{ W^{(n_{i})} \}_{i \in \nn}$
(see Proposition~\ref{prop_1.8}).

Let $\widetilde{\hh}$ is the partition of space $A$, induced
by the sequence $\{\widetilde{W}^{(m_{j})} \}$. Then it follows from
Proposition~\ref{prop_3.3} that the partition $\widetilde{\hh}$
is refinement of the partition $\hh$. And it is possible only if
$\widetilde{\hh} = \hh$. Now from Proposition~\ref{prop_3.5} we receive
$$
\Phi (\{ m_{j} \; | \; j \in \nn \}) =
\Phi (\{ n_{i} \; | \; i \in \nn \}) \,.
$$

{\bf 2.} {}
Let now $\{\widetilde{W}^{(m_{j})} \}_{j \in \nn}$ is a certain
regular sequence of periodic partitions of the dynamical system
$(A, g)$. Assume $\widetilde{\hh}$ is the partition of space $A$, induced
by the sequence $\{\widetilde{W}^{(m_{j})} \}$.

We conclude from the first part of Proposition~\ref{prop_3.7} and
Remark~\ref{rem_2.8} that $\Phi (\{ n_{i} \; | \; i \in \nn \}) \geq
\Phi (\{ m_{j} \; | \; j \in \nn \})$.

Odometer $(A, g)$ is the minimal dynamical system (see.
Remark~\ref{rem_3.2}), so this dynamical system is indecomposable and
regular sequences $\{ W^{(n_{i})} \}$ and
$\{\widetilde{W}^{(m_{j})} \}$ are compatible (see Remark~\ref{rem_1.12}).

Applying Corollary~\ref{prop_3.2} we conclude, that
the partition $\hh$ of space $A$ is refinement of the partition
$\widetilde{\hh}$, with $\hh = \widetilde{\hh}$ if and
only if
$$
\Phi (\{ m_{j} \; | \; j \in \nn \}) =
\Phi (\{ n_{i} \; | \; i \in \nn \}) = \Phi (\EssP{A, g}) \,.
$$

Remind, that a family $\{ U_{\alpha} \}_{\alpha \in \lambda}$ of open
subsets of topological space $X$ is its basis of topology,
when the following conditions (see~\cite{Kelley}) are fulfilled:
\begin{itemize}
    \item[(a)] if $U_{\alpha} \cap U_{\beta} \neq \emptyset$ for certain
        $\alpha$, $\beta \in \Lambda$, then there exists $\gamma \in \Lambda$,
        such that $U_{\gamma} \subseteq U_{\alpha} \cap U_{\beta}$;
    \item[(b)] for every $x \in X$ and any open neighborhood $U$ of
        $x$ there exists $\alpha \in \Lambda$, such that $x \in U_{\alpha}
        \in U$.
\end{itemize}

Remark, that for any regular sequence of periodic
partitions condition (a) is always fulfilled (it follows immediately from
Definition~\ref{defn_1.6} and Remark~\ref{rem_1.10}).

Let $\Phi (\{ m_{j} \; | \; j \in \nn \}) \lneqq
\Phi (\{ n_{i} \; | \; i \in \nn \})$. Then partitions $\hh$ and $\widetilde{\hh}$
do not coincide (see Corollary~\ref{corr_3.2}) and there exist two points
$x_{1}$, $x_{2} \in A$, $x_{1} \neq x_{2}$, contained in the same element of
the partition $\widetilde{\hh}$. Hence, for every
$\widetilde{W}^{(m_{j})}_{s_{j}}$, $s_{j} \in \zz_{m_{j}}$, $j \in \nn$,
either $\widetilde{W}^{(m_{j})}_{s_{j}} \cap \{ x_{1}, x_{2} \} = \emptyset$
or $\{ x_{1}, x_{2} \} \subseteq \widetilde{W}^{(m_{j})}_{s_{j}}$, and
the condition (b) is not fulfilled.

Let now $\Phi (\{ m_{j} \; | \; j \in \nn \}) =
\Phi (\{ n_{i} \; | \; i \in \nn \})$. In this case partitions $\hh$ and
$\widetilde{\hh}$ coincide.

Let $x \in X$ and $U$ is an open neighborhood of $x$. On definition of
the partition $\widetilde{\hh}$ there exists the unique sequence
$\{\alpha_{j} \in \zz_{m_{j}} \}_{j \in \nn}$, such that
$$
\{ x \} = H (x) = \widetilde{H} (x) =
\bigcap_{i \in \nn} \widetilde{W}^{(m_{j})}_{\alpha_{j}} \,.
$$
In this case all sets from intersection in the right-hand part of the equality
are compact (being closed subsets of the compact space $A$) and
$\widetilde{W}^{(m_{j+1})}_{\alpha_{j+1}} \subseteq
\widetilde{W}^{(m_{j})}_{\alpha_{j}}$, $j \in \nn$
(see Definition~\ref{defn_1.6}). We apply Lemma~\ref{lemma_3.2}
and conclude that there exists $k \in \nn$, for which $x \in
\widetilde{W}^{(m_{j_{k}})}_{\alpha_{j_{k}}} \subseteq U$.
\end{proof}

From Proposition~\ref{prop_A}, Corollary~\ref{corr_3.4},
Proposition~\ref{prop_3.7}, Remark~\ref{rem_2.6} and
Proposition~\ref{prop_2.2} we receive the following statement.

\begin{theorem} \label{theorem_A}
Let $(X, f)$ is a dynamical system with Hausdorff compact phase
space, $(A, g)$ is an odometer.

The following statements are equivalent:
\begin{itemize}
    \item[(i)] there exists a projection $\pi: (X, f) \rightarrow (A, g)$;
    \item[(ii)] the inequality $\Phi (\EssP{A, g}) \leq
        \Phi (\EssP{X, f})$ is fulfilled.
\end{itemize}
\end{theorem}

In order to formulate the following statement, we require two
definitions.

Let $(X, f)$ is a dynamical system with compact metric phase
space $(X, \rho)$.

\begin{defn}
Points $x$, $y \in X$, $x \neq y$, are {\em distal}, if
there exists $\delta > 0$, such that $\rho (f^{n} (x), f^{n} (y)) > \delta$ for
every $n \in \zz$.

The dynamical system $(X, f)$ is called {\em distal}, if any pair of
points $x$, $y \in X$, $x \neq y$, is distal.
\end{defn}

\begin{defn}
The dynamical system $(X, f)$ is called {\em equicontinuous},
if the family of maps $\{ f^{n} \}_{n \in \zz}$ is equicontinuous
under the metrics $\rho$, that is if for every
$\varepsilon > 0$ there exists $\delta > 0$, such that if $\rho (x, y) <
\delta$ for some $x$, $y \in X$, then $\rho (f^{n} (x), f^{n} (y)) <
\varepsilon$ for every $n \in \zz$.
\end{defn}

\begin{rem}
It is easy to see, that distality and equicontinuity of dynamical
system $(X, f)$ do not depend (by virtue of compactness of $X$) on
a choice of metric function,
which generates the given topology on $X$. In other words, distality and
equicontinuity are topological properties of dynamical
system $(X, f)$ with a metrizable compact phase space $X$.
\end{rem}

\begin{theorem}[see.~\cite{Skau}]\label{theorem_3.1}
Let $(\Gamma, f)$ is a minimal dynamical system on Cantor set
$\Gamma$.

Then the following conditions are equivalent:
\begin{itemize}
    \item[1]. d. s. $(\Gamma, f)$ is topologically conjugate with an odometer;
    \item[2]. d. s. $(\Gamma, f)$ is distal;
    \item[3]. d. s. $(\Gamma, f)$ is equicontinuous.
\end{itemize}
\end{theorem}

\begin{proof}
Equivalence of conditions 2. and 3. for dynamical systems with zero-dimensional
compact phase space immediately follows from results,
obtained in paper~\cite{Ellis}.

Verify implication 1. $ \Rightarrow $ 3.

Let d. s. $(\Gamma, f)$ is topologically conjugate with the help of a homeomorphism
$h: \Gamma \rightarrow A$ with an odometer $(A, g)$, which is generated
by admissible sequence $\{ n_{i} \}_{i \in \nn}$. We transfer
the natural metric $\dist$ from the space $A$ on $\Gamma$ with the help of the
relations
$$
\rho (x, y) = \dist (h (x), h (y)) \,, \quad x, y \in \Gamma \,.
$$
Remark~\ref{rem_3.3} consequences, that the map $f$ is
isometric under the metric $\rho$. Thus, d. s. $(\Gamma,
f)$ is equicontinuous.

Prove implication 3. $\Rightarrow$ 1.

We fix the metrics $\rho: \Gamma \times \Gamma \rightarrow \rr_{+}$. Let
the dynamical system $(\Gamma, f)$ is equicontinuous under
the metrics $\rho$.

Let $x \in \Gamma$ and $V$ is an open-closed neighborhood of $x$. Since
the closed sets $V$ and $\Gamma \setminus V$ are disjoint and $\Gamma$
is compact, then
$$
\rho (V, \Gamma \setminus V) = \varepsilon > 0 \,.
$$
There exists $\delta > 0$, such that for any $y_{1}$, $y_{2} \in \Gamma$ the implication is valid
$$
(\rho (y_{1}, y_{2}) < 2 \delta) \Rightarrow
(\rho (f^{n} (y_{1}), f^{n} (y_{2})) < \varepsilon \,, \quad n \in \zz) \,.
$$

Let $U = U_{\delta} (x)$. Then $\Diam U < 2 \delta$ and $\Diam f^{n} (U) <
\varepsilon$ for all $n \in \zz$. Hence, for every $n \in \zz$
either $V \cap f^{n} (U) = \emptyset$ or $f^{n} (U) \subseteq V$.

Since the dynamical system $(\Gamma, f)$ is minimal, then there exists $k \in
\nn$, such that $f^{k} (x) \in U$. Then $f^{k} (U) \cap U \neq \emptyset$.
It is clear that since $f$ is the homeomorphism, we have
$$
f^{k+n} (U) \cap f^{n} (U) \neq \emptyset \,, \quad n \in \zz \,.
$$

Check, that $f^{kn} (U) \subseteq V$, $n \in \zz$.

We shall carry out verification on an induction for negative $n$.

{\em Basis of induction.} {}
Since $\emptyset \neq f^{-k} (U) \cap U \subseteq f^{-k} (U) \cap V$,
then $f^{-k} (U) \subseteq V$.

{\em Step of induction.} {}
Let $f^{-ki} (U) \subseteq V$ for some $i \in \nn$. Then
$\emptyset \neq f^{-k (i+1)} (U) \cap f^{-ki} (U) \subseteq f^{-k (i+1)} (U)
\cap V$ and $f^{-k (i+1)} (U) \subseteq V$.

On the induction we conclude, that $f^{kn} (U) \subseteq V$ for all $n < 0$.

The proof of this inclusion for positive $n$ is done similarly.

So,
$$
\bigcup_{n \in \zz} f^{kn} (x) \subseteq \bigcup_{n \in \zz} f^{kn} (U)
\subseteq V \,.
$$
We conclude from Lemma~\ref{lemma_1.5} that there exist
$m \in \EssP{\Gamma, f}$ and periodic partition $W^{(m)}$
of dynamical system $(\Gamma, f)$ of length $m$, such that
$x \in W^{(m)}_{0} \subseteq V$.

We fix $x \in \Gamma$ and sequence $\{\varepsilon_{i} \geq 0
\}_{i \in \nn}$, such that $\varepsilon_{i} \rightarrow 0$ when $i
\rightarrow \infty$. For every $\varepsilon_{i}$ we can find
$\delta_{i} > 0$, such that for all $y_{1}$, $y_{2} \in \Gamma$
$$
(\rho (y_{1}, y_{2}) < \delta_{i}) \Rightarrow
(\rho (f^{n} (y_{1}), f^{n} (y_{2})) < \varepsilon_{i} \,,
\quad n \in \zz) \,.
$$

Now we shall construct on an induction a coherent regular sequence
of periodic partitions $\{ W^{(n_{i})} \}_{i \in \nn}$ of dynamical
system $(\Gamma, f)$, such that for every $i \in \nn$
\begin{equation} \label{eq_14}
\Diam W^{(n_{i})}_{s_{i}} < \varepsilon_{i} \,, \quad
s_{i} \in \zz_{n_{i}} \,.
\end{equation}

The space $\Gamma$ is zero-dimensional, therefore there exists a sequence
$\{ V_{i} \ni x \}_{i \in \nn}$ of open-closed sets, such that
$\Diam V_{i} < \delta_{i}$, $i \in \nn$.

{\em Basis of induction.} {}
Find periodic partition $W^{(n_{1})}$ of dynamical system
$(\Gamma, f)$, such that $x \in W^{(n_{1})}_{0} \subseteq V_{1}$. Then
$\Diam W^{(n_{1})}_{0} < \delta_{1}$ and $\Diam W^{(n_{1})}_{s_{1}} = \Diam
f^{s_{1}} (W^{(n_{1})}_{0}) < \varepsilon_{1}$, $s_{1} \in \zz_{n_{1}}$,
according to the choice of $\delta_{1}$.

{\em Step of induction.} {}
Assume the family $\{ W^{(n_{i})} \}_{i = 1}^{k}$ of
periodic partitions of dynamical system $(\Gamma, f)$ is already constructed,
such that $n_{i} $ divides $n_{i+1} $, $i = 1, \ldots, k-1 $,
$$
W^{(n_{1})}_{0} \supseteq \ldots \supseteq W^{(n_{k})}_{0} \ni x
$$
(we conclude from Corollary~\ref{corr_1.5} that every two periodic
partitions from this family are compatible), and which complies with the
relations~(\ref{eq_14}).

Designate $\widetilde{V}_{k+1} = V_{k+1} \cap W^{(n_{k})}_{0} \ni x$.
Obviously, $\Diam \widetilde{V}_{k+1} < \delta_{k+1}$.

Find periodic partition $W^{(n_{k+1})}$ of the dynamical system
$(\Gamma, f)$, such that $x \in W^{(n_{k+1})}_{0} \subseteq
\widetilde{V}_{k+1} \subseteq W^{(n_{k})}_{0}$.

On one hand $\Diam W^{(n_{k+1})}_{0} \leq \Diam \widetilde{V}_{k+1} <
\delta_{k+1}$, hence $\Diam W^{(n_{k+1})}_{s_{1}} = \Diam
f^{s_{1}} (W^{(n_{k+1})}_{0}) < \varepsilon_{k+1} $, $s_{k+1} \in
\zz_{n_{k+1}}$.

On the other hand, $n_{k}$ divides $n_{k+1}$ and periodic partitions
$W^{(n_{k+1})}$ and $W^{(n_{k})}$ are compatible by Corollary~\ref{corr_1.5}.

On an induction we receive the coherent sequence of periodic partitions
$\{ W^{(n_{i})} \}_{i \in \nn}$ of dynamical system $(\Gamma, f)$, all
elements of which satisfy to relation~(\ref{eq_14}).

Construct now the partition $\hh$ of space $\Gamma$ on the sequence
$\{ W^{(n_{i})} \}$ and the projection $F: (\Gamma, f) \rightarrow (A, g)$.

Let $y \in \Gamma$, $(\alpha_{i}) = F(y) \in A$. Remark that for
every $k \in \nn$ we have
$$
\Diam H (y) \leq
\Diam \Bigl( \bigcap_{i=1}^{k} W^{(n_{i})}_{\alpha_{i}} \Bigr) =
\Diam W^{(n_{k})}_{\alpha_{k}} \leq \varepsilon_{k}
\underset{k \rightarrow \infty} {\rightarrow} 0 \,.
$$

Hence, $H (y) = \{ y \} $ for every $y \in \Gamma$ and the projection map
$pr : \Gamma \rightarrow \Gamma / \hh$ is bijective.
Therefore the map $F = (\fact F) \circ pr : \Gamma \rightarrow A$ is also
bijective. Since $\Gamma$ is compact set, then $F$ is the homeomorphism, which
conjugates dynamical systems $(\Gamma, f)$ and $(A, g)$ (see commutative
diagram~(\ref{eq_12})).
\end{proof}

\begin{rem} \label{rem_3.7}
When we defined odometer constructed on a regular sequence
$\{ n_{i} \}_{i \in \nn}$, we required, that this sequence should be
{\em unlimited}. Actually this requirement can be written in the following
way
$$
\Phi (\{ n_{i} \; | \; i \in \nn \}) \in \Sigma \setminus \Phi_{0} (\nn) \,.
$$

Let us look, what will change, if $\Phi (\{ n_{i} \; | \; i \in \nn \}) =
\Phi_{0} (m) \in \Phi_{0} (\nn)$. In this case $A = \LimInv_{i \rightarrow
\infty} \zz_{n_{i}} = \zz_{m} = \{ 0, 1, \ldots, m-1 \}$, $\vec{e} = 1 \in
\zz_{m}$, and the dynamical system $(A, g)$ consists from unique
periodic trajectory of length $m$.

We expand definition of odometers, including in it the case $\Phi (\{ n_{i}
\; | \; i \in \nn \}) \in \Phi_{0} (\nn)$.

It is trivially checked, that everything we said in subsections~\ref{subsect_3.2}
and~\ref{subsect_3.3}, except for theorem~\ref{theorem_3.1}, is valid for
our new definition.
\end{rem}

\begin{theorem}[see~\cite{Glimm, Down_1, Block}] \label{theorem_3.2}
{\bf 1.} {}
For every $N \in \Sigma$ there exists an odometer $(A, g)$, such that
$\Phi (\EssP{A, g}) = N$.

{\bf 2.} {}
Odometers $(A_{1}, g_{1})$ and $(A_{2}, g_{2})$ are topologically conjugate
if and only if $\Phi (\EssP{A_{1}, g_{1}}) = \Phi (\EssP{A_{2},
g_{2}})$.
\end{theorem}

\begin{proof}
{\bf 1.} {}
Let $N \in \Sigma$. We can find admissible set $\overline{A} \subseteq
\aa$, such that $\Phi (\overline{A}) = N$ (see Proposition~\ref{prop_2.1}).
Fix also regular sequence $\{ n_{i} \in \overline{A} \}_{i \in
\nn}$, for which $\Phi (\{ n_{i} \; | \; i \in \nn \}) =
\Phi (\overline{A})$ (see Proposition~\ref{prop_2.2}). Construct on
the sequence $\{ n_{i} \}$ odometer $(A, g)$. From
Proposition~\ref{prop_3.7} we conclude, that $\Phi (\EssP{A, g}) =
\Phi (\overline{A}) = N$.

{\bf 2.} {} (a)
Assume that odometers $(A_{1}, g_{1})$ and $(A_{2}, g_{2})$ are topologically conjugate
with the help of a homeomorphism $h : A_{1} \rightarrow A_{2}$.

We have two projections $h : (A_{1}, g_{1}) \rightarrow (A_{2}, g_{2})$ and
$h^{-1} : (A_{2}, g_{2}) \rightarrow (A_{1}, g_{1})$. From
Corollary~\ref{corr_3.4} we conclude that $\Phi (\EssP{A_{1}, g_{1}}) =
\Phi (\EssP{A_{2}, g_{2}})$.

(b)
Let now $\Phi (\EssP{A_{1}, g_{1}}) = \Phi (\EssP{A_{2}, g_{2}}) = N
\in \Sigma$.

We find regular sequence $\{ n_{i} \}_{i \in \nn}$, for
which $\Phi (\{ n_{i} \; | \; i \in \nn \}) = N$ (see above) and construct on
it odometer $(A, g)$.

Since $\Phi (\EssP{A, g}) = \Phi (\EssP{A_{1}, g_{1}}) = \Phi (\EssP{A_{2},
g_{2}})$, then we receive from Proposition~\ref{prop_2.1} and Remark~\ref{rem_2.6}
the relation $\EssP{A, g} = \EssP{A_{1}, g_{1}} = \EssP{A_{2}, g_{2}}$. It follows from
Lemma~\ref{lemma_3.1} that $n_{i} \in \EssP{A, g}$, $i \in \nn$.
Consequently, $n_{i} \in \EssP{A_{k}, g_{k}}$, $k = 1, 2$, $i \in \nn$.

We fix coherent sequence $\{ W^{(n_{i})} \}_{i \in \nn}$ of periodic
partitions of dynamical system $(A_{1}, g_{1})$ and construct
on it the partition $\hh$ of space $A_{1}$.

Similarly, we take coherent sequence $\{
\widetilde{W}^{(n_{i})} \}_{i \in \nn}$ of periodic partitions of
dynamical system $(A_{2}, g_{2})$ and partition $\widetilde{\hh}$ of the
space $A_{2}$, induced by this sequence.

Consider the commutative diagram
$$
\begin{CD}
( A_{1}, g_{1}) @>{F}>> (A, g)
    @<{\widetilde{F}}<< (A_{2}, g_{2}) \\
@V{pr}VV                                 @|
                         @VV{\widetilde{pr}}V \\
( A_{1}/\hh, \fact g_{1})  @>{\fact F}>> (A, g)
    @<{\fact \widetilde{F}}<< (A_{2}/\widetilde{\hh}, \fact g_{2})
\end{CD}
$$

We know already, that all maps in the lower line of this diagram are
isomorphisms in category $\kk_{0}$.

Maps $pr$ and $\widetilde{pr}$ are one-to-one on
Proposition~\ref{prop_3.7}. Since spaces $A_{1}$ and $A_{2}$ are compact, then
$pr$ and $\widetilde{pr}$ are isomorphisms in category $ \kk_{0} $.

From what has been said it follows, that the morphism
$$
\widetilde{F}^{-1} \circ F = \widetilde{pr}^{-1} \circ
( \fact \widetilde{F})^{-1} \circ (\fact F) \circ pr
$$
is isomorphism and dynamical systems $(A_{1}, g_{1})$ and $(A_{2},
g_{2})$ are topologically conjugate.
\end{proof}

 %*********************************************************************%

\begin{rem} \label{rem_3.8}
Let $(X, f)$, $(Y, g) \in \kk_{0}$, $h : (X, f) \rightarrow (Y,
g)$ is a morphism. If dynamical system $(Y, g)$ is minimal,
then $h$ is the projection.
\end{rem}

Really, we fix $y \in Y$. Under the Birkgoff theorem we have
$\Cl{\Orb_{g} (h (y))} = Y$. On the other hand, since $X$ is compact set,
then $h(X)$ is the closed subset of space $Y$ and, obviously,
$h(X) \supseteq \Orb_{g}(h (y)) = h(\Orb_{f} (y))$.

\begin{prop} \label{prop_3.8}
Let $(A, g)$ is an odometer, $h : (A, g) \rightarrow (A, g)$ is a morphism.
Then $h$ is isomorphism.
\end{prop}

\begin{proof}
We fix regular sequence $\{ n_{i} \in \EssP{A, g}
\}_{i \in \nn}$, such that $\Phi (\{ n_{i} \; | \; i \in \nn \}) =
\Phi (\EssP{A, g})$ (see Remark~\ref{rem_2.6} and
Proposition~\ref{prop_2.2}).

Fix regular sequence of periodic partitions
$\{ W^{(n_{i})} \}_{i \in \nn}$.

$h$ is epimorphism according to Remark~\ref{rem_3.8}, therefore from
Proposition~\ref{prop_3.6} follows, that for every $i \in \nn$ the
family of sets $\widetilde{W}^{(n_{i})} =
\{\widetilde{W}^{(n_{i})}_{s_{i}} = h^{-1} (W^{(n_{i})}_{s_{i}})
\}_{i \in \zz_{n_{i}}}$ is periodic partition of
the dynamical system $(A, g)$ of length $n_{i}$.

Odometer $(A, g)$ is the minimal dynamical system, therefore $ (A,
g)$ is indecomposable. We apply Corollary~\ref{corr_1.1} and conclude that
$\{\widetilde{W}^{(n_{i})} \}_{i \in \nn}$ is regular
sequence of periodic partitions of $(A, g)$.

From Proposition~\ref{prop_3.7} it follows, that each of the families
of sets $\{ W^{(n_{i})}_{s_{i}} \; | \; s_{i} \in \zz_{n_{i}},\,
i \in \nn \}$ and
$\{\widetilde{W}^{(n_{i})}_{s_{i}} \; | \; s_{i} \in \zz_{n_{i}},\,
i \in \nn \}$ is basis of the topology of space $A$.
Therefore $h$ is continuous one-to-one map. Since
$A$ is compact set, then $h$ is homeomorphism.
\end{proof}

\begin{corr} \label{corr_3.5}
Let $(Y_{1}, h_{1})$ and $(Y_{2}, h_{2})$ are two dynamical systems,
which are topologically conjugate with some odometers.

If the objects $(Y_{1}, h_{1})$, $(Y_{2}, h_{2}) \in \Ob \kk_{0}$ are isomorphic,
then any morphism $\alpha : (Y_{1}, h_{1}) \rightarrow (Y_{2}, h_{2})$
is isomorphism.
\end{corr}

\begin{proof}
{\bf 1.} {} Let dynamical system $(Y_{1}, h_{1})$ is topologically
conjugate with an odometer $(A, g)$. Let $ \rho : (Y_{1}, h_{1}) \rightarrow
( Y_{1}, h_{1})$ is a morphism. Then $\rho$ is isomorphism.

Really, we fix isomorphism $ \phi : (Y_{1}, h_{1}) \rightarrow (A,
g)$ and consider morphism $\phi \circ \rho \circ \phi^{-1} : (A, g)
\rightarrow (A, g)$. According to Proposition~\ref{prop_3.8} $\phi \circ \rho
\circ \phi^{-1}$ is isomorphism. Then also $\rho = \phi^{-1} \circ (\phi \circ
\rho \circ \phi^{-1}) \circ \phi$ is isomorphism.

{\bf 2.} {}
Assume that there exists an isomorphism $\psi : (Y_{1}, h_{1}) \rightarrow
( Y_{2}, h_{2})$.

Morphism $\chi = \alpha \circ \psi^{-1} : (Y_{2}, h_{2}) \rightarrow
( Y_{2}, h_{2})$ is isomorphism (see above). Hence also
$\chi \circ \psi = \alpha \circ \psi^{-1} \circ \psi = \alpha$ is
isomorphism.
\end{proof}

Let $(A, +)$ is an adic group. Let
$g: A \rightarrow A $,
$$
g : \vec{a} \mapsto \vec{a} + \vec{e} \, \quad \vec{a} \in A \;.
$$
Let $(A, g)$ is the relevant odometer.

We fix $ \vec{a} $, $ \vec{b} \in A$. Consider the map
$$
h_{\vec{a}, \, \vec{b}}: A \rightarrow A \,,
$$
$$
h_{\vec{a}, \, \vec{b}} : \vec{c} \mapsto \vec{c} + (\vec{b} -
\vec{a}) \,, \quad \vec{c} \in A \;.
$$

\begin{prop} \label{prop_3.9}
$h_{\vec{a}, \, \vec{b}} \circ g = g \circ h_{\vec{a}, \, \vec{b}}$.
\end{prop}

\begin{proof}
This is the obvious corollary of commutability of group $(A, g)$.
\end{proof}

\begin{rem} \label{rem_3.9}
Let dynamical system $(X, f)$ is minimal, $h_{1}$,
$h_{2} : (X, f) \rightarrow (Y, g)$ are two morphisms, such that
$h_{1} (x) = h_{2} (x)$ for some $x \in X$.

Then $h_{1} = h_{2}$.
\end{rem}

Really, for every $y = f^{n} (x) \in \Orb_{f} (x)$ we have
$h_{1} (y) = h_{1} \circ f^{n} (x) = g^{n} \circ h_{1} (x) =
g^{n} \circ h_{2} (x) = h_{2} \circ f^{n} (x) = h_{2} (y)$. Therefore
$h_{1} |_{\Orb_{f} (x)} = h_{2} |_{\Orb_{f} (x)}$. Since
$X = \Cl{\Orb_{f} (x)}$ under Birkgoff theorem, then $h_{1} = h_{2}$.

\medskip

Thus, combining Propositions~\ref{prop_3.8}, \ref{prop_3.9} and
Remark~\ref{rem_3.9}, we receive
\begin{corr} \label{corr_3.6}
Let $(A, g)$ is a dynamical system, topologically conjugate with
an odometer. For any pair of points $x$, $y \in A$ there exists
unique morphism $h_{x, \, y} : (A, g) \rightarrow (A, g)$,
such that $h_{x, \, y} (x) = y$, and this morphism is
isomorphism.
\end{corr}

\subsection{Deviation. One categorial
construction.}\label{subsect_3.4}

Let $\LL$ is a category.

\begin{defn} \label{defn_3.5}
We say that $\LL$ has property $\LU$ (Lifting
Upstairs), if for every objects $A$, $B \in \Ob {\LL}$ and for
arbitrary morphisms $\alpha \in H_{\LL} (A, B)$ and $e_{B} \in
H_{\LL} (B, B) \cap \Iso{\LL}$ there exists morphism $e_{A} \in
H_{\LL} (A, A) \cap \Iso{\LL}$, such that
$$
e_{B} \circ \alpha = \alpha \circ e_{A} \;.
$$
\end{defn}

\begin{defn} \label{defn_3.6}
We say that $\LL$ has property $\LD$ (Lifting
Downstairs), if for any objects $A $, $B \in \Ob {\LL}$ and for
arbitrary morphisms $\alpha \in H_{\LL} (A, B)$ and $f_{A} \in
H_{\LL} (A, A) \cap \Iso{\LL}$ there exists morphism $f_{B} \in
H_{\LL} (B, B) \cap \Iso{\LL}$, such that
$$
\alpha \circ f_{A} = f_{B} \circ \alpha \;.
$$
\end{defn}

Let $\LL$ is a category. For each pair of objects $A $, $B \in
\Ob {\LL}$ we shall define binary relation $\sim$ on the set $H_{\LL} (A, B)$.
Say that $\alpha \sim \beta$, $\alpha$,
$\beta \in H_{\LL} (A, B)$, if there exist such $e_{A} \in
H_{\LL} (A, A) \cap \Iso{\LL}$ and $e_{B} \in
H_{\LL} (B, B) \cap \Iso{\LL}$, that
$$
\alpha \circ e_{A} = e_{B} \circ \beta \;.
$$

It is easy to see, that $\sim$ is the equivalence relation. An equivalence
class of a morphism $\alpha$ shall be designated by $[\alpha]$.

\begin{prop} \label{prop_3.10}
Assume that category $\LL$ has one of properties $\LU$ or $\LD$.

Then the category $\Bar{\LL}$ is correctly defined, for which the objects are
same with the objects of category $\LL$ and for any pair of objects
$A$, $B \in \Bar{\LL}$ a set of morphisms $H_{\Bar{\LL}} (A, B)$
is the set of equivalence classes of morphisms from $H_{\LL} (A,
B)$.
\end{prop}

\begin{proof}
Assume, that the category $\LL$ has property $\LU$.

It is trivially checked, that $\Bar{\LL}$ satisfies to properties 1)
and 2) of category.

In order to define correctly composition of morphisms in
$\Bar{\LL}$, we shall prove that for any triple of objects $A$, $B$,
$C \in \LL$ and for any morphisms $\alpha \in H_{\LL} (A, B)$,
$\beta \in H_{\LL} (B, C)$ the equality is fulfilled

\begin{equation}\label{eq_15}
[\beta \circ \alpha] = [\beta] \circ [\alpha] = \{ \beta' \circ \alpha'
\; | \; \alpha'  \in [\alpha], \, \beta'  \in [\beta] \} \;.
\end{equation}

Let $\alpha'  \in [\alpha]$, $\beta'  \in [\beta]$.
Then there exist such isomorphisms
$e_{A} \in H_{\LL} (A, A) \cap \Iso{\LL}$,
$e_{B}$, $f_{B} \in H_{\LL} (B, B) \cap \Iso{\LL}$ and
$f_{C} \in H_{\LL} (C, C) \cap \Iso{\LL}$, that
$$
\beta' \circ \alpha' = (f_{C} \circ \beta \circ f_{B}^{-1}) \circ
(e_{B} \circ \alpha \circ e_{A}^{-1}) = f_{C} \circ \beta \circ
(f_{B}^{-1} \circ e_{B}) \circ \alpha \circ e_{A}^{-1} \;.
$$
Obviously, $f_{B}^{-1} \circ e_{B} \in H_{\LL} (B, B) \cap \Iso{\LL}$.
From property $\LU$ we conclude, that there exists $g_{A} \in
H_{\LL} (A, A) \cap \Iso{\LL}$, for which
$(f_{B}^{-1} \circ e_{B}) \circ \alpha = \alpha \circ g_{A}$. Hence,
$$
\beta' \circ \alpha' = f_{C} \circ \beta \circ (f_{B}^{-1} \circ e_{B})
\circ \alpha \circ e_{A}^{-1} = f_{C} \circ \beta \circ \alpha \circ
(g_{A} \circ e_{A}^{-1}) \;,
$$
$\beta' \circ \alpha' \in [\beta \circ \alpha]$ and
$[\beta \circ \alpha] \supseteq \{ \beta' \circ \alpha' \; | \;
\alpha'  \in [\alpha], \, \beta'  \in [\beta] \}$.

Inversely, let $\gamma \in [\beta \circ \alpha]$. It means, that for
some $e_{A} \in H_{\LL} (A, A) \cap \Iso{\LL}$ and
$e_{C} \in H_{\LL} (C, C) \cap \Iso{\LL}$ the relation takes place
$$
\gamma = e_{C} \circ (\beta \circ \alpha) \circ e_{A}^{-1} =
(e_{C} \circ \beta) \circ (\alpha \circ e_{A}^{-1}) \;.
$$
Obviously, $\alpha \circ e_{A}^{-1} = 1_{B} \circ \alpha \circ e_{A}^{-1}
\in [\alpha]$ and $e_{C} \circ \beta = e_{C} \circ \beta \circ 1_{B}^{-1}
\in [\beta]$. Hence, $[\beta \circ \alpha] \subseteq
\{ \beta' \circ \alpha' \; | \; \alpha'  \in [\alpha], \,
\beta'  \in [\beta] \}$.

So, we have established, that partial multiplication of equivalence classes
of morphisms does not depend on a choice of representatives,
hence it is defined correctly.

Associativity of multiplication of morphisms in $\Bar{\LL}$ follows from
associativity of multiplication of morphisms in $\LL$.

To complete the proof it suffices to note, that for
any $A \in \Ob {\LL}$ the unit morphism of object $A$ in
$\Bar{\LL}$ is $[1_{A}]$.

If the category $\LL$ has property $\LD$, proof
is conducted similarly.
\end{proof}

\begin{rem}
The described above construction is a special case of so-called
{\em factor--category} (see~\cite{MacLane}).
\end{rem}

\subsection{Main properties of odometers (continuation).} \label{subsect_3.5}

Let $(A_{1}, g_{1})$, $(A_{2}, g_{2})$ are dynamical systems,
topologically conjugate with odometers, $\pi_{1}$, $\pi_{2} :
(A_{1}, g_{1}) \rightarrow (A_{2}, g_{2})$ and $h : (A_{2}, g_{2})
\rightarrow (A_{2}, g_{2})$ are morphisms.

Designate by $\Fcal$ the set of all morphisms $f : (A_{1}, g_{1})
\rightarrow (A_{1}, g_{1})$, such that the diagram is commutative
$$
\begin{CD}
(A_{1}, g_{1}) @>{f}>> (A_{1}, g_{1}) \\
@V{\pi_{1}}VV @VV{\pi_{2}}V \\
(A_{2}, g_{2}) @>>{h}> (A_{2}, g_{2})
\end{CD}
$$

It follows from Proposition~\ref{prop_3.8} that $h \in \Iso{\kk_{0}}$ and
$\Fcal \subseteq \Iso{\kk_{0}}$.

\begin{prop} \label{prop_3.11}
The set $\Fcal$ is not empty.

For any $y \in A_{2}$ and $x_{1} \in \pi_{1}^{-1} (y)$ the
equality is valid
$$
\Fcal = \{ h_{x_{1}, \, x_{2}} \; | \; x_{2} \in \pi_{2}^{-1} (h (y)) \}
\;.
$$
\end{prop}

Proof of Proposition~\ref{prop_3.11} is based on the following
Lemmas.

\begin{lemma} \label{lemma_3.3}
Let $(X, f)$ is an indecomposable dynamical system, $(A, g)$ is
dynamical system, topologically conjugate with an odometer,
$\pi_{1}$, $\pi_{2} : (X, f) \rightarrow (A, g)$ are projections.

Then the partitions $\zer \pi_{1}$ and $\zer \pi_{2}$ of space $X$
coincide.
\end{lemma}

\begin{lemma} \label{lemma_3.4}
Let $(X, f)$, $(A, g)$, $\pi_{1}$, $\pi_{2} : (X, f) \rightarrow
(A, g)$ are the same, as in Lemma~\ref{lemma_3.3}.

For any morphism $h : (X, f) \rightarrow (X, f)$
the continuous map $\fact \pi_{2} \circ h : A \rightarrow A$ is well defined
to comply with the relation $g \circ (\fact \pi_{2} \circ h) =
( \fact \pi_{2} \circ h) \circ g$.

Therefore, the commutative diagram is valid
$$
\begin{CD}
(X, f) @>{h}>> (X, f) \\
@V{\pi_{1}}VV @VV{\pi_{2}}V \\
(A, g) @>>{\fact \pi_{2} \circ h}> (A, g)
\end{CD}
$$
\end{lemma}

\begin{proof}[of Lemma~\ref{lemma_3.3}]
We fix a regular sequence
$\{ n_{i} \in \Phi (\EssP{A, g}) \}_{i \in \nn} $, such that
$\Phi (\{ n_{i} \; | \; i \in \nn \}) = \Phi (\EssP{A, g})$
(see Remark~\ref{rem_2.6} and Proposition~\ref{prop_2.2}).

We construct regular sequence $\{ W^{(n_{i})} \}_{i \in
\nn}$ of periodic partitions of dynamical system $(A, g)$.
According to Proposition~\ref{prop_3.7} the family of sets
$\{ W^{(n_{i})}_{s_{i}} \; | \; s_{i} \in \zz_{n_{i}} \,, \; i
\in \nn \}$ is basis of the topology of space $A$.

It follows from Proposition~\ref{prop_3.6} that for every
$i \in \nn$ the systems of sets
$\widetilde{W}^{(n_{i})} = \{ \widetilde{W}^{(n_{i})}_{s_{i}} =
\pi_{1}^{-1} (W^{(n_{i})}_{s_{i}}) \; | \;
s_{i} \in \zz_{n_{i}} \}$ and
$\widehat{W}^{(n_{i})} = \{ \widehat{W}^{(n_{i})}_{s_{i}} =
\pi_{2}^{-1} (W^{(n_{i})}_{s_{i}}) \; | \;
s_{i} \in \zz_{n_{i}} \}$ are periodic partitions of
dynamical system $(X, f)$ of length $n_{i}$.

Dynamical system $(X, f)$ is indecomposable, therefore we conclude from
Corollary~\ref{corr_1.1} and Remark~\ref{rem_1.12}
that $\{ \widetilde{W}^{(n_{i})} \}_{i \in \nn}$ and
$\{ \widehat{W}^{(n_{i})} \}_{i \in \nn}$ are compatible regular
sequences of periodic partitions of $(X, f)$.

Let $\widetilde{\hh}$ and $\widehat{\hh}$ are the partitions of space
$X$, induced respectively by sequences
$\{ \widetilde{W}^{(n_{i})} \}_{i \in \nn}$ and
$\{ \widehat{W}^{(n_{i})} \}_{i \in \nn}$.
Then partitions $\widetilde{\hh}$ and
$\widehat{\hh}$ coincide by Corollary~\ref{corr_3.2}.

The family of sets $\{ \widetilde{W}^{(n_{i})}_{s_{i}}
\; | \; s_{i} \in \zz_{n_{i}} \,, \; i \in \nn \}$ is
pre-image of the basis of topology $\{ W^{(n_{i})}_{s_{i}} \; | \;
s_{i} \in \zz_{n_{i}} \,, \; i \in \nn \}$ under action of the projection
$\pi_{1}$, therefore partitions $\zer \pi_{1}$ and $\widetilde{\hh}$
coincide.

Similarly, the partitions $\zer \pi_{2}$ and $\widehat{\hh}$ coincide.
\end{proof}

\begin{proof}[of Lemma~\ref{lemma_3.4}]
Let $h : (X, f) \rightarrow (X, f)$ is a morphism.

Since dynamical system $(A, g)$ is minimal, then morphism
$\pi_{2} \circ h : (X, f) \rightarrow (A, g)$ is the projection.
From Lemma~\ref{lemma_3.3} it follows, that $\zer \pi_{1} =
\zer \pi_{2} \circ h$. To complete the proof it remains
to apply Lemma~\ref{lemma_2} to morphisms $\varphi_{1} = \pi_{1}$ and
$\varphi_{2} = \pi_{2} \circ h$.
\end{proof}

\begin{proof}[of Proposition~\ref{prop_3.11}] 
Dynamical system $(A_{1}, g_{1})$ is minimal, therefore it
is indecomposable.

We fix $y \in A_{2}$ and $x_{1} \in \pi_{1}^{-1} (y)$. Let
$x_{2} \in \pi_{2}^{-1} (h (y))$. Consider morphism
$h_{x_{1}, \, x_{2}} : (A_{1}, g_{1}) \rightarrow (A_{1}, g_{1})$.
We conclude from Lemma~\ref{lemma_3.4} that morphism
$\fact \pi_{2} \circ h_{x_{1}, \, x_{2}}: (A_{2}, g_{2})
\rightarrow (A_{2}, g_{2})$ is correctly defined and the sequence
of equalities is valid
$$
\fact \pi_{2} \circ h_{x_{1}, \, x_{2}} (y) =
(\fact \pi_{2} \circ h_{x_{1}, \, x_{2}}) \circ \pi_{1} (x_{1}) =
\pi_{2} \circ h_{x_{1}, \, x_{2}} (x_{1}) =
\pi_{2} (x_{2}) = h (y) \;.
$$
Therefore, it follows from Corollary~\ref{corr_3.6} that
$\fact \pi_{2} \circ h_{x_{1}, \, x_{2}} = h$ and $\Fcal \supseteq
\{ h_{x_{1}, \, x_{2}} \; | \; x_{2} \in \pi_{2}^{-1} (h (y)) \}$.

On the other hand, for any $f \in \Fcal$ the relation
$f (x_{1}) \in f (\pi_{1}^{-1} (y)) = \pi_{2}^{-1} (h (y))$ should be fulfilled
(see Lemma~\ref{lemma_3.3}). Corollary~\ref{corr_3.6}
consequences that $f = h_{x_{1}, \, f (x_{1})}$ and $\Fcal \subseteq
\{ h_{x_{1}, \, x_{2}} \; | \; x_{2} \in \pi_{2}^{-1} (h (y)) \}$.
\end{proof}

Consider the complete subcategory $\Acal$ of $\kk_{0}$,
objects of which are all dynamical systems, topologically
conjugate with odometers.

\begin{corr} \label{corr_3.7}
The category $\Acal$ has properties $\LU$ and $\LD$.
\end{corr}

\begin{proof}
{\bf 1.}
Condition $\LU$ follows from Proposition~\ref{prop_3.11},
with $\pi_{1} = \pi_{2} = \alpha$, and from Corollary~\ref{corr_3.6}.

{\bf 2.}
Condition $\LD$ is the consequence of Lemma~\ref{lemma_3.4},
where $(X, f)$ and $(A, g)$ are topologically conjugate with odometers and
$\pi_{1} = \pi_{2} = \alpha$, and from Corollary~\ref{corr_3.6}.
\end{proof}

We fix a skeleton $\Acal_{0}$ of category $\Acal$.

\begin{rem} \label{rem_3.10}
From Theorem~\ref{theorem_3.2} and Remark~\ref{rem_3.5} we conclude
that for every $N \in \Sigma$ category $\Acal_{0}$ contains
exactly one object $(A, g)$, such that $\Phi (\EssP{A, g}) = N$.
\end{rem}

\begin{rem} \label{rem_3.11}
Theorem~\ref{theorem_A} consequences that for any two objects $(A_{1},
g_{1})$ and $(A_{2}, g_{2})$ of category $\Acal_{0}$ the following statements
are equivalent:
\begin{itemize}
    \item[(i)] $\Phi (\EssP{A_{1}, g_{1}}) \geq \Phi (\EssP{A_{2},
        g_{2}})$;
    \item[(ii)] $H_{\Acal_{0}} ((A_{1}, g_{1}), (A_{2}, g_{2}))
        \neq \emptyset$.
\end{itemize}
\end{rem}

We take advantage of Proposition~\ref{prop_3.10} and construct category $\Bar{\Acal_{0}}$
by the category $\Acal_{0}$.

From Lemma~\ref{lemma_3.4} and Corollary~\ref{corr_3.6} we obtain
\begin{corr} \label{corr_3.8}
Let $(A_{1}, g_{1})$, $(A_{2}, g_{2}) \in \Ob{\Acal_{0}}$,
$\pi_{1}$, $\pi_{2} \in H_{\Acal_{0}} ((A_{1}, g_{1}),
(A_{2}, g_{2}))$. Then $[\pi_{1}] = [\pi_{2}]$.
\end{corr}

Now from Remark~\ref{rem_3.11} we have

\begin{corr} \label{corr_3.9}
Let $(A_{1}, g_{1})$, $(A_{2}, g_{2}) \in \Ob{\Bar{\Acal_{0}}}$.

Then
\begin{itemize}
    \item[(i)] if $\Phi (\EssP{A_{1}, g_{1}}) \geq
        \Phi (\EssP{A_{2}, g_{2}})$, then the set
        $H_{\Bar{\Acal_{0}}} ((A_{1}, g_{1}) (, A_{2}, g_{2}))$
        contains exactly one element;
    \item[(ii)] $H_{\Bar{\Acal_{0}}} ((A_{1}, g_{1}),
        (A_{2}, g_{2})) = \emptyset$, otherwise.
\end{itemize}
\end{corr}

By partially ordered set $\Sigma$ we construct category
$\LL (\Sigma)$, which objects are the elements of the set
$\Sigma$, and morphisms are all possible pairs of elements $(M, N)$,
such that $M \geq N$. For any two elements $M$, $N \in \Sigma$
the set $H_{\LL (\Sigma)} (M, N)$ consists of one morphism
$(M, N)$ if $M \geq N$, otherwise this set is empty.

Remark~\ref{rem_3.10} and corollary~\ref{corr_3.9} give us
the following
\begin{theorem} \label{theorem_3.3}
Correspondence $\Psi_{0} : \Ob{\Bar{\Acal_{0}}} \rightarrow
\Ob{\LL (\Sigma)}$,
$$
\Psi_{0} : (A, g) \mapsto \Phi (\EssP{A, g}) \,, \quad
(A, g) \in \Ob{\Bar{\Acal_{0}}} \,,
$$
is uniquely extended to functor $\Psi : \Bar{\Acal_{0}}
\rightarrow \LL (\Sigma)$.

The functor $\Psi$ sets the isomorphism of categories $\Bar{\Acal_{0}}$
and $\LL(\Sigma)$.
\end{theorem}

\section{Expansions of odometers.}

\subsection{General case.}\label{subsect_3.6}

Let $(X, f)$ is a dynamical system with Hausdorff compact phase
space.

Remind, that for any dynamical system $(A, g) \in \Ob{\Acal}$
category $\kk_{0}$ contains a morphism $h : (X, f) \rightarrow (A, g)$
if and only if $\Phi (\EssP{A, g}) \leq \Phi (\EssP{X, f})$
(see Theorem~\ref{theorem_A}).

Consider a complete subcategory $\Acal (X, f)$ of $\Acal$, objects of
which are dynamical systems $(A, g) \in \Ob{\Acal}$, such that
$\Phi (\EssP{A, g}) \leq \Phi (\EssP{X, f})$.

It follows from Remark~\ref{rem_3.10} that the category $\Acal (X, f)$ has
properties $\LU$ and $\LD$.

We fix a skeleton $\Acal_{0} (X, f)$ of $\Acal (X, f)$ and construct category
$\Bar{\Acal_{0} (X, f)}$, taking an advantage of Proposition~\ref{prop_3.10}.

In the same way as we did it for Corollary~\ref{corr_3.9}, we prove
\begin{prop} \label{prop_3.12}
Let $(A_{1}, g_{1})$, $(A_{2}, g_{2}) \in \Ob{\Bar{\Acal_{0} (X, f)}}$.

Then
\begin{itemize}
    \item[(i)] if $\Phi (\EssP{A_{1}, g_{1}}) \geq
        \Phi (\EssP{A_{2}, g_{2}})$, then the set
        $H_{\Bar{\Acal_{0} (X, f)}} ((A_{1}, g_{1}), (A_{2}, g_{2}))$
        contains exactly one element;
    \item[(ii)] $H_{\Bar{\Acal_{0} (X, f)}} ((A_{1}, g_{1}),
        (A_{2}, g_{2})) = \emptyset$, otherwise.
\end{itemize}
\end{prop}

Consider the subset $\Sigma (X, f) = \{ N \in \Sigma \; | \; N \leq
\Phi (\EssP{X, f}) \}$ of $(\Sigma, \leq)$ and construct category
$\LL (\Sigma (X, f))$ by this partially ordered set.

Similarly to Theorem~\ref{theorem_3.3} we prove
\begin{theorem} \label{theorem_3.4}
Correspondence $\Psi_{0} (X, f): \Ob{\Bar{\Acal_{0} (X, f)}} \rightarrow
\Ob{\LL (\Sigma (X, f))}$,
$$
\Psi_{0} (X, f) : (A, g) \mapsto \Phi (\EssP{A, g}) \,, \quad
( A, g) \in \Ob{\Bar{\Acal_{0} (X, f)}} \,,
$$
is uniquely extended to functor $\Psi (X, f) : \Bar{\Acal_{0} (X, f)}
\rightarrow \LL (\Sigma (X, f))$.

The functor $\Psi (X, f)$ sets the isomorphism of categories $\Bar{\Acal_{0} (X, f)}$ and
$\LL (\Sigma (X, f))$.
\end{theorem}

We define category $\BB (X, f)$ as follows:
\begin{itemize}
    \item[---] objects of $\BB (X, f)$ are the pairs $(h, (A, g))$, where
        $h \in \Mor_{\kk_{0}} ((X, f),
        (A, g))$ and $(A, g) \in \Ob{\Acal (X, f)}$;
    \item[---] for any two $(h_{1}, (A_{1}, g_{1}))$,
        $(h_{2}, (A_{2}, g_{2})) \in \Ob{\BB (X, f)}$ the set
        $H_{\BB (X, f)} ((h_{1}, (A_{1}, g_{1})), (h_{2}, (A_{2}, g_{2})))$
        consists of all triples $((h_{1}, (A_{1}, g_{1})), \pi,
        (h_{2}, (A_{2}, g_{2})))$, such that $\pi \in
        H_{\Acal (X, f)} ((A_{1}, g_{1}), (A_{2}, g_{2}))$ and
        $h_{2} = \pi \circ h_{1}$.
\end{itemize}

The correctness of this definition is checked immediately.

\begin{lemma} \label{lemma_3.5}
Let $(h_{1}, (A_{1}, g_{1}))$, $(h_{2}, (A_{2}, g_{2})) \in
\Ob{\BB (X, f)}$.

Then
\begin{itemize}
    \item[(i)] the set
        $H_{\BB (X, f)} ((h_{1}, (A_{1}, g_{1})),
        (h_{2}, (A_{2}, g_{2})))$ contains no more than one element;
    \item[(ii)] if $H_{\BB (X, f)} ((h_{1}, (A_{1}, g_{1})),
        (h_{2}, (A_{2}, g_{2}))) \neq \emptyset$, then
        $\Phi (\EssP{A_{1}, g_{1}}) \geq \Phi (\EssP{A_{2}, g_{2}})$.
\end{itemize}
\end{lemma}

\begin{proof}
{\bf (i)}
Let $\pi_{1}$, $\pi_{2} : (A_{1}, g_{1}) \rightarrow (A_{2}, g_{2})$ are
two morphisms, such that $h_{2} = \pi_{1} \circ h_{1} = \pi_{2} \circ
h_{1}$.

Let $x \in X$. Then $h_{2} (x) = \pi_{1} (h_{1} (x)) =
\pi_{2} (h_{1} (x))$. Since dynamical system $(A_{1}, g_{1})$
is minimal, then it follows from Remark~\ref{rem_3.9} that $\pi_{1} =
\pi_{2}$.

{\bf (ii)}
This statement immediately follows from Theorem~\ref{theorem_A}.
\end{proof}

\begin{lemma} \label{lemma_3.6}
Let $(h, (A, g)) \in \Ob{\BB (X, f)}$, $N \in \Sigma (X, f)$ and
$\Phi (\EssP{A, g}) \leq N$.

Then there exists $(h_{1}, (A_{1}, g_{1})) \in \Ob{\BB (X, f)}$, such that
\begin{itemize}
    \item[(i)] $\Phi (\EssP{A_{1}, g_{1}}) = N$;
    \item[(ii)] $H_{\BB (X, f)} ((h_{1}, (A_{1}, g_{1})),
        (h, (A, g))) \neq \emptyset$.
\end{itemize}
\end{lemma}

\begin{proof}
On the condition of Lemma and by definition of the set $\Sigma (X, f)$ we have
following inequalities $\Phi (\EssP{A, g}) \leq N \leq \Phi (\EssP{X, f})$.

Proposition~\ref{prop_2.1} consequences that there exists unique
regular set $Q \in \nn$, such that $\Phi (Q) = N$. Since
the sets $\EssP{A, g}$ and $\EssP{X, f}$ are regular
(see Remark~\ref{rem_2.6}), then from Lemma~\ref{lemma_2.1} we receive the
inclusions $\EssP{A, g} \subseteq Q \subseteq \EssP{X, f}$.

We use Proposition~\ref{prop_2.2} and find regular
sequences $\{ n_{i} \in \EssP{A, g} \}_{i \in \nn}$ and
$\{ m_{j} \in Q \}_{j \in \nn}$, such that
$\Phi (\{ n_{i} \; | \; i \in \nn \}) = \Phi (\EssP{A, g})$ and
$\Phi (\{ m_{j} \; | \; j \in \nn \}) = \Phi (Q) = N$.

Fix regular sequence $\{ V^{(n_{i})} \}_{i \in \nn}$ of
periodic partitions of dynamical system $(A, g)$. Consider
pre-images $W^{(n_{i})} = \{ W^{(n_{i})}_{s_{i}} =
h^{-1} (V^{(n_{i})}_{s_{i}}) \}_{s_{i} \in \zz_{n_{i}}}$, $i \in \nn$, of
periodic partitions of the sequence $\{ V^{(n_{i})} \}_{i \in
\nn}$. We conclude from Proposition~\ref{prop_3.6}, Definition~\ref{defn_1.6} and
Corollaries~\ref{corr_1.5} and~\ref{corr_1.6} that
$\{ W^{(n_{i})} \}_{i \in \nn}$ is the regular sequence of
periodic partitions of dynamical system $(X, f)$.

Let $\hh$ is the partition of space $X$, induced
by the sequence $\{ W^{(n_{i})} \}$. In accord with
Proposition~\ref{prop_3.7} the family of sets
$\{ V^{(n_{i})}_{s_{i}} \; | \; s_{i} \in \zz_{n_{i}}, \, i \in \nn \}$
is basis of the topology of space $A$. Therefore, $\hh$
coincides with the partition $\zer h$.

We take advantage of Proposition~\ref{prop_1.8} and Remark~\ref{rem_1.6}
and construct coherent sequence
$\{ \widetilde{W}^{(m_{j})} \}_{j \in \nn}$ of periodic partitions of
dynamical system $(X, f)$, which is compatible with the sequence
$\{ W^{(n_{i})} \}$.

Let $\widetilde{\hh}$ is the partition of space $X$, induced
by the sequence $\{ \widetilde{W}^{(m_{j})} \}$. Taking into account
Remark~\ref{rem_3.4}, we have the commutative diagram
(see relation~(\ref{eq_12}))
\begin{equation} \label{eq_16}
\begin{CD}
(X, f) @>{F}>> (A_{1}, g_{1}) \\
@V{pr}VV     @| \\
(X / \widetilde{\hh}, \widetilde{f}) @>{\fact F}>> (A_{1}, g_{1})
\end{CD}
\end{equation}
In this diagram $(A_{1}, g_{1})$ is an odometer, $F$ is a projection, $\fact F$ is
homeomorphism and
$\Phi (\EssP{X / \widetilde{\hh}, \widetilde{f}}) =
\Phi (\EssP{A_{1}, g_{1}}) = \Phi (\{ m_{j} \; | \; j \in \nn \}) = N$.

From Corollary~\ref{corr_3.2} we conclude that the partition $\widetilde{\hh} =
\zer F$ is refinement of the partition $\hh = \zer h$ of space $X$.
Designate $h_{1} = F$. Applying now Lemma~\ref{lemma_2} to projections
$\varphi_{1} = h_{1}$ and $\varphi_{2} = h$ we consequence that there exists $\psi
\in H_{\kk_{0}} ((A_{1}, g_{1}), (A, g))$, such that $h = \psi \circ
h_{1}$.
\end{proof}

We define the ``forgetful'' functor $\Theta : \BB (X, f) \rightarrow
\Acal (X, f)$ with the help of the relations
$$
\Theta : (h, (A, g)) \mapsto (A, g) \,, \quad
(h, (A, g)) \in \Ob{\BB (X, f)} \;;
$$
$$
\Theta : ((h_{1}, (A_{1}, g_{1})), \pi, (h_{2}, (A_{2}, g_{2}))) \mapsto
\pi \,, \quad
((h_{1}, (A_{1}, g_{1})), \pi, (h_{2}, (A_{2}, g_{2})))
\in \Mor (\BB (X, f)) \;.
$$

Consider pre-image $\BB' (X, f)$ of the skeleton $\Acal_{0} (X, f)$ under
action of $\Theta$. The easy immediate verification shows
that $\BB' (X, f)$ is complete subcategory of $\BB (X, f)$.

\begin{rem} \label{rem_3.13}
Let $ \LL' $ is a complete subcategory of a category $\LL$. Let $A$, $B \in \Ob
\LL'$. By definition we have $H_{\LL'} (A, B) = H_{\LL} (A, B)$.

Hence, the objects $A$ and $B$ are isomorphic in $ \LL'$ if and only if
they are isomorphic in $\LL$.
\end{rem}

Let $\BB_{0}' (X, f)$ is a skeleton of $ \BB' (X, f) $. Obviously,
$\BB_{0}' (X, f)$ is the complete subcategory of $\BB (X, f)$.

\begin{prop} \label{prop_3.13}
Category $\BB_{0}' (X, f)$ is the skeleton of the category $\BB (X, f)$.
\end{prop}

\begin{proof}
It follows from definition of the category $\BB_{0}' (X, f)$ and from Remark~\ref{rem_3.13}
that this subcategory contains no more than one
representative from every class of isomorphic objects of the category $\BB (X,
f)$.

We shall prove that for every $(h, (A, g)) \in \Ob \BB (X, f)$ there exists
an object contained in $\BB_{0}' (X, f)$, which is isomorphic to $(h, (A, g))$.

Consider $(A, g) = \Theta (h, (A, g)) \in \Ob \Acal (X, f)$. Subcategory
$\Acal_{0} (X, f)$ is the skeleton of $ \Acal(X, f)$, therefore
there exists exactly one object $(A', g') \in \Ob \Acal_{0} (X, f)$,
which is isomorphic to $(A, g)$. Let $\rho : (A, g) \rightarrow (A', g')$ is an
isomorphism.

Designate $h' = \rho \circ h : (X, f) \rightarrow (A', g')$. Consider
the object $(h', (A', g'))$ of category $\BB (X, f)$. Obviously,
$(h', (A', g')) \in \Ob \BB_{0}' (X, f)$ and $((h, (A, g)), \rho,
(h', (A', g'))) \in \Iso \BB (X, f)$.
\end{proof}

We shall introduce now binary relation $\preceq$ on the class $\Ob{\BB_{0}' (X, f)}$.
We say that $(h_{1}, (A_{1}, g_{1})) \preceq (h_{2}, (A_{2}, g_{2}))$
if $H_{\BB (X, f)} ((h_{2}, (A_{2}, g_{2})) (, h_{1}, (A_{1}, g_{1})))
\neq \emptyset$.

\begin{prop} \label{prop_3.14}
Relation $\preceq$ is the partial order and the following statements
hold true:
\begin{itemize}
    \item[(i)] each element $(h, (A, g)) \in \Ob \BB_{0}' (X, f)$
        is majorized by certain element $(h', (A', g')) \in
        \Ob \BB_{0}' (X, f)$, such that $\Phi (\EssP{A', g'}) =
        \Phi (\EssP{X, f})$;
    \item[(ii)] an element $(h, (A, g)) \in \Ob \BB_{0}' (X, f)$ is
        maximal under order relation~$\preceq$ if and only if
        $\Phi (\EssP{A, g}) = \Phi (\EssP{X, f})$.
\end{itemize}
\end{prop}

Before we proceed with the proof of Proposition~\ref{prop_3.14}, let us
prove

\begin{lemma} \label{lemma_3.7}
Let $(h, (A, g)) \in \Ob \BB_{0}' (X, f)$, $N \in \Sigma (X, f)$ and
$\Phi (\EssP{A, g}) \leq N$.

There exists such object $(h', (A', g')) \in \Ob \BB_{0}' (X, f)$, that
\begin{itemize}
    \item[(i)] $\Phi (\EssP{A', g'}) = N$;
    \item[(ii)] $H_{\BB (X, f)} ((h', (A', g')), (h, (A, g))) \neq
        \emptyset$.
\end{itemize}
\end{lemma}

\begin{proof}
According to Lemma~\ref{lemma_3.6} there exists $(h_{1}, (A_{1}, g_{1}))
\in \Ob \BB (X, f)$, such that $\Phi (\EssP{A_{1}, g_{1}}) = N$ and
$H_{\BB (X, f)} ((h_{1}, (A_{1}, g_{1})), (h, (A, g))) \neq \emptyset$.

We conclude from Proposition~\ref{prop_3.13} that there exists an object
$(h', (A', g')) \in \Ob \BB_{0}' (X, f)$, which is isomorphic to
$(h_{1}, (A_{1}, g_{1}))$. Lemma~\ref{lemma_3.5} consequences that
$\Phi (\EssP{A', g'}) = \Phi (\EssP{A_{1}, g_{1}}) = N$. In addition, we have
$H_{\BB (X, f)} ((h', (A', g')), (h, (A, g))) \neq \emptyset$ on
the construction.
\end{proof}

\begin{proof}[of Proposition~\ref{prop_3.14}]
For every object $B \in \Ob \BB_{0}' (X, f)$ there exists unit
morphism $1_{B}$ on definition of category, therefore relation
$\preceq$ is reflexive.

Composition of any morphisms $\alpha \in H_{\BB (X, f)} (B, B') $ and
$\beta \in H_{\BB (X, f)} (B', B'') $ is element of the set
$H_{\BB (X, f)} (B, B'')$, hence relation $\preceq$ is transitive.

Let $\alpha \in H_{\BB (X, f)} (B, B')$, $\beta \in H_{\BB (X, f)} (B', B)$
for certain $B$, $B' \in \Ob \BB_{0}' (X, f)$. According
to Lemma~\ref{lemma_3.5} we have $H_{\BB (X, f)} (B, B) = \{ 1_{B} \}$,
$H_{\BB (X, f)} (B', B') = \{ 1_{B'} \}$. Hence, $\beta \circ \alpha =
1_{B}$, $\alpha \circ \beta = 1_{B'}$ and the objects $B$ and $B'$ are
isomorphic in category $\BB (X, f)$. Since $\BB_{0}' (X, f)$ is
the complete subcategory of $\BB (X, f)$, then $B$ and $B'$ are isomorphic
in $\BB_{0}' (X, f)$ (see. Remark~\ref{rem_3.13}). From what has
been said we conclude, that the relation $\preceq$ is antisymmetric.

So, $\preceq$ is the partial order relation.

Statement~(i) of Proposition~\ref{prop_3.14} follows from
Lemma~\ref{lemma_3.7} and Theorem~\ref{theorem_A}.

We shall prove now statement~(ii).

Let $(h, (A, g))$, $(h', (A', g')) \in \Ob \BB_{0}' (X, f)$,
$\Phi (\EssP{A, g}) = \Phi (\EssP{X, f})$ and $(h, (A, g)) \preceq
(h', (A', g'))$. From Lemma~\ref{lemma_3.5} and Theorem~\ref{theorem_A}
we conclude that $\Phi (\EssP{A, g}) \leq \Phi (\EssP{A', g'}) \leq
\Phi (\EssP{X, f})$. Hence, $\Phi (\EssP{A', g'}) =
\Phi (\EssP{X, f}) = \Phi (\EssP{A, g})$ and any morphism $\rho : (A', g')
\rightarrow (A, g)$ is isomorphism (see Theorem~\ref{theorem_3.2}
and corollary~\ref{corr_3.5}).

Since the set $H_{\BB (X, f)} ((h', (A', g')), (h, (A, g)))$ is not empty on
our supposition, then the objects $(h', (A', g'))$, $(h, (A, g)) \in
\Ob \BB_{0}' (X, f)$ are isomorphic. Hence, $(h', (A', g')) =
(h, (A, g))$ (see Proposition~\ref{prop_3.13}) and $(h, (A, g))$ is the
maximal element under relation $\preceq$.

Let now $(h, (A, g))$ is a maximal element under
$\preceq$. The equality $\Phi (\EssP{A, g}) = \Phi (\EssP{X, f})$ follows from
the statement~(i) of Proposition~\ref{prop_3.14} and Theorem~\ref{theorem_A}.
\end{proof}

On definition we have $\Theta (\Ob \BB_{0}' (X, f)) = \Ob \Acal_{0} (X, f) =
\Ob \Bar{\Acal_{0} (X, f)}$, therefore the map is correctly defined
$$
\Lambda_{0} = \Psi_{0} \circ \Theta : \Ob \BB_{0}' (X, f) \rightarrow
\Ob \LL (\Sigma (X, f)) = (\Sigma (X, f), \leq) \;.
$$

From Lemma~\ref{lemma_3.5} and Proposition~\ref{prop_3.14} we obtain

\begin{corr} \label{corr_3.10}
The map $\Lambda_{0}$ preserves the order relation.

Pre-image $\Lambda_{0}^{-1} (\Phi (\EssP{X, f}))$ of the greatest element
of the set $(\Sigma (X, f), \leq)$ coincides with the class of all maximal
elements from $(\Ob \BB_{0}' (X, f), \preceq)$.
\end{corr}

From definition of the category $\LL (\Sigma (X, f))$, Lemma~\ref{lemma_3.5}
and Corollary~\ref{corr_3.10} we get

\begin{corr} \label{corr_3.11}
Map $\Lambda_{0} : \Ob \BB_{0}' (X, f) \rightarrow
\Ob \LL (\Sigma (X, f))$ is uniquely extended to the functor
$$
\Lambda : \BB_{0}' (X, f) \rightarrow \LL (\Sigma (X, f)) \;.
$$

For any two objects $B$, $B' \in \Ob \BB_{0}' (X, f)$ the map
$$
\Lambda_{B, \, B'} : H_{\BB_{0}' (X, f)} (B, B') \rightarrow
H_{\LL (\Sigma (X, f))} (\Lambda (B), \Lambda (B'))
$$
is injective.
\end{corr}

\begin{rem} \label{rem_3.14}
1) {}
It follows from Lemma~\ref{lemma_3.5}, that
the equality $\Phi (\EssP{A_{1}, g_{1}}) =
\Phi (\EssP{A_{2}, g_{2}})$ is the necessary condition
for two objects $(h_{1}, (A_{1}, g_{1}))$, $(h_{2}, (A_{2}, g_{2})) \in
\Ob \BB (X, f)$ to be isomorphic.

2) {}
From Corollary~\ref{corr_3.5} we obtain the statement: if
$\Phi (\EssP{A_{1}, g_{1}}) = \Phi (\EssP{A_{2}, g_{2}})$, then objects
$(h_{1}, (A_{1}, g_{1}))$, $(h_{2}, (A_{2}, g_{2})) \in \Ob \BB (X, f)$
are isomorphic if and only if at least one of the sets
$$
H_{\BB (X, f)} ((h_{1}, (A_{1}, g_{1})), (h_{2}, (A_{2}, g_{2}))) \;,
$$
$$
H_{\BB (X, f)} ((h_{2}, (A_{2}, g_{2})), (h_{1}, (A_{1}, g_{1})))
$$
is not empty.
\end{rem}

\subsection{Indecomposable dynamical systems.}\label{subsect_3.7}

We have natural desire to ``compare'' somehow categories
$\BB (X, f)$ and $\LL (\Sigma (X, f))$.

In what follows we shall see, that if dynamical system $(X, f)$ is indecomposable,
then the category $\LL (\Sigma (X, f))$ is isomorphic to a skeleton of the category $\BB (X, f)$
(and the isomorphism is set by the functor $\Lambda$).

In the case, when dynamical system $(X, f)$ is not indecomposable,
generally speaking it is not clear how to ``compare'' categories $\BB (X, f)$ and
$\LL (\Sigma (X, f))$, as shows following

\begin{lemma} \label{lemma_3.8}
Let $(h, (A, g)) \in \Ob \BB_{0}' (X, f)$, $N \in \Sigma (X, f)$ and
$\Phi (\EssP{A, g}) \lneqq N$.

Object $(h', (A', g ')) \in \Ob \BB_{0}' (X, f)$, which satisfies
to Lemma~\ref{lemma_3.7}, is defined uniquely if and only if
dynamical system $(X, f)$ is indecomposable.
\end{lemma}

\begin{proof}
{\bf 1)}{}
Assume, that dynamical system $(X, f)$ is indecomposable.

Let $(h_{1}, (A_{1}, g_{1}))$, $(h_{2}, (A_{2}, g_{2})) \in
\Ob \BB_{0}' (X, f)$, $\Phi (\EssP{A_{1}, g_{1}}) = \Phi (\EssP{A_{2}, g_{2}})
= N$ and $(h, (A, g)) \preceq (h_{i}, (A_{i}, g_{i}))$, $i = 1, 2$.

Since $\Phi (\EssP{A_{1}, g_{1}}) = \Phi (\EssP{A_{2}, g_{2}})$, then
dynamical systems $(A_{1}, g_{1})$ and $(A_{2}, g_{2})$ are topologically
conjugate. We fix isomorphism $\rho : (A_{1}, g_{1}) \rightarrow
(A_{2}, g_{2})$.

We conclude from Remark~\ref{rem_3.8} and Lemma~\ref{lemma_3.4} that
the following commutative diagram holds true
$$
\begin{CD}
(X, f)          @=                                  (X, f)                       @=            (X, f) \\
@V{h_{1}}VV                                    @V {\rho^{-1} \circ h_{2}}VV       @VV {h_{2}} V \\
(A_{1}, g_{1}) @>>{\fact (\rho^{-1} \circ h_{2})}>  (A_{1}, g_{1})  @>>{\rho}>      (A_{2}, g_{2})
\end{CD}
$$

In accord with Corollary~\ref{corr_3.5} the map $\rho \circ
\fact (\rho^{-1} \circ h_{2}) : (A_{1}, g_{1}) \rightarrow (A_{2}, g_{2})$
is isomorphism.

$\BB_{0}' (X, f)$ is the skeleton of $\BB (X, f)$, therefore
$(h_{1}, (A_{1}, g_{1})) = (h_{2}, (A_{2}, g_{2}))$.

{\bf 2)}{}
Let now dynamical system $(X, f)$ is not indecomposable.

Designate $\Phi (\EssP{A, g}) = M \in \Sigma (X, f)$. Since $M \lneqq N$
on condition of Lemma, then there exists prime $p \in \ss$, such that
$M_{p} \lneqq N_{p}$. This implies, that $M_{p} \neq \infty$.
Let $M_{p} = k$. Then $N_{p} \geq k+1$.

We fix regular sequence $\{ n_{i} \}_{i \in \nn}$, such
that $\Phi (\{ n_{i} \; | \; i \in \nn \}) = M$ (see.
Proposition~\ref{prop_2.2}). It follows from definition of function $\Phi$ that
\begin{itemize}
    \item[---] $n_{i} = p^{k_{i}} a_{i}$, $k_{i} \leq k$, $\gcd{a_{i}}{p}
        = 1$, $i \in \nn$;
    \item[---] there exists such $i_{0} \in \nn$, that $k_{i_{0}} = k$.
\end{itemize}
Sequence $\{ n_{i} \}_{i \in \nn}$ is regular, so
$k_{i} = k$ for all $i \geq i_{0}$. Without loss of generality
(see Corollary~\ref{corr_2.1}), we can consider that
$$
n_{i} = p^{k} a_{i} \,, \quad \gcd{a_{i}}{p} = 1 \,, \quad i \in \nn \;.
$$
Once again we take advantage of Corollary~\ref{corr_2.1} and we shall consider, that
$n_{1} = p^{k}$.

Fix regular sequence $\{ m_{j} \}_{j \in \nn}$, such
that $\Phi (\{ m_{j} \; | \; j \in \nn \}) = N $(see the beginning
of proof of Lemma~\ref{lemma_3.6}).
Sequence $\{ m_{j} \}_{j \in \nn}$ is regular and $N_{p} \geq
k+1$, therefore there exists $j_{0} \in \nn$, such that $m_{j}$ is divided by
$p^{k+1}$ for every $j \geq j_{0}$.

Again using Corollary~\ref{corr_2.1} we assume, that $m_{j}$
is divided by $p^{k+1}$ for all $j \in \nn$ and $m_{1} = p^{k+1}$.

We fix regular sequence $\{ V^{(n_{i})} \}_{i \in \nn}$ of
periodic partitions of dynamical system $(A, g)$. Consider
pre-images $W^{(n_{i})} = \{ W^{(n_{i})}_{s_{i}} =
h^{-1} (V^{(n_{i})}_{s_{i}}) \}_{s_{i} \in \zz_{n_{i}}}$, $i \in \nn$, of
periodic partitions of the sequence $\{ V^{(n_{i})} \}_{i \in \nn}$.

Repeating argument from the proof of Lemma~\ref{lemma_3.6} we conclude,
that $\{ W^{(n_{i})} \}_{i \in \nn}$ is the regular sequence of
periodic partitions of dynamical system $(X, f)$ and the partition $\hh$ of
space $X$, which is induced by this sequence, coincides with
the partition $\zer h$.

Since $N \leq \Phi (\EssP{X, f})$ on condition of Lemma, then there exists a
coherent regular sequence $\{ U^{(m_{j})} \}_{j \in \nn}$ of
periodic partitions of dynamical system $(X, f)$, which is compatible with the
sequence $\{ W^{(n_{i})} \}_{i \in \nn}$
(see Proposition~\ref{prop_1.8}).

We have assumed, that the dynamical system $(X, f)$ is not
indecomposable. Therefore, there exist partition $X = X_{1} \coprod X_{2}$ of
space $X$ on two proper disjoint invariant closed subsets $X_{1}$ and
$X_{2}$.

Families of sets $\{ P^{(m_{1})}_{s_{1}} = U^{(m_{1})}_{s_{1}}
\cap X_{1} \}_{s_{1} \in \zz_{m_{1}}}$ and $\{ Q^{(m_{1})}_{s_{1}} =
U^{(m_{1})}_{s_{1}} \cap X_{2} \}_{s_{1} \in \zz_{m_{1}}}$ are
periodic partitions of dynamical systems $(X_{1}, f |_{X_{1}})$ and
$(X_{2}, f |_{X_{2}})$, respectively (see Remark~\ref{rem_1.3} and
the proof of Proposition~\ref{prop_1.3}). Therefore the family of sets
$$
\widetilde{U}^{(m_{1})}_{s_{1}} = P^{(m_{1})}_{s_{1}} \cup
f^{n_{1}} (Q^{(m_{1})}_{s_{1}}) = P^{(m_{1})}_{s_{1}} \cup
f^{p^{k}} (Q^{(m_{1})}_{s_{1}}) \,, \quad s_{1} \in \zz_{m_{1}} \;,
$$
is the periodic partition of dynamical system $(X, f)$ of length
$m_{1} = p^{k+1}$ (see the proof of Proposition~\ref{prop_1.3}).

Periodic partitions $U^{(m_{1})}$ and $W^{(n_{1})}$ are compatible and
$n_{1}$ divides $m_{1}$, hence there exists $\tau \in \zz_{n_{1}}$,
such that $U^{(m_{1})}_{0} = P^{(m_{1})}_{0} \cup Q^{(m_{1})}_{0}
\subseteq W^{(n_{1})}_{\tau}$ (see Corollary~\ref{corr_1.6}).

Remark, that since $f : X \rightarrow X$ is homeomorphism, then
$$
f^{n_{1}} (Q^{(m_{1})}_{0}) = f^{n_{1}} (Q^{(m_{1})}_{0} \cap
W^{(n_{1})}_{\tau}) = f^{n_{1}} (Q^{(m_{1})}_{0}) \cap
f^{n_{1}} (W^{(n_{1})}_{\tau}) = f^{n_{1}} (Q^{(m_{1})}_{0}) \cap
W^{(n_{1})}_{\tau} \subseteq W^{(n_{1})}_{\tau} \;.
$$
Hence $\widetilde{U}^{(m_{1})}_{0} \subseteq W^{(n_{1})}_{\tau}$.
From Corollary~\ref{corr_1.5} we conclude, that periodic partitions
$\widetilde{U}^{(m_{1})}$ and $W^{(n_{1})}$ are compatible.

The inductive application of Proposition~\ref{prop_1.7} and
Remark~\ref{rem_1.6} gives us a coherent sequence
$\{ \widetilde{U}^{(m_{j})} \}_{j \in \nn}$ of periodic partitions of
dynamical system $(X, f)$, which is compatible with the sequence
$\{ W^{(n_{i})} \}_{i \in \nn}$.

Let $\Tfrak$ and $\widetilde{\Tfrak}$ are the partitions of space $X$,
induced by sequences $\{ U^{(m_{j})} \}_{m_{j} \in \nn}$ and
$\{ \widetilde{U}^{(m_{j})} \}_{m_{j} \in \nn}$, respectively.

Iterating argument from the proof of Lemma~\ref{lemma_3.6}, we shall find
$(h_{1}, (A_{1}', g_{1}'))$,
$(h_{2}, (A_{2}', g_{2}')) \in \Ob \BB (X, f)$, such that
$\Tfrak = \zer h_{1}$, $\widetilde{\Tfrak} = \zer h_{2}$ and
$$
h_{\BB (X, f)} ((h_{i}, (A_{i}', g_{i}')), (h, (A, g))) \neq
\emptyset \,, \quad i = 1, 2 \;.
$$

Find $ (\pi_{1}, (A_{1}, g_{1}))$, $(\pi_{2}, (A_{2}, g_{2})) \in
\Ob \BB_{0}' (X, f)$, which are isomorphic to objects
$(h_{1}, (A_{1}', g_{1}'))$ and
$(h_{2}, (A_{2}', g_{2}'))$ respectively.
Obviously,
$$
h_{\BB (X, f)} ((\pi_{i}, (A_{i}, g_{i})), (h, (A, g))) \neq
\emptyset \,, \quad i = 1, 2 \;.
$$

In order to verify the inequality $(\pi_{1}, (A_{1}, g_{1}))
\neq (\pi_{2}, (A_{2}, g_{2}))$ it suffices to show, that objects
$(h_{1}, (A_{1}', g_{1}'))$ and
$(h_{2}, (A_{2}', g_{2}'))$ are not isomorphic. For this purpose we shall
take an advantage of Remark~\ref{rem_3.14} (remind,
that $\Phi (\EssP{A_{1}', g_{1}'}) =
\Phi (\EssP{A_{2}', g_{2}'}) = N$ on the construction, hence
dynamical systems $(A_{1}', g_{1}')$ and
$(A_{2}', g_{2}')$ are topologically conjugate).

Let us check the equality
\begin{equation} \label{eq_17}
H_{\BB (X, f)} ((h_{1}, (A_{1}', g_{1}')),
( h_{2}, (A_{2}', g_{2}'))) = \emptyset \;.
\end{equation}

Assume, that this equality is not valid and there exists $\alpha :
(h_{1}, (A_{1}', g_{1}')) \rightarrow (h_{2}, (A_{2}', g_{2}'))$.
Therefore $\widetilde{\alpha} = \Theta (\alpha) :
(A_{1}', g_{1}') \rightarrow (A_{2}', g_{2}')$
is isomorphism (see Corollary~\ref{corr_3.5}). Hence,
$\zer (\widetilde{\alpha} \circ h_{1}) = \zer h_{1}$. Since $h_{2} =
\widetilde{\alpha} \circ h_{1}$ by definition, then $\widetilde{\Tfrak}
= \zer h_{2} = \zer (\widetilde{\alpha} \circ h_{1}) = \zer h_{1} =
\Tfrak$.

On the other hand, sequences $\{ U^{(m_{j})} \}_{j \in \nn}$ and
$\{ \widetilde{U}^{(m_{j})} \}_{j \in \nn}$ are not compatible.
Really, on the construction we have $\emptyset \neq U^{(m_{1})}_{0} \cap
\widetilde{U}^{(m_{1})}_{0} = P^{(m_{1})}_{0} \subseteq X_{1}$,
Hence
$$
\bigcup_{n \in \zz} f^{n} (U^{(m_{1})}_{0} \cap
\widetilde{U}^{(m_{1})}_{0}) \subseteq X_{1}
$$
and periodic partitions $U^{(m_{1})}$ and $\widetilde{U}^{(m_{1})}$
are not compatible. Therefore Corollary~\ref{corr_3.2} consequences that
$\Tfrak \neq \widetilde{\Tfrak}$.

The obtained contradiction proves the equality~(\ref{eq_17}).

The equality
$$
H_{\BB (X, f)} ((h_{2}, (A_{2}', g_{2}')),
(h_{1}, (A_{1}', g_{1}'))) = \emptyset
$$
is proved similarly.
\end{proof}

\begin{rem} \label{rem_3.15}
Obviously, the set $(\Sigma, \leq)$ has the least element
$E = (E_{p} = 0)_{p \in \ss} = \Phi_{0} (1)$.

Hence, for any dynamical system $(X, f)$ the element
$E \in \Ob \LL (\Sigma (X, f))$ is the right zero of category
$\LL (\Sigma (X, f))$.

In the category $\Acal_{0} (X, f)$ dynamical
system $(\{ pt \}, Id) $ with a phase space consisting of one point corresponds
to the element $E$ in the sense that $\Phi (\EssP{\{ pt \}, Id})$ = E.
\end{rem}

\begin{prop}
The category $\BB_{0}' (X, f)$ has the right zero $0_{R} \in \Ob
\BB_{0}' (X, f)$.
\end{prop}

\begin{proof}
Obviously, there exists exactly one projection $\pi_{0} : X \rightarrow
\{ pt \}$ and it satisfies to the relation $\pi_{0} \circ f = Id \circ
\pi_{0}$. Designate $0_{R} = (\pi_{0}, (\{ pt \}, Id))$.

Let $(h, (A, g)) \in \Ob \BB_{0}' (X, f) $. Obviously,
the projection $\pi : A \rightarrow \{ pt \}$ is uniquely defined and relations are fulfilled
$\pi \circ g = Id \circ \pi$ and $\pi \circ h = \pi_{0} :
(X, f) \rightarrow (\{ pt \}, Id)$.
\end{proof}

\begin{rem} \label{rem_3.16}
Since any two different objects of category $\BB_{0}' (X, f)$ are not
isomorphic (on definition of skeleton of a category), then the right zero is defined
uniquely.
\end{rem}

Now we can extract from Lemma~\ref{lemma_3.8} the following statements.

\begin{corr} \label{corr_3.12}
Let $\Phi (\EssP{X, f}) \neq E$, $N \in \Sigma (X, f)$ and $N \neq E$.

There exists $(h, (A, g)) \in \Ob \BB_{0}' (X, f)$, such that
$\Phi (\EssP{A, g}) = N$ and the following statements are equivalent:
\begin{itemize}
    \item[(i)] object $(h, (A, g)) \in \Ob \BB_{0}' (X, f)$, such that
        $\Phi (\EssP{A, g}) = N$, is defined uniquely;
    \item[(ii)] dynamical system $(X, f)$ is indecomposable.
\end{itemize}
\end{corr}

\begin{proof}
We apply Lemmas~\ref{lemma_3.7} and~\ref{lemma_3.8} to object $0_{R} \in
\Ob \BB_{0}' (X, f)$ and number $N \in \Sigma (X, f)$.
\end{proof}

\begin{rem} \label{rem_3.17}
In other words Corollary~\ref{corr_3.12} can be formulated as follows:
\begin{itemize}
    \item[---] map $\Lambda_{0} : \Ob \BB_{0}' (X, f) \rightarrow
        \Ob \LL (\Sigma (X, f)) $ is surjective;
    \item[---] if $ \Phi (\EssP{X, f}) \neq E$, then injectivity of
        $\Lambda_{0}$ is equivalent to that the dynamical
        system $(X, f)$ is indecomposable.
\end{itemize}
\end{rem}

\begin{corr} \label{corr_3.13}
If dynamical system $(X, f)$ is indecomposable, then for any two objects
$(h_{1}, (A_{1}, g_{1}))$, $(h_{2}, (A_{2}, g_{2})) \in
\Ob \BB_{0}' (X, f)$ the inequality
$$
H_{\BB (X, f)} ((h_{2}, (A_{2}, g_{2})),
(h_{1}, (A_{1}, g_{1}))) \neq \emptyset
$$
is fulfilled if and only if
$\Phi (\EssP{A_{1}, g_{1}}) \leq \Phi (\EssP{A_{2}, g_{2}})$.
\end{corr}

\begin{proof}
Let $(h_{1}, (A_{1}, g_{1}))$, $ (h_{2}, (A_{2}, g_{2})) \in
\Ob \BB_{0}' (X, f)$.

If $H_{\BB (X, f)} ((h_{2}, (A_{2}, g_{2})),
(h_{1}, (A_{1}, g_{1}))) \neq \emptyset$, then Lemma~\ref{lemma_3.5}
consequences that $\Phi (\EssP{A_{1}, g_{1}}) \leq \Phi (\EssP{A_{2}, g_{2}})$.

Let now $\Phi (\EssP{A_{1}, g_{1}}) \leq \Phi (\EssP{A_{2}, g_{2}})$.
We take advantage of Lemma~\ref{lemma_3.7} and find $(h_{2}', (A_{2}', g_{2}'))
\in \Ob \BB_{0}' (X, f)$, such that $\Phi (\EssP{A_{2}', g_{2}'}) =
\Phi (\EssP{A_{2}, g_{2}})$ and
$$
H_{\BB (X, f)} ((h_{2}', (A_{2}', g_{2}')),
(h_{1}, (A_{1}, g_{1}))) \neq \emptyset \;.
$$
From Corollary~\ref{corr_3.12} we conclude, that $(h_{2}', (A_{2}', g_{2}'))
= (h_{2}, (A_{2}, g_{2}))$.
\end{proof}

\begin{lemma} \label{lemma_3.9}
Let $\Phi (\EssP{X, f}) \neq E$. Then the following statements
are equivalent:
\begin{itemize}
    \item[(i)] category $\BB_{0}' (X, f)$ has the left zero $0_{L}
        \in \Ob \BB_{0}' (X, f)$;
    \item[(ii)] dynamical system $(X, f)$ is indecomposable.
\end{itemize}
\end{lemma}

\begin{proof}
{\bf 1)} {}
Let dynamical system $(X, f)$ is indecomposable.

From Corollary~\ref{corr_3.12} it follows, that there exists a unique object
$(h, (A, g)) \in \Ob \BB_{0}' (X, f)$, such that $\Phi (\EssP{A, g}) =
\Phi (\EssP{X, f})$.

Theorem~\ref{theorem_A} guarantees that for any $(h', (A', g'))
\in \Ob \BB_{0}' (X, f)$ the inequality $\Phi (\EssP{A', g'}) \leq
\Phi (\EssP{A, g})$ is fulfilled. Now Corollary~\ref{corr_3.13} and
Lemma~\ref{lemma_3.5} show, that for every $(h', (A', g')) \in
\Ob \BB_{0}' (X, f)$ the set
$$
H_{\BB (X, f)} ((h, (A, g)), (h', (A', g')))
$$
contains exactly one element, hence $(h, (A, g))$ is the left zero
of $\BB_{0}' (X, f)$.

{\bf 2)} {}
Let dynamical system $(X, f)$ is not indecomposable.

Assume, that there exists a left zero $(h, (A, g))$ of category
$\BB_{0}' (X, f)$. From Theorem~\ref{theorem_A} and Lemma~\ref{lemma_3.5}
we conclude, that $\Phi (\EssP{A, g}) = \Phi (\EssP{X, f})$.

Corollary~\ref{corr_3.12} implies, that there exists $(h', (A', g'))
\in \Ob \BB_{0}' (X, f)$, such that $\Phi (\EssP{A', g'}) =
\Phi (\EssP{X, f})$ and $(h', (A', g')) \neq (h, (A, g))$.

Since $\BB_{0}' (X, f)$ is the skeleton of $\BB (X, f)$
(see Proposition~\ref{prop_3.13}), then objects $(h', (A', g'))$ and
$(h, (A, g))$ are not isomorphic in $\BB (X, f)$. From
Remark~\ref{rem_3.14} we conclude, that
$$
H_{\BB (X, f)} ((h, (A, g)), (h', (A', g'))) = \emptyset
$$
and the object $(h, (A, g))$ can not be left zero of
$\BB_{0}' (X, f)$.

The obtained contradiction finishes the proof.
\end{proof}

\begin{theorem} \label{theorem_3.5}
Let $\Phi (\EssP{X, f}) \neq E$, $\BB_{0} (X, f)$ is a skeleton of the category
$\BB (X, f)$.

The following statements are equivalent:
\begin{itemize}
    \item[(i)] dynamical system $(X, f)$ is indecomposable;
    \item[(ii)] categories $\BB_{0} (X, f)$ and
        $\LL (\Sigma (X, f))$ are isomorphic.
\end{itemize}
\end{theorem}

\begin{proof}
On definition $\LL (\Sigma (X, f))$ has left zero
$\Phi (\EssP{X, f}) \in \Ob \LL (\Sigma (X, f))$. Hence,
the existence of left zero in the category $\BB_{0} (X, f)$ is the
necessary condition for categories $\BB_{0} (X, f)$ and
$\LL (\Sigma (X, f))$ to be isomorphic.

Subcategory $\BB_{0}' (X, f)$ is the skeleton of $\BB (X, f)$
(see Proposition~\ref{prop_3.13}), therefore it is isomorphic to the skeleton
$\BB_{0} (X, f)$. If dynamical system $(X, f)$ is not
indecomposable, then we conclude from Lemma~\ref{lemma_3.9}, that
there is no left zero in the category $\BB_{0} (X, f)$  and it can not
be isomorphic to $\LL (\Sigma (X, f))$.

Assume now that $(X, f)$ is an indecomposable dynamical system.
We shall prove, that the functor $ \Lambda : \BB_{0}' (X, f) \rightarrow
\LL (\Sigma (X, f))$ (see Corollary~\ref{corr_3.11})
sets the isomorphism of categories $\BB_{0}' (X, f)$ and $\LL (\Sigma (X, f))$.

From Remark~\ref{rem_3.17} the map $\Lambda_{0} :
\Ob \BB_{0}' (X, f) \rightarrow \Ob \LL (\Sigma (X, f))$ appears
to be bijective.

Let $M$, $N \in \Ob \LL (\Sigma (X, f))$. Corollary~\ref{corr_3.13}
consequences that inequalities
$$
H_{\LL (\Sigma (X, f))} (M, N) \neq \emptyset \quad \mbox{ and } \quad
H_{\BB_{0}' (X, f)} (\Lambda_{0}^{-1} (M), \Lambda_{0}^{-1} (N)) \neq
\emptyset
$$
are equivalent. Since each of the sets $H_{\LL (\Sigma (X, f))} (M, N)$
and $H_{\BB_{0}' (X, f)} (\Lambda_{0}^{-1} (M), \Lambda_{0}^{-1} (N))$ contain
no more than one element (see definition of category
$\LL (\Sigma (X, f))$ and Lemma~\ref{lemma_3.5}), then
$$
\Lambda_{\Lambda_{0}^{-1} (M), \, \Lambda_{0}^{-1} (N)} :
H_{\BB_{0}' (X, f)} (\Lambda_{0}^{-1} (M), \Lambda_{0}^{-1} (N))
\rightarrow
H_{\LL (\Sigma (X, f))} (M, N)
$$
is bijective map for any pair $M$, $N \in \Ob \LL (\Sigma (X, f))$.

Thus, $\Lambda : \BB_{0}' (X, f) \rightarrow \LL (\Sigma (X, f))$
sets the isomorphism of categories $\BB_{0}' (X, f)$ and $\LL (\Sigma (X, f))$.
Therefore, categories $\BB_{0} (X, f)$ and $\LL (\Sigma (X, f))$ are also
isomorphic.
\end{proof}

\subsection{Expansions of odometers and almost periodic
points} \label{subsect_3.8}

\begin{prop} \label{prop_3.16}
Let $(A, g) \in \Ob \Acal (X, f)$.

Suppose there exists projection $\pi : (X, f) \rightarrow (A, g)$,
such that for a certain $x \in X$ the equality $\pi^{-1} (\pi (x))
= \{ x \} $ is fulfilled.

Then the following statements hold true:
\begin{itemize}
    \item[(i)] $x \in X$ is the almost periodic point of dynamical
        system $(X, f)$;
    \item[(ii)] for any $y \in X$ the inclusion is valid
        $\Cl{\Orb_{f} (x)} \subseteq \alpha (y) \cap \omega (y)$ (hence,
        $\Cl{\Orb_{f} (x)}$ is the unique minimal set of dynamical
        system $(X, f)$, in particular $(X, f)$ is indecomposable);
    \item[(iii)] $\Phi (\EssP{A, g}) = \Phi (\EssP{X, f})$;
    \item[(iv)] for any projection $\pi' : (X, f) \rightarrow (A', g')$,
        $(A', g') \in \Ob \Acal$, the following conditions are equivalent:
        \begin{itemize}
            \item[{\bfseries a)}] $(\pi')^{-1} (\pi' (x)) = \{ x \}$,
            \item[{\bfseries b)}] $\Phi (\EssP{A', g'}) = \Phi (\EssP{X, f})$;
        \end{itemize}
    \item[(v)] if dynamical system $(X, f)$ is minimal, then for
        every almost periodic point $y \in X$ of dynamical system
        $(X, f)$ the equality $\pi^{-1} (\pi (y)) = \{ y \}$ is fulfilled.
\end{itemize}
\end{prop}

Before we proceed to prove Proposition~\ref{prop_3.16} and extract
corollaries from it, we shall prove three lemmas.

\begin{lemma} \label{lemma_3.10}
Let $(A_{1}, g_{1})$, $(A_{2}, g_{2}) \in \Ob \Acal$,
$h : (A_{1}, g_{1}) \rightarrow (A_{2}, g_{2})$ is a morphism.

If there exists $x \in A_{1}$, such that $h^{-1} (h (x)) = \{ x \}$,
then $h \in \Iso A$.
\end{lemma}

\begin{proof}
Since $(A_{2}, g_{2})$ is minimal dynamical system
(see Remark~\ref{rem_3.2}), then $h$ is projection
(see Remark~\ref{rem_3.8}).

Let $y \in A_{1}$. Designate by $H (y)$ the element of partition
$\zer h$ of space $A_{1}$, which contains $y$.

We know from Corollary~\ref{corr_3.6} that there exists unique
isomorphism $h_{y, \, x} : (A_{1}, g_{1}) \rightarrow (A_{1}, g_{1})$, such
that $h_{y, \, x} (y) = x$. Consider a projection $\widetilde{h} =
h \circ h_{y, \, x} : (A_{1}, g_{1}) \rightarrow (A_{2}, g_{2})$. Designate
by $\widetilde{H} (y)$ the element of partition $\zer \widetilde{h}$
of space $A_{1}$, which contains $y$.

Let $z = h (x) \in A_{2}$. It is clear, that
$$
\widetilde{h}^{-1} (z) = h_{y, \, x}^{-1} (h^{-1} (z)) = h_{y, \, x}^{-1} (x) =
\{ y \} = \widetilde{H} (y) \;.
$$
Dynamical system $(A_{2}, g_{2})$ is minimal, hence it
is indecomposable. We conclude from Lemma~\ref{lemma_3.3} that $\zer h =
\zer \widetilde{h}$. Hence, $h^{-1} (h (y)) = H (y) = \widetilde{H} (y)
= \{ y \}$.
\end{proof}

\begin{lemma} \label{lemma_3.11}
Let dynamical system $(X, f)$ is minimal and $x \in X$ is almost
periodic point of this dynamical system.

Then there exist dynamical system $(A, g) \in \Acal (X, f)$ and projection
$\pi : (X, f) \rightarrow (A, g)$, such that $\Phi (\EssP{A, g}) =
\Phi (\EssP{X, f})$ and $\pi^{-1} (\pi (x)) = \{ x \}$.
\end{lemma}

\begin{proof}
We fix regular sequence $\{ n_{i} \}_{i \in \nn}$, such that
$\Phi (\{ n_{i} \; | \; i \in \nn \}) = \Phi (\EssP{X, f})$, and build regular
sequence $\{ W^{(n_{i})} \}_{i \in \nn}$ of periodic partitions of
dynamical system $(X, f)$. Without loss of generality we can
suppose that $x \in W^{(n_{i})}_{0}$, $i \in \nn$
(see Remark~\ref{rem_1.6}).

Let $\hh$ is the partition of space $X$, induced
by the sequence $\{ W^{(n_{i})} \}_{i \in \nn}$ and let
$$
F : (X, f) \rightarrow (X/\hh, \Bar {f}) = (A, g)
$$
is the projection to dynamical system $(A, g) \in \Acal (X, f)$
(see relation~(\ref{eq_12}), Remark~\ref{rem_3.4} and
Proposition~\ref{prop_3.2}). Then $\Phi (\EssP{A, g}) =
\Phi (\EssP{X, f})$.

Since $\zer F = \hh$, then in order to complete the proof it is enough
to verify the equality
$$
\{ x \} = \bigcap_{i \in \nn} W^{(n_{i})}_{0} \;.
$$

Assume, that this equality is invalid and there exists $y \neq x$, such
that
$$
y \in \bigcap_{i \in \nn} W^{(n_{i})}_{0} \;.
$$
The space $X$ is Hausdorff, therefore there exists a closed neighborhood
$U \subseteq X \setminus \{ y \}$ of $x$. Since $x$ is the almost
periodic point, then we conclude from Lemma~\ref{lemma_1.5} that there exist
$m \in \EssP{X, f}$ and periodic partition $\widetilde{W}^{(m)}$ of
dynamical system $(X, f)$, such that $x \in
\widetilde{W}^{(m)}_{0} \subseteq U$.

Consider the stationary regular sequence
$\{ m_{j} = m \}_{j \in \nn}$. From Lemma~\ref{lemma_2.1} it follows, that
$\Phi_{0} (m) = \Phi (\{ m_{j} \; | \; j \in \nn \}) \leq
\Phi (\EssP{X, f}) = \Phi (\{ n_{i} \; | \; i \in \nn \})$. Therefore
Proposition~\ref{prop_2.3} consequences that there exists $k \in \nn$, such that
$m_{1} = m$ divides $n_{k}$.

The dynamical system $(X, f)$ is minimal, hence it is indecomposable.
Corollary~\ref{corr_1.1} implies, that periodic partitions
$\widetilde{W}^{(m)}$ and $W^{(n_{k})}$ are compatible. Therefore
(see Corollary~\ref{corr_1.6}) the inclusions hold true $x \in
W^{(n_{k})}_{0} \subseteq \widetilde{W}^{(m)}_{0} \subseteq U
\subseteq X \setminus \{ y \}$ and
$$
y \notin \bigcap_{i \in \nn} W^{(n_{i})}_{0} \;.
$$

The obtained contradiction completes the proof.
\end{proof}

\begin{lemma} \label{lemma_3.12}
Let $x \in Y_{1}$ is an almost periodic point of a dynamical system
$(Y_{1}, h_{1})$. Let $\pi : (Y_{1}, h_{1}) \rightarrow (Y_{2}, h_{2})$
is a projection.

Then $y = \pi (x) \in Y_{2}$ is an almost periodic point of the dynamical
system $(Y_{2}, h_{2})$.
\end{lemma}

\begin{proof}
Let $U \subseteq Y_{2}$ is an open neighborhood of $y$. Since
$\pi : Y_{1} \rightarrow Y_{2}$ is continuous map, then
$V = \pi^{-1}(U)$ is the open neighborhood of the almost periodic point $x$.
Hence, there exists $n(V) \in \nn$, such that
$$
\bigcup_{k \in \zz} h_{1}^{k n (V)} (x) \subseteq V \;.
$$
Then
$$
\bigcup_{k \in \zz} h_{2}^{k n (V)} (y) =
\bigcup_{k \in \zz} h_{2}^{k n (V)} \circ \pi (x) =
\bigcup_{k \in \zz} \pi \circ h_{1}^{k n (V)} (x) \subseteq
\pi (V) = U \;.
$$
From arbitrariness in the choice of a neighborhood $U$ it follows, that
$y = \pi (x)$ is almost periodic point of dynamical system $(Y_{2}, h_{2})$.
\end{proof}

\begin{proof}[of Proposition~\ref{prop_3.16}]
{\bfseries (i)} {}
We fix regular sequence $\{ n_{i} \}_{i \in \nn}$, such
that $\Phi (\{ n_{i} \; | \; i \in \nn \}) = \Phi (\EssP{A, g})$. Construct
regular sequence $\{ V^{(n_{i})} \}_{i \in \nn}$ of
periodic partitions of dynamical system $(A, g)$. Without loss of generality
we can assume, that $z = \pi (x) \in V^{(n_{i})}_{0}$,
$i \in \nn$ (see Remark~\ref{rem_1.6}).

According to Proposition~\ref{prop_3.7} the family of sets
$\{ V^{(n_{i})}_{s_{i}} \; | \; s_{i} \in \zz_{n_{i}} \, \; i \in \nn \}$
is basis of the topology of space $A$. Therefore,
$$
\{ z \} = \bigcap_{i \in \nn} V^{(n_{i})}_{0} \;.
$$
In addition since the sequence $\{ V^{(n_{i})} \}_{i \in \nn}$ is
regular, then $V^{(n_{i+1})}_{0} \subseteq V^{(n_{i})}_{0}$, $i \in \nn$
(see Corollary~\ref{corr_1.6}).

Consider pre-images $\{ W^{(n_{i})}_{s_{i}} =
\pi^{-1} (V^{(n_{i})}_{s_{i}}) \; | \; s_{i} \in \zz_{n_{i}} \}$, $i \in
\nn$, of periodic partitions $\{ V^{(n_{i})}_{s_{i}} \; | \;
s_{i} \in \zz_{n_{i}} \}$, $i \in \nn$.

It is clear, that $W^{(n_{i+1})}_{0} \subseteq W^{(n_{i})}_{0}$,
$i \in \nn$, and
$$
\{ x \} = \pi^{-1} (z) = \bigcap_{i \in \nn} W^{(n_{i})}_{0} \;.
$$
We consequence from Proposition~\ref{prop_3.6} and Corollary~\ref{corr_1.5},
that $\{ W^{(n_{i})} \}_{i \in \nn}$ is the regular sequence of
periodic partitions of dynamical system $(X, f)$.

Let $U \subseteq X$ is an open neighborhood of $x$. Since
$X$ is compact, then all sets $W^{(n_{i})}_{0}$, $i \in \nn$, are compact.
Applying Lemma~\ref{lemma_3.2}, we find $k \in \nn$, such that $x \in
W^{(n_{k})}_{0} \subseteq U$.

On definition of periodic partition we have
$$
\bigcup_{m \in \zz} f^{m n_{k}} (x) \subseteq
\bigcup_{m \in \zz} f^{m n_{k}} (W^{(n_{k})}_{0}) =
W^{(n_{k})}_{0} \subseteq U \;.
$$

By virtue of arbitrariness in the choice of a neighborhood $U$
the point $x$ is almost periodic.

\medskip

{\bfseries (ii)} {}
Let $y \in X$, $t = \pi (y) \in A$. Dynamical system $(A, g)$
is minimal (see Remark~\ref{rem_3.2}), consequently $\alpha (t) = \omega (t) =
\Cl{\Orb_{g} (t)} = A$. Hence there exists monotonic unlimited
sequence of numbers $\{ n_{i} \in \zz \}_{i \in \nn}$,
such that $z = \pi (x) = \lim_{i \rightarrow \infty} g^{n_{i}} (t)$.

Consider the sequence $\{ f^{n_{i}} (y) \in X \}_{i \in \nn}$.
Space $X$ is compact, therefore this sequence has at least
one limit point $x' \in X $. Passing to a
subsequence we can assume, that $x' =
\lim_{i \rightarrow \infty} f^{n_{i}} (y)$. Hence,
$x' \in \omega (y)$.

On the other hand, $\pi \circ f^{n_{i}} (y) = g^{n_{i}} \circ \pi (y) =
g^{n_{i}} (t)$, therefore
$$
\pi (x') = \lim_{i \rightarrow \infty} \pi \circ f^{n_{i}} (x) =
\lim_{i \rightarrow \infty} g^{n_{i}} (t) = \pi (x)
$$
and $x = x' \in \omega (y)$. Since $\omega (y)$ is closed invariant
set of dynamical system $(X, f)$, then $\Cl{\Orb_{f} (x)} \subseteq
\omega (y)$.

The relation $\Cl{\Orb_{f} (x)} \subseteq \alpha (y)$ is proved
similarly.

\medskip

{\bfseries (iii)} {}
We consider $(\pi, (A, g)) \in \Ob \BB (X, f)$. It follows from Theorem~\ref{theorem_A},
that $\Phi (\EssP{A, g}) \leq \Phi (\EssP{X, f})$. Apply
Lemma~\ref{lemma_3.6} and find $(\pi', (A', g')) \in \Ob \BB (X, f)$, such
that $\Phi (\EssP{A', g'}) = \Phi (\EssP{X, f})$ and exists
$$
h \in H_{\BB (X, f)} ((\pi', (A', g')), (\pi, (A, g))) \;.
$$
In other words, there exist $(A', g') \in \Acal (X, f)$ and
$h : (A', g') \rightarrow (A, g)$, such that $\Phi (\EssP{A', g'}) =
\Phi (\EssP{X, f})$ and $\pi = h \circ \pi'$.

We conclude from Remark~\ref{rem_3.8}, that map $\pi' : X \rightarrow
A'$ is surjective. Hence for any subset $B \subseteq A'$
the equality $\pi' ((\pi')^{-1} (B)) = B$ is fulfilled.

Let $z' = \pi' (x) \in A'$. Then $h (z') = h \circ \pi' (x) = z$. The
equalities hold true
$(\pi')^{-1} (h^{-1} (z)) = \pi^{-1} (z) = \{ x \}$, hence
$$
h^{-1} (h (z')) = h^{-1} (z) = \pi' ((\pi')^{-1} (h^{-1} (z))) = \{ \pi' (x) \}
= \{ z' \} \;.
$$

Applying Lemma~\ref{lemma_3.10} we conclude, that $h : (A', g') \rightarrow
(A, g)$ is isomorphism. Hence, $\Phi (\EssP{A, g}) = \Phi (\EssP{A', g'})
= \Phi (\EssP{X, f})$.

\medskip

{\bfseries (iv)} {}
We have already proved in (iii) the implication {\bfseries a)} $\Rightarrow$
{\bfseries b)}.

Let $(A', g') \in \Acal (X, f)$, $\Phi (\EssP{A', g'})
= \Phi (\EssP{X, f})$ and $\pi' : (X, f) \rightarrow (A', g')$ is a projection.

Dynamical system $(X, f)$ is indecomposable (see. (ii)), therefore from
Lemma~\ref{lemma_3.8} it follows (we shall remind, that $\BB_{0}' (X, f)$
is the skeleton of $\BB (X, f)$), that objects $(\pi', (A', g'))$ and
$(\pi, (A, g))$ of the category $\BB (X, f)$ are isomorphic. Hence
there exists isomorphism $\rho : (A, g) \rightarrow (A', g')$, such that
$\pi' = \rho \circ \pi$. In addition it is obvious that $\pi = \rho^{-1}
\circ \pi'$. Applying Lemma~\ref{lemma_2} we conclude, that partitions
$\zer \pi$ and $\zer \pi'$ coincide.

The sets $\pi^{-1} (\pi (x))$ and $(\pi')^{-1} (\pi' (x))$ are elements of the
partitions $\zer \pi$ and $\zer \pi'$ respectively, and contain the point $x$.
Hence, $(\pi')^{-1} (\pi' (x)) = \pi^{-1} (\pi (x)) = \{ x \}$.

\medskip

{\bfseries (v)} {}
Let dynamical system $(X, f)$ is minimal and $y \in X$ is almost
periodic point of this dynamical system.

Applying Lemma~\ref{lemma_3.11}, we find projection $\pi' : (X, f)
\rightarrow (A', g')$, $(A', g') \in \Acal (X, f)$, such that
$(\pi')^{-1} (\pi') (y) = \{ y \}$. Then changing roles of projections
$\pi$ and $\pi'$, we shall receive from (iv) that $(\pi)^{-1} (\pi) (y)
= \{ y \}$.
\end{proof}

\begin{defn} \label{defn_3.7}
Dynamical system $(Y_{1}, h_{1})$ is called {\em almost
one-to-one expansion} of a dynamical system $(Y_{2}, h_{2})$,
if there exist projection $\pi : (Y_{1}, h_{1}) \rightarrow
(Y_{2}, h_{2})$ and a dense subset $Q \subseteq Y_{1}$,
such that $\pi^{-1} (\pi (y)) = \{ y \}$ for any $y \in Q$.
\end{defn}

\begin{corr}[see~\cite{Block}]\label{corr_3.14}
Dynamical system $(X, f)$ is almost one-to-one
expansion of an odometer if and only if $X = \Cl{\Orb_{f} (x)}$
for an almost periodic point $x \in X$.
\end{corr}

\begin{proof}
Let $\pi : (X, f) \rightarrow (Y, h)$ is a projection and $\pi^{-1} (\pi (x))
= \{ x \}$ for a certain $x \in X$.

Remark, that since $f : X \rightarrow X$ and $h : Y \rightarrow Y$
are the homeomorphisms, then for every $n \in \zz$ equalities hold true
$$
\{ x \} = \pi^{-1} (\pi (x))
= \pi^{-1} (h^{-n} \circ \pi \circ f^{n} (x))
= (h^{n} \circ \pi)^{-1} (\pi \circ f^{n} (x)) =
$$
$$
= (\pi \circ f^{n})^{-1} (\pi \circ f^{n} (x))
= f^{-n} (\pi^{-1} (\pi (f^{n} (x))))
\;,
$$
Hence
$$
\{ f^{n} (x) \} = f^{n} \circ f^{-n} (\pi^{-1} (\pi (f^{n} (x)))) =
\pi^{-1} (\pi (f^{n} (x))) \;.
$$

In other words, for every $y \in \Orb_{f} (x)$
equality $\pi^{-1} (\pi (y)) = \{ y \}$ is fulfilled.

If the dynamical system $(X, f)$ is minimal, then $X =
\Cl{\Orb_{f} (x)}$ and $(X, f)$ is almost one-to-one expansion
of the dynamical system $(Y, h)$.

{\bf 1.} {}
Let $X = \Cl{\Orb_{f} (x)}$ for some almost-periodic
point $x \in X$. Then dynamical system $(X, f)$ is minimal by
Birkgoff Theorem and applying Lemma~\ref{lemma_3.11} and argument
mentioned above, we conclude that $(X, f)$ is almost one-to-one
expansion of an odometer.

{\bf 2.} {}
Let dynamical system $(X, f)$ is almost
one-to-one expansion of an odometer $(A, g)$.

We fix a projection $\pi : (X, f) \rightarrow (A, g)$, such that
for each point $y$ from some everywhere dense set $Q
\subseteq X$ the equality $\pi^{-1} (\pi (y)) = \{ y \}$ is fulfilled.

From Proposition~\ref{prop_3.16}, item (i), we conclude, that any
point $y \in Q$ is an almost-periodic point of $(X, f)$. Hence, for every
$y \in Q$ the set $\Cl{\Orb_{f} (y)}$ is minimal set of dynamical
system $(X, f)$.

From here we consequence that either $\Cl{\Orb_{f} (y_{1})} =
\Cl{\Orb_{f} (y_{2})}$ or $\Cl{\Orb_{f} (y_{1})} \cap
\Cl{\Orb_{f} (y_{2})} = \emptyset$ for arbitrary $y_{1}$,
$y_{2} \in Q$.

We fix $x \in Q$. For any $y \in Q$ we have inclusions (see.
Proposition~\ref{prop_3.16}, item (ii)) $\Cl{\Orb_{f} (x)}
\subseteq \alpha (y) \cap \omega (y) \subseteq \Cl{\Orb_{f} (y)}$.
Hence, $\Cl{\Orb_{f} (x)} = \Cl{\Orb_{f} (y)}$, in particular
$y \in \Cl{\Orb_{f} (x)}$. Since $Q$ is dense in $X$, then $X = \Cl{Q}
= \Cl{\Orb_{f} (x)}$.
\end{proof}

\begin{corr} \label{corr_3.15}
Let $(Y_{1}, h_{1})$ is almost one-to-one expansion of an odometer,
$\pi : (Y_{1}, h_{1}) \rightarrow (Y_{2}, h_{2})$ is a projection.

Then $(Y_{2}, h_{2})$ is almost one-to-one expansion of an odometer.
\end{corr}

\begin{proof}
Since map $\pi : Y_{1} \rightarrow Y_{2}$ is surjective and
continuous, then for any dense subset $Q$
of space $Y_{1}$ its image $\pi (Q)$ is dense in
$Y_{2}$.

We conclude from Corollary~\ref{corr_3.14} that there exists almost
periodic point $x \in Y_{1}$, such that $Y_{1} =
\Cl{\Orb_{h_{1}} (x)}$. From Lemma~\ref{lemma_3.12} it follows, that
$\pi (x)$ is an almost-periodic point of the dynamical system $(Y_{2},
h_{2})$, and $\Cl{\Orb_{h_{2}} (\pi (x))} = Y_{2}$ (see above).

Again applying Corollary~\ref{corr_3.14}, we come to a conclusion, that
$(Y_{2}, h_{2})$ is almost one-to-one expansion of
odometer.
\end{proof}

\begin{corr} \label{corr_3.16}
Let $(X, f)$ is almost one-to-one expansion of certain
odometer, $(A, g) \in \Acal (X, f)$ and $\pi : (X, f) \rightarrow
( A, g)$ is a projection.

Let $Q = \{ y \in X \; | \; \pi^{-1} (\pi (y)) = \{ y \} \} $.

If $\Phi (\EssP{A, g}) = \Phi (\EssP{X, f})$, then $Q$ coincides with
the set of all almost periodic points of the dynamical system $(X, f)$.

If $\Phi (\EssP{A, g}) \neq \Phi (\EssP{X, f})$, then $Q = \emptyset$.
\end{corr}

\begin{proof}
From Corollary~\ref{corr_3.14} and Birkgoff theorem we conclude, that
the dynamical system $(X, f)$ is minimal.

The first part of the current statement follows from Lemma~\ref{lemma_3.11}
and Proposition~\ref{prop_3.16}, items~(i), (iv) and~(v).

The second part follows from Proposition~\ref{prop_3.16}, item~(iii).
\end{proof}

\begin{corr}[see~\cite{Block}]\label{corr_3.17}
Let $(X, f)$ is a transitive dynamical system.

The dynamical system $(X, f)$ is topologically conjugate with an odometer
if and only if each point of space $X$ is an almost periodic point of
dynamical system $(X, f)$.
\end{corr}

\begin{proof}
{\bf 1.} {}
Let each point of space $X$ is an almost-periodic point of $(X, f)$.
Since dynamical system $(X, f)$ is transitive, then on definition
there exists $x \in X$, such that $X = \Cl{\Orb_{f} (x)}$. Applying
Corollary~\ref{corr_3.14} we conclude, that $(X, f)$ is almost
one-to-one expansion of an odometer (in particular, dynamical
system $(X, f)$ is minimal under Birkgoff theorem). We take advantage
of Lemma~\ref{lemma_3.11} and find $(A, g) \in \Acal (X, f)$ and
projection $\pi : (X, f) \rightarrow (A, g)$, such that
$\pi^{-1} (\pi (x)) = \{ x \}$.

From Corollary~\ref{corr_3.16} we consequence, that map
$\pi : X \rightarrow A$ is bijective. Space $X$
is compact, hence $\pi$ is the homeomorphism and $\pi : (X, f)
\rightarrow (A, g)$ is isomorphism in the category $\kk_{0}$.

{\bf 2.} {}
Let $(X, f) \in \Acal$. Consider the unit morphism $Id :
(X, f) \rightarrow (X, f)$. Since map $Id : X
\rightarrow X$ is bijective, then we conclude from
Proposition~\ref{prop_3.16}, item~(i), that every
point of space $X$ is an almost-periodic point of
dynamical system $(X, f)$.
\end{proof}

\begin{corr}[see~\cite{Block}]\label{corr_3.18}
Let dynamical system $(A, g)$ is topologically conjugate with an odometer,
$\pi : (A, g) \rightarrow (X, f) $ is a projection.

Then dynamical system $(X, f)$ is topologically conjugate with an odometer.
\end{corr}

\begin{proof}
The dynamical system $(A, g)$ is minimal (see.
Remark~\ref{rem_3.2}), therefore it is transitive. Hence
the dynamical system $(X, f)$ is transitive too (see proof of
Corollary~\ref{corr_3.15}).

From Corollary~\ref{corr_3.17} and Lemma~\ref{lemma_3.12} we conclude,
that each point of the space $X$ is an almost-periodic point of
dynamical system $(X, f)$.

In order to complete the proof it suffices to apply
Corollary~\ref{corr_3.17} once again.
\end{proof}

\end{document}